\documentclass{amsart}
\usepackage{amssymb}
\usepackage{graphics}
\usepackage[hyperref]{degt}

\theoremstyle{definition}
\addtheorem{Algorithm}{algorithm}

\makeatletter
\def\HyPsd@CatcodeWarning#1{}
\makeatother

\def\dts{\mathbin{..}}
\def\ds{\QOPNAME{ds}}
\def\gd{\QOPNAME{gd}}

\let\graph\Gamma
\def\pencil{\Cal P}
\let\fiber\Sigma

\def\mm{M}
\def\bm{\bar M}
\def\Ms{\bar M_{\sec}}
\def\Mg{\bar M_\graph}
\def\Mm{\bar M}
\def\mr{M_\R}
\def\mc{M_\C}
\def\mhyp{M_{\mathrm{h}}}
\def\mpar{M_{\mathrm{p}}}
\def\md{\QOPNAME{md}}

\def\Eg{E_{\mathrm{gen}}}

\def\mdag{^\dag}
\def\mstar{^*}
\def\mplus{^+}
\def\mreal{^\ddag}

\def\rr{^{(r)}}
\def\rr{^{\prime{\cdot}\prime}}

\def\CM{\Cal M}
\def\CN{\Cal N}

\def\onefam#1{\multispan2 \hss$\dim=#1$\hss}

\def\minitab#1{\hbox to\hsize\bgroup\rm
 \def\-##1{\setbox0\hbox{$00$}\hbox to\wd0{\hss$##1$\hss}}%
 \let\\\cr
 \hss\vbox\bgroup\halign\bgroup\strut\hss$##$&&#1\hss$##$\hss\cr}
\def\endminitab{\crcr\egroup\egroup\hss\egroup}

\def\minitab#1{\vcenter\bgroup\rm
 \def\-##1{\setbox0\hbox{$00$}\hbox to\wd0{\hss$##1$\hss}}%
 \let\\\cr
 \halign\bgroup\strut\hss$##$&&#1\hss$##$\hss\cr}
\def\endminitab{\crcr\egroup\egroup}
\def\otherwise{\text{else}}
\def\twenty#1{\hidewidth\le#1\hidewidth}

\let\LE=\Lambda
\let\QF=\Phi
\let\TC=\Psi
\let\QC=\Theta
\let\SC=\Delta

\def\bX{\bold{X}}
\def\bY{\bold{Y}}
\def\bQ{\bold{Q}}

\def\br{\bar r}
\def\bs{\bar s}
\def\be{\bar e}

\def\Fano{\Cal F}
\def\Fn{\QOPNAME{Fn}}
\def\val{\QOPNAME{val}}
\def\girth{\QOPNAME{girth}}

\def\L{\bold L}
\def\prop{\frak P}

\def\dual{^\vee}
\let\kk=k
\let\*=*
\let\gen=G

\def\upssrm#1{^{\scriptscriptstyle\mathrm{#1}}}
\def\0{}
\def\A{\upssrm{A}}
\def\E{\upssrm{E}}
\def\K{\upssrm{K}}

\def\St{\QOPNAME{star}}

\makeatletter
\def\p@ncil#1#2#3{\count0=#1\relax
 \ifnum\count0=0\else\ifnum\count0=1\relax#3#2\else#3\the\count0#2\fi\fi}
\def\trigpencil#1[#2,#3]{\p@ncil{#2}{\tA_2}{}%
 \p@ncil{#3}{\bA_1}{+}%
 }
\def\quadpencil#1[#2,#3,#4]{\p@ncil{#2}{\tA_3}{}%
 \p@ncil{#3}{\bA_2}{+}%
 \p@ncil{#4}{\bA_1}{+}%
 }
\def\pentpencil#1[#2,#3,#4,#5]{\p@ncil{#2}{\tA_4}{}%
 \p@ncil{#3}{\bA_3}{+}%
 \p@ncil{#4}{\bA_2}{+}%
 \p@ncil{#5}{\bA_1}{+}%
 }
\def\starpencil#1[#2,#3,#4,#5,#6,#7]{\p@ncil{#2}{\tD_4}{}%
 \p@ncil{#3}{\tA_5}{+}%
 \p@ncil{#4}{\bA_4}{+}%
 \p@ncil{#5}{\bA_3}{+}%
 \p@ncil{#6}{\bA_2}{+}%
 \p@ncil{#7}{\bA_1}{+}%
 }
\def\trig@tem#1[#2],#3]#4{(\trigpencil#2])^{#3}%
 \ifx\relax#4\else\expandafter\expandafter\expandafter\trig@tem\fi}
\def\trigs#1{\expandafter\trig@tem#1\relax}
\def\quad@tem#1[#2],#3]#4{(\quadpencil#2])^{#3}%
 \ifx\relax#4\else\expandafter\expandafter\expandafter\quad@tem\fi}
\def\quads#1{\expandafter\quad@tem#1\relax}
\def\pent@tem#1[#2],#3]#4{(\pentpencil#2])^{#3}%
 \ifx\relax#4\else\expandafter\expandafter\expandafter\pent@tem\fi}
\def\pents#1{\expandafter\pent@tem#1\relax}
\def\star@tem#1[#2],#3]#4{(\starpencil#2])^{#3}%
 \ifx\relax#4\else\expandafter\expandafter\expandafter\star@tem\fi}
\def\stars#1{\expandafter\star@tem#1\relax}
\makeatother

\title{Lines on smooth polarized $K3$-surfaces}

\author{Alex Degtyarev}

\address{%
Department of Mathematics\\
Bilkent University\\
06800 Ankara, TURKEY}

\email{
degt@fen.bilkent.edu.tr}

\thanks{%
The author was partially supported by the T\"{U}B\DOTaccent{I}TAK grant 116F211}



\begin{document}

\begin{abstract}
For each integer $D\ge3$,
we give a sharp bound on the number of lines contained in a
smooth complex $2D$-polarized
$K3$-surface in $\mathbb{P}^{D+1}$.
In the two most interesting cases of
sextics in $\mathbb{P}^4$ and octics in
$\mathbb{P}^5$, the bounds are $42$ and $36$, respectively, as conjectured in
an earlier paper.
\end{abstract}

\maketitle

%\vspace{-12pt}
%\centerline{\bf\color{red} This is a preliminary incomplete version.
%Full version is coming soon!}

\section{Introduction}\label{S.intro}

All algebraic varieties considered in the paper are over~$\C$.

\subsection{The line counting problem}\label{s.problem}
The paper deals with a very classical algebra-geometric problem, \viz. counting
straight lines in a projective surface.
We confine ourselves to the smooth $2D$-polarized $K3$-surfaces
$X\subset\Cp{D+1}$, $D\ge3$, and obtain sharp upper bounds on the
number of lines. Our primary interest are sextics in~$\Cp4$ ($D=3$) and octics
in~$\Cp5$ ($D=4$); however, the same approach gives us a complete answer for
all higher degrees/dimensions as well.

%The principal goal
%%of this paper
%is establishing sharp upper bounds on the
%maximal
%number of lines contained in a smooth $2D$-polarized $K3$-surface
%$X\subset\Cp{D+1}$, $D\ge3$.

Recall (see~\cite{Saint-Donat}) that a projectively normal smooth
$K3$-surface $X\subset\Cp{D+1}$ has projective degree~$2D$. Given a smooth
embedding $\Gf\:X\into\Cp{D+1}$, the
\emph{polarization} is the pull-back $\Gf^*\Cal{O}_{\Cp{D+1}}(1)$
regarded as a class $h\in H_2(X;\Z)$; one has $h^2=2D$.
Since each line $l\subset X$ is a
$(-2)$-curve, it is uniquely determined by its homology class
$[l]\in\NS(X)\subset H_2(X;\Z)$; to simplify the notation, we identify lines
and their classes.
In particular, the set of lines is finite;
its dual incidence graph $\Fn X$ is called the \emph{Fano graph},
and
the primitive sublattice $\Fano_h(X)\subset H_2(X;\Z)$ generated over~$\Q$
by~$h$ and all lines $l\in\Fn X$ is called the \emph{Fano
configuration} of~$X$. This polarized lattice, subject to certain
restrictions (see \autoref{th.K3}),
defines an equilinear family of
$2D$-polarized $K3$-surfaces $X\subset\Cp{D+1}$; on the other hand, the
lattice
$\Fano_h(X)$ is usually recovered from the graph $\Fn X$ (see
\autoref{s.graphs}).

The case $D=2$, \ie, that of spatial quartics, is very classical and well
known. The
%(only)
quartic
$X_{64}$ maximizing the number of lines was
constructed by F.~Schur~\cite{Schur:quartics} as early as in 1882.
The problem kept reappearing here and there ever since, but no significant
progress had been made until 1943, when B.~Segre~\cite{Segre} published a
paper asserting that 64 is indeed the maximal number of lines in a smooth
spatial quartic. Recently, S.~Rams and M.~Sch\"{u}tt~\cite{rams.schuett}
discovered a gap in Segre's argument (but not the statement);
they patched the proof and extended it to algebraically closed fields of all
characteristics other than~$2$ or~$3$.
Finally, in~\cite{DIS},
we gave an alternative proof and a number of refinements of
the statement, including the complete classification of all quartics carrying
more than 52 lines. In particular, Schur's quartic is \emph{the only one} with the
maximal number 64 of lines.

All extremal
(carrying more than 52 lines)
quartics found in~\cite{DIS} are projectively rigid, as
they are the so-called \emph{singular $K3$-surfaces}.
Recall that a $K3$-surface~$X$ is \emph{singular} if its Picard
rank is maximal, $\rho(X):=\rank\NS(X)=20$. (This unfortunate
term is not to be confused with singular
\vs. smooth
projective models of surfaces.)
Up to isomorphism, an abstract singular
$K3$-surface is determined by the oriented isomorphism type of
its \emph{transcendental lattice}
\[*
T:=\NS(X)^\perp\subset H_2(X;\Z),
\]
which is a positive definite even integral lattice of rank~$2$
(see \autoref{s.lattice}; the orientation is given by the class
$[\Go]\in T\otimes\C$ of a holomorphic $2$-form on~$X$);
we use the notation $X(T)$, see \autoref{s.notation}.
%for the singular $K3$-surface determined by a lattice~$T$
%and the single line notation $[a,b,c]$ for the lattice $\Z u+\Z v$,
%$u^2=a$, $u\cdot v=b$, $v^2=c$.
As a follow-up to~\cite{DIS} (and also motivated by~\cite{Shimada:X56}),
in~\cite{degt:singular.K3} I tried to study smooth projective models of
singular $K3$-surfaces $X(T)$ of small discriminant $\det T$.
Unexpectedly, it was discovered that Schur's quartic~$X_{64}$ can
alternatively be characterized as the only smooth spatial
model minimizing this discriminant: one has $X_{64}\cong X([8,4,8])$
(see \autoref{s.notation} for the notation)
and $\det T\ge55$ for any other smooth quartic $X(T)\subset\Cp3$.

After quartics, next most popular projective models of $K3$-surfaces are
\emph{sextics} $X\subset\Cp4$
%(regular complete intersection of a quadric and a cubic)
and \emph{octics} $X\subset\Cp5$, and the results of~\cite{degt:singular.K3}
extend to these two classes: if a singular $K3$-surface $X(T)$ admits a
smooth sextic or octic model, then $\det T\ge39$ or~$32$, respectively.
In view of the
alternative characterization of Schur's quartic~$X_{64}$ discovered
in~\cite{degt:singular.K3},
this classification,
followed by a study of the models,
%of smooth sextic and octic models
%of singular $K3$-surfaces
suggested a conjecture that
a smooth sextic $X\subset\Cp4$ (respectively, octic $X\subset\Cp5$)
may have at most~$42$
(respectively,~$36$) lines. This conjecture
%, among other results,
is proved in the present paper (the cases $D=3,4$ in
\autoref{th.main}),
even though the
original
motivating observation
that discriminant minimizing singular $K3$-surfaces maximize the number of
lines
fails
already for degree~$10$ surfaces $X\subset\Cp6$
(see \autoref{th.min.det}).

Each sextic in~$\Cp4$
is a regular complete intersection of a quadric and a cubic. I am not aware
of any previously known interesting examples of large configurations of lines
in such surfaces. The maximal number~$42$ of lines is attained at a
$1$-parameter family containing the discriminant minimizing surfaces
$X(2,1,20])$ and $X([6,3,8])$.

Most octics in~$\Cp5$ are also regular complete intersections:
they are the so-called \emph{triquadrics}, \ie,
intersections of three quadrics.
The most well-known example is the \emph{Kummer family}
\[
\sum_{i=0}^5z_i=\sum_{i=0}^5a_iz_i=\sum_{i=0}^5a_i^2z_i=0,\quad a_i\in\C,
\label{eq.Kummer}
\]
whose generic members contain $32$ lines: famous Kummer's $16_6$
configuration (see, \eg, Dolgachev~\cite{Dolgachev:book}).
There are four other, less symmetric, configurations of $32$ lines, either
rigid or realized by $2$-parameter families. However, $32$ is not the
maximum: there also are configurations with $33$, $34$, and $36$ lines (see
\autoref{tab.main} on \autopageref{tab.main}).

The space of octics contains a
%codimension~$1$ subspace
divisor composed by
surfaces
that need one cubic defining equation
(see~\cite{Saint-Donat} and \autoref{s.polarized} below).
We call these octics \emph{special} and show that they do stand out in
what concerns the line counting problem.
Large Fano graphs of special octics
(\vs. triquadrics)
are described by \autoref{th.special}.

\subsection{Common notation}\label{s.notation}
We use the following notation for particular integral lattices
of small rank (see \autoref{s.lattice}):
\roster*
\item
$\bA_p$, $p\ge1$, $\bD_q$, $q\ge4$, $\bE_6$, $\bE_7$, $\bE_8$
are the \emph{positive definite} root
lattices generated by the indecomposable root systems of the same name
(see~\cite{Bourbaki:Lie});
\item
$[a]:=\Z u$, $u^2=a$, is a lattice of rank~$1$;
\item
$[a,b,c]:=\Z u+\Z v$, $u^2=a$, $u\cdot v=b$, $v^2=c$, is a lattice of
rank~$2$; when it is positive definite, we assume that $0<a\le c$ and
$0\le2b\le a$: then, $u$ is a shortest vector, $v$ is a next shortest one,
and the triple $(a,b,c)$ is unique;
\item
$\bU:=[0,1,0]$ is the \emph{hyperbolic plane}.
%: the unimodular even lattice of rank~$2$.
\endroster
Besides,
$L(n)$, $n\in\Z$, is the lattice obtained by the scaling of
a given lattice~$L$, \ie, multiplying
all values of the
form by a fixed integer~$n\ne0$.

To simplify the statements, we let the girth of a forest equal to infinity,
so that
the inequality $\girth(\graph)\ge m$ means that $\graph$ has no cycles
of length less than~$m$.

When describing lists of integers,
$a\dts b$
%means the full range $a,a+1,\dots,b-1,b$.
means the full range $\Z\cap[a,b]$.

We denote by $\{a\bmod d\}$
the arithmetic progression $a+nd,\ldots$, $n\ge0$, and use the shortcut
$\{a_1,\ldots,a_r\bmod d\}:=\bigcup_i\{a_i\bmod d\}$ for finite unions.

\subsection{Principal results}\label{s.results}
Given a field $K\subset\C$ and an integer $D\ge2$, let $\mm_K(D)$ be the
maximal number of lines defined over~$K$ that a smooth $2D$-polarized
$K3$-surface $X\subset\Cp{D+1}$ defined over~$K$ may have.

Most principal results of the paper are collected in \autoref{tab.main};
%and stated more precisely in
precise statements are found in the
several theorems below.
(For the sake of completeness, we also
cite some results of~\cite{DIS} concerning quartics, \ie, $D=2$.)
\table\rm
\caption{Exceptional configurations (see \autoref{th.main})}\label{tab.main}
\def\MMi#1{\expandafter\MMii#1\endMM}
\def\MMii#1,#2\endMM{$#1,#2$}
\def\-{\rlap{$\mdag$}}
\def\*{\rlap{$\mstar$}}
\def\r{\rlap{$\mreal$}}
\def\2#1{\!{\le}2#1\!}
\def\={\relax\afterassignment\ul\count0=}
\def\ul{\underline{\the\count0}}
\def\u#1{\mathbf{#1}}
\hbox to\hsize{\hss\vbox{\halign{\strut\quad\hss$#$\hss\quad&
 \hss$#$\hss\quad&\hss$#$\hss\quad&\hss$#$\hss\quad&\hss$#$\hss\quad&
 \hss$#$\hss\quad&\hss$#$\hss\quad&\hss$#$\hss\quad&$#$\hss\quad\cr
\noalign{\hrule\vspace{2pt}}
D&\mm,\bm&\!\mr\!&\graph&\ls|\Aut\graph|&\det&(r,c)&
 \ls|\Aut X|&T:=\Fano_h(\graph)^\perp\cr
\noalign{\vspace{1pt}\hrule\vspace{2pt}}
 2  &64,52&56&\bX_{64}    &4608&\=48&(1,0)&1152&[8,4,8]\cr
           &&&\bX_{60}'   & 480&  60&(1,0)& 120&[4,2,16]\cr
           &&&\bX_{60}''  & 240&  55&(0,1)& 120&[4,1,14]\*\cr
           &&&\bX_{56}    & 128&  64&(0,1)&  32&[8,0,8]\*\cr
           &&&\bY_{56}\r  &  64&  64&(1,0)&  32&[2,0,32]\cr
           &&&\bQ_{56}    & 384&  60&(1,0)&  96&[4,2,16]\cr
           &&&\bX_{54}    & 384&  96&(1,0)&  48&[4,0,24]\cr
           &&&\bQ_{54}    &  48&  76&(1,0)&   8&[4,2,20]\cr
 3  &42,36&42&\TC_{42}\r  & 432&\=39&\onefam1  &\bA_2\oplus[-18]\cr
           &&&\TC_{38}\r  &  32&  48&\onefam1  &\bA_2\oplus[-24]\cr
 4  &36,30&34&\QC_{36}'   &  64&\=32&(1,0)&  16&[4,0,8]\cr
           &&&\QC_{36}''  & 576&  36&(1,0)& 144&[6,0,6]\cr
           &&&\QC_{34}'\r &  96&    &\onefam1  &\bU(2)\oplus[12]\cr
           &&&\QC_{33}    & 192&  80&(1,0)&  24&[8,4,12]\cr
           &&&\TC_{33}    &6912&  36&\onefam2  &\bU^2(3)\cr
           &&&\QC_{32}\0  &  96&  60&(1,0)&  24&[4,2,16]\cr
           &&&\QC_{32}'\r & 384&   &\onefam2  &\bU(2)\oplus[-4]\oplus[4]\cr
           &&&\QC_{32}''\r& 512&    &\onefam2  &\bU(2)\oplus\bU(4)\cr
           &&&\QC_{32}''' & 768&  36&\onefam2  &\bU^2(3)\cr
           &&&\QC_{32}\K\r&23040&\=32&\onefam3 &\bU^2(2)\oplus[-4]\cr
 5\-&30,28&28&\QF_{30}'   & 240& 100&(1,0)&  60&[10,0,10]\cr
           &&&\QF_{30}''  &  40&  75&(1,0)&  10&[10,5,10]\cr
           &&&\QF_{30}''' &  24&  36&(1,0)&  12&[4,2,10]\cr
 6  &36,28&36&\QF_{30}'   & 240&  60&(1,0)&  60&[8,2,8]\cr
           &&&\QF_{36}''\r&1440&\=15&(1,0)& 720&[2,1,8]\cr
 7  &30,26&26&\QF_{30}'   & 240&\=20&(1,0)& 120&[4,2,6]\cr
           &&&\SC_{27}''  &  72&  99&(1,0)&  18&[6,3,18]\cr
 8  &32,24&24&\SC_{32}'   &2304&\=12&(1,0)&1152&[4,2,4]\cr
           &&&\SC_{26}    &  24&  60&(1,0)&   6&[4,2,16]\cr
           &&&\SC_{25}''  &  16&  60&(1,0)&   4&[4,2,16]\cr
           &&&\QF_{25}    &  80&    &\onefam1  &\bU(5)\oplus[4]\cr
 9  &25,24&25&\QF_{25}\r  &  80&\=15&\onefam1  &\bU(5)\oplus[2]\cr
10\-&25,24&25&\SC_{25}'   &  24&  96&(0,1)&   6&[4,0,24]\cr
           &&&\QF_{25}\r  &  80&    &\onefam2  &\bU(5)\oplus\bU\cr
           &&&\LE_{25}    & 144& 140&(1,0)&  12&[12,2,12]\cr
14  &28,24&28&\LE_{28}\r  & 336& \=7&(1,0)& 336&[2,1,4]\cr
\noalign{\vspace{1pt}\hrule}
\crcr}}\hss}
\endtable
In the table, we list:
\roster*
\item
the degrees $h^2=2D$ for which extra information is available
(marking with a $\mdag$ the two values that are
special in the sense of \autoref{th.min.det}),
\item
the bounds $\mm:=\mm(D)$ and $\bm:=\bm(D)$ used in \autoref{th.main}
(according to which, $M(D)=\mc(D)$ for all values of~$D$ in the table),
and
\item
the maximal number $\mr:=\mr(D)$ of real lines (see \autoref{th.main.real}).
\endroster
Then, for each value of~$D$, we list, line-by-line, all Fano graphs
$\graph:=\Fn X$ containing more than $\bm(D)$ lines
(the notation is explained below), marking with a~$\mreal$ those
realized by real lines in real surfaces (see \autoref{th.main.real}) and indicating
\roster*
\item
the order $\ls|\Aut\graph|$ of the full automorphism group of~$\graph$ and
\item
the transcendental lattices
$T:=\Fano_h(\graph)^\perp=\NS(X)^\perp\in H_2(X;\Z)$ of generic smooth
$2D$-polarized $K3$-surfaces~$X$ with $\Fn X\cong\graph$
(marking with a $\mstar$ the lattices resulting in pairs of complex
conjugate equilinear families).
\endroster
For the rigid configurations ($\rank T=2$), we list, in addition,
\roster*
\item
the determinant $\det T$ (underlining the minimal ones, see
\autoref{th.min.det}),
\item
the numbers
$(r,c)$ of, respectively, real and pairs of complex conjugate projective
isomorphism classes of surfaces~$X$ with $\Fn X\cong\graph$, and
\item
the order $\ls|\Aut X|$ of the
%(common)
group of projective automorphisms of~$X$.
\endroster
Each of the few non-rigid configurations~$\graph$ appearing in the table
is realized by a single connected
equilinear deformation family~$\CM_D(\graph)$; we indicate the dimension
%modulo the projective group $\PGL(D+2,\C)$.
$\dim\bigl(\CM_D(\graph)/\!\PGL(D+2,\C)\bigr)=\rank T-2$ and,
when known, the minimum of the discriminants
$\det T$ of the singular $K3$-surfaces $X(T)\in\CM_D(\graph)$.

If $D=2$, we use the notation for the extremal configurations introduced
in~\cite{DIS}; otherwise, we refer to the isomorphism classes of the
Fano graphs introduced and discussed in more details
elsewhere in the paper.
In both cases, the subscript
is the number of lines in the configuration.
For technical reasons,
we subdivide Fano graphs into several
classes (see \autoref{s.taxonomy}) and study them separately,
obtaining more refined bounds for each class.
Thus, a Fano graph $\graph:=\Fn X$ and the
%corresponding
configuration $\Fano_h(\graph):=\Fano_h(X)$ are called
\roster*
\item
\emph{triangular} (the $\TC_*$-series, see \autoref{th.trig}),
if $\girth(\graph)=3$; all extremal quartics also fall into this class,
\item
\emph{quadrangular} (the $\QC_*$-series, see \autoref{th.quad}),
if $\girth(\graph)=4$,
\item
\emph{pentagonal} (the $\QF_*$-series, see \autoref{th.a4}),
if $\girth(\graph)=5$,
\item
\emph{astral} (the $\SC_*$-series, see \autoref{th.astral}),
if $\girth(\graph)\ge6$ and $\graph$ has a vertex~$v$ of valency
$\val v\ge4$.
\endroster
All other graphs are \emph{locally elliptic} (the $\LE_*$-series, see
\autoref{s.locally.elliptic}), \ie, one has $\val v\le3$ for each
vertex $v\in\graph$.

In our notation for particular graphs/configurations,
the subscript always stands for
the number of vertices/lines.
The precise description of all graphs ``named'' in the paper is
available
electronically
(in the form of \texttt{GRAPE} records)
in~\cite{degt:Fano.graphs};
in most cases, the implicit reference to~\cite{degt:Fano.graphs} is, in fact,
\emph{the definition} of the graph.

\subsection{The bounds}\label{s.bounds}
Geometrically, apart from the spatial quartics, the two most interesting projective
models of $K3$-surfaces are sextics in $\Cp4$ and octics (especially
triquadrics) in $\Cp5$. However, it turns out that the structure of the Fano
graphs simplifies dramatically when
the degree $h^2=2D$ grows, and one can easily obtain the sharp bounds
$\mc(D)$ and $\mr(D)$ for all values of~$D$.
Below, these bounds are stated for
$D$ small and $D\to\infty$;
%postponing till \autoref{s.locally.elliptic}
the rather erratic precise values of
\[*
21\le\mc(D)\le24,\ D>20
\quad\text{and}\quad
19\le\mr(D)\le24,\ D>16
\]
are postponed till \autoref{s.locally.elliptic}
(see Corollaries~\ref{cor.le.mc} and~\ref{cor.le.mr}, respectively).

%$21\le\mc(D)\le24$ for $D>20$ (see \autoref{cor.le.mc})
%and $19\le\mr(D)\le24$ for $D>16$ (see \autoref{cor.le.mr}).

For the next theorem, given a
degree $h^2=2D$, let
$\mm:=\mm(D)$ and $\bm:=\bm(D)$ be as in \autoref{tab.main} or
$\mm(D)=\bm(D)=24$ if $D$ is not found in the table.
As in~\cite{DIS}, in addition to the upper bound
$\ls|\Fn X|\le\mm(D)$, we give a complete classification of all large
(close to maximal)
configurations of lines.

\theorem[see \autoref{proof.main}]\label{th.main}
Let $X\subset\Cp{D+1}$ be a smooth $2D$-polarized $K3$-surface.
Then one has $\ls|\Fn X|\le\mm(D)$.
More precisely, $\ls|\Fn X|\le\bm(D)$
unless $\graph:=\Fn X$ is one of the exceptional graphs \rom(configurations\rom)
listed in \autoref{tab.main}.

Both bounds are sharp
%\rom(and, in particular, $\mc(D)=M$\rom)
whenever $2\le D\le20$ or $D\in\{1,2,4,6,7,10\bmod12\}$\rom;
in particular, one has $\mc(D)=\mm(D)$ for these values.
\endtheorem

\addendum[see \autoref{proof.main}]\label{ad.main}
The complete set $\{\ls|\Fn X|\}$ of values taken by the line count
of a smooth $2D$-polarized $K3$-surface $X\subset\Cp{D+1}$ is as follows\rom:
\roster*
\item
if $D=2$, then $\{\ls|\Fn X|\}=\{0\dts52,54,56,60,64\}$
\rom(see~\cite{DIS}\rom)\rom;
\item
if $D=3$, then $\{\ls|\Fn X|\}=\{0\dts36,38,42\}$\rom;
\item
if $D=4$, then $\{\ls|\Fn X|\}=\{0\dts30,32,33,34,36\}$.
\endroster
%If $D=4$, the extremal
%Fano graphs~$\QC_*^*$ are
%realized by triquadrics,\mnote{\todo: to edit the proof}
%but $\TC_{33}$ is not.
%If $D\ge13$, any configuration maximizing the number of lines is locally
%elliptic.
%
%If $D\ge17$, any configuration with more than $20$ lines is locally
%elliptic.
\endaddendum

As mentioned at the end of \autoref{s.problem}, the Fano graphs of special
octics differ from those of triquadrics. They are described by the following
theorem, which implies,
in particular, that the graphs $\QC_*^*$ in \autoref{tab.main} are
realized \emph{only} by triquadrics, whereas $\TC_{33}$ is realized only by
special octics.

\theorem[see \autoref{proof.special}]\label{th.special}
Let $X\subset\Cp5$ be a smooth special octic, and
assume that $\ls|\Fn X|\ge20$. Then one of the
following statements holds\rom:
\roster
\item
$\girth(\Fn X)=3$ and
$\Fn X\cong\TC_{33}$, $\TC_{29}$, $\TC_{27}'$, $\TC_{27}''$
or $\ls|\Fn X|\le25$, or
\item
$\girth(\Fn X)=4$ and $\Fn X\cong\QC_{21}$ has a biquadrangle \rom(see
\autoref{def.biquad}\rom).
\endroster
Conversely, if $X\subset\Cp5$ is a smooth octic such that either
$\girth(\Fn X)=3$ or $\Fn X$ has a biquadrangle, then $X$ is special.
\endtheorem

The threshold $\ls|\Fn X|\ge20$ in \autoref{th.special} is optimal: there
are graphs with $19$ lines
realized by both triquadrics and special octics (see \autoref{prop.dual}).

For completeness, in \autoref{s.hyperelliptic} we discuss also ``lines'' in
the hyperelliptic models
\[*
X\to\Sigma_{D-2}\into\Cp{D+1},\quad D\ge3,
\]
which turn out very simple and quite similar to special octics (see
\autoref{th.hyperelliptic}; the case $D=2$ is considered
in~\cite{degt:singular.K3}).
The conjectural bound of $144$ lines in double planes $X\to\Cp2$
(see~\cite{degt:singular.K3}) is left to a subsequent paper.

As a by-product of the classification of large Fano graphs, we obtain bounds
on the maximal number $\mr(D)$ of real lines in a real
$K3$-surface.

\theorem[see \autoref{proof.main.real}]\label{th.main.real}
The number $\mr(D)$ is as given by \autoref{tab.main}. If $D$ is not found
in the table, then $\mr(D)\le24$\rom;
this bound is sharp for $D\le16$.
\endtheorem

\subsection{Asymptotic bounds}
If $D$ is large,
%($D>620$),
all Fano graphs become very simple, \viz. disjoint
unions of (affine) Dynkin diagrams, and we have
\[*
\limsup_{D\to\infty}\mc(D)=24,\qquad
\limsup_{D\to\infty}\mr(D)=21.
\]
More precisely, we have the following two asymptotic statements.

\theorem[see \autoref{proof.D>>0}]\label{th.D>>0}
Let $X\subset\Cp{D+1}$ be a smooth $2D$-polarized $K3$-surface, and assume
that $D\gg0$. Then, either
%Assume that $D\gg0$. Then, for any smooth $2D$-polarized $K3$-surface
%$X\in\Cp{D+1}$, either
\roster*
\item
$\Fn X$ is a disjoint union of Dynkin diagrams, and $\ls|\Fn X|\le19$,
or
\item
all lines in~$X$ are fiber components of an elliptic pencil,
and $\ls|\Fn X|\le24$.
\endroster
The precise sharp bound $\mc(D)$, $D\gg0$, is as
follows\rom:
\roster*
\item
if $D=1\bmod3$ or $D=2\bmod4$, then $\mc(D)=24$\rom;
\item
otherwise, if $D=0\bmod4$ or $D=2\bmod5$, then $\mc(D)=22$\rom;
\item
in all other cases, $\mc(D)=21$.
\endroster
\endtheorem

\theorem[see \autoref{proof.D>>0.real}]\label{th.D>>0.real}
If $D\gg0$, the sharp bound $\mr(D)$ is as follows\rom:
\roster*
\item
if $D=0\bmod3$, $D=1\bmod5$, or $D=0\bmod7$, then $\mr(D)=21$\rom;
\item
otherwise, if $D\notin\{17,113\bmod120\}$, then $\mr(D)=20$\rom;
\item
in all other cases, $\mr(D)=19$.
\endroster
\endtheorem

The conclusion of \autoref{th.D>>0} holds to full extent (all Fano graphs are
elliptic or parabolic) for $D>620$ (see \autoref{rem.min.parabolic}). The
expression for $\mc(D)$ is valid for $D>45$ (see \autoref{cor.le.mc}), and the
expression for $\mr(D)$ in \autoref{th.D>>0.real} is valid for $D>257$ (\cf.
\autoref{cor.le.mr}).

The structure of
elliptic pencils carrying many lines is
given by \autoref{cor.parabolic} and \autoref{tab.parabolic} on
\autopageref{tab.parabolic}.

\subsection{Discriminant minimizing surfaces}\label{s.minimal}
Given a degree $h^2=2D$, one can pose a question on the minimal discriminant
\[*
\md(D):=
\min\bigl\{\det T\bigm|\text{$\rank T=2$ and there is a smooth model
 $X(T)\subset\Cp{D+1}$}\bigr\}
\]
%$\det T$
of a singular $K3$-surface $X(T)$ admitting a smooth model
%$X(T)\subset\Cp{D+1}$.
of degree $2D$.
In~\cite{degt:singular.K3}, this question was answered
for $D=2,3,4$, and the conjecture on the maximal number of lines in sextic
and octic models was based on the assumption that the surfaces
minimizing the discriminant
should also maximize the number of lines. Ironically, although the statement
of the conjecture does hold, its motivating assumption fails for the very
next value $D=5$.

The proof of the next theorem is omitted as it repeats,
almost literally, the same
computation as in~\cite{degt:singular.K3}.

\theorem\label{th.min.det}
For each degree $h^2=2D$ as in \autoref{tab.main}, there is a unique genus
of discriminant minimizing transcendental lattices~$T_0$. If $D\ne5,10$, the
discriminant minimizing surfaces have
smooth models $X(T_0)\subset\Cp{D+1}$ maximizing
the number of lines \rom(underlined in the table\rom).
For the two exceptional values, one has\rom:
\roster*
\item
$\md(5)=32$ and $X([2,0,16])\subset\Cp6$ has $28$ lines
\rom(see $\QF_{28}$ in~\cite{degt:Fano.graphs}\rom)\rom;
\item
$\md(10)=15$ and $X([4,1,4])\subset\Cp{11}$ has $24$
\rom(see $\SC_{24}$ in~\cite{degt:Fano.graphs}\rom)
or $20$ lines.
\done
\endroster
\endtheorem

Unless $D=3$ or $9$, the discriminant minimizing transcendental lattices~$T_0$
are shown in \autoref{tab.main}. If $D=3$, one has $T_0=[2,1,20]$ or
$[6,1,8]$ and
$X(T_0)\in\CM_3(\TC_{42})$
%$\Fn X(T_0)\cong\TC_{42}$
(see~\cite{degt:singular.K3});
if $D=9$, then $T_0=[2,1,8]$ and
$X(T_0)\in\CM_9(\QF_{25})$.

\subsection{Contents of the paper}
In \autoref{S.prelim}, after introducing the necessary definitions and results concerning
lattices, polarized lattices and their relation to $K3$-surfaces, and graphs,
we state the arithmetical reduction of the line counting problem that is used
throughout the paper.
An important statement is \autoref{prop.hyperbolic}: it rules out
most hyperbolic graph for most degrees.
As an application, in \autoref{s.valency} we prove
several sharp bounds on the valency of a line in a smooth polarized
$K3$-surface.

In \autoref{S.elliptic}, we discuss parabolic graphs (essentially, affine
Dynkin diagrams) and their relation to elliptic pencils on a
$K3$-surface.
As a counterpart of \autoref{prop.hyperbolic}, we prove
\autoref{prop.parabolic} stating the periodicity
of the set of geometric degrees for each parabolic graph.
In \autoref{proof.D>>0}, we prove Theorems~\ref{th.D>>0}
and~\ref{th.D>>0.real}.

%The rest of the paper deals mainly with hyperbolic Fano graphs.
Having the conceptual part settled, we start the classification of hyperbolic
Fano graphs.
We divide them into groups, according to the minimal fiber, and
study case-by-case:
locally elliptic (see \autoref{S.locally.elliptic}), astral (see
\autoref{S.astral}), pentagonal (see \autoref{S.pent}),
quadrangular (see \autoref{S.quad}), and
triangular (see \autoref{S.trig}). For each group, we obtain a finer
classification of large graphs in small degrees, which is stated as a
separate theorem.
Combining these partial statements, we prove Theorems~\ref{th.main}
and~\ref{th.main.real} in \autoref{proof.main}.
%The principal results of the paper (Theorems~\ref{th.main}
%and~\ref{th.main.real}) are proved in \autoref{proof.main}.
%by combining the
%bound sound for each class.

In \autoref{S.astral}, we discuss also special octics and hyperelliptic
models of $K3$-surfaces.

In \autoref{S.algorithms}, we outline a few technical details
concerning the
algorithms used in the proofs when enumerating large Fano graphs.

\subsection{Acknowledgements}
I am grateful to S{\l}awomir Rams and Ichiro Shimada for fruitful and
motivating discussions. My special gratitude goes to Dmitrii Pasechnik
for his indispensable help in
``identifying'' some of the Fano graphs.
This paper was revised and finalized during my stay at the
Abdus Salam International Centre for Theoretical Physics,
Trieste; I would like to thank this institution and its
friendly staff for their warm hospitality and excellent working conditions.
%\mnote{\todo: ICTP}
%discovered in the paper.

\section{Preliminaries}\label{S.prelim}

We start with an arithmetical reduction of the line counting problem (which,
in fact, is well known, see Theorems~\ref{th.K3} and~\ref{th.real}) and
describe the technical tools used to detect geometric configurations of
lines. Then, Fano graphs are subdivided into \emph{elliptic},
\emph{parabolic}, and \emph{hyperbolic}. The first two classes are not very
interesting, as they may contain at most $19$ or $24$ lines, respectively.
(Parabolic graphs are treated in more detail in \autoref{S.elliptic}.) Thus,
the rest of the paper deals mainly with hyperbolic Fano graphs. The very first
results, obtained in this section, are a degree bound (see
\autoref{prop.hyperbolic}) and valency bounds (see
Propositions~\ref{prop.valency} and~\ref{prop.valency.trig}).

\subsection{Lattices}\label{s.lattice}
A \emph{lattice} is a finitely generated free abelian group~$L$ equipped with
a symmetric bilinear form $b\:L\otimes L\to\Z$ (which is usually understood
and omitted from the notation). We abbreviate $x\cdot y:=b(x,y)$ and
$x^2:=b(x,x)$. In this paper, all lattices are \emph{even}, \ie,
$x^2=0\bmod2$ for all $x\in L$.

The \emph{determinant} $\det L$ is the determinant of the Gram matrix of~$b$
in any integral basis, and the \emph{kernel} $\ker L$ (as opposed to the
kernel $\Ker\Gf$ of a homomorphism~$\Gf$) is the subgroup
\[*
\ker L:=L^\perp=\bigl\{x\in L\bigm|\text{$x\cdot y=0$ for all $y\in L$}\bigr\}.
\]
An integral lattice~$L$ is \emph{unimodular} if $\det L=\pm1$; it is
\emph{nondegenerate} if $\det L\ne0$ or, equivalently, $\ker L=0$.

The \emph{inertia indices} $\Gs_{\pm,0}(L)$ are the classical inertia indices
of the quadratic form $L\otimes\R$. Clearly, $\Gs_0(L)=\rank\ker L$.
The inertia index $\Gs_+(L)$ (the sum
$\Gs_+(L)+\Gs_0(L)$) is the dimension of maximal positive definite
(respectively, semidefinite) subspaces in $L\otimes\R$; a similar
statement holds for $\Gs_-(L)$ and negative subspaces. Hence,
these quantities are monotonous:
\[
\text{if $S\subset L$, then $\Gs_\pm(S)\le\Gs_\pm(L)$ and
 $\Gs_\pm(S)+\Gs_0(S)\le\Gs_\pm(L)+\Gs_0(L)$}.
\label{eq.monotonous}
\]
A lattice~$L$ is called \emph{hyperbolic} if $\Gs_+(L)=1$.

An \emph{extension} of a lattice~$L$ is any overlattice $S\supset L$. If $L$
is nondegenerate, there is a canonical inclusion of~$L$ to the
\emph{dual group}
\[*
L\dual=\bigl\{x\in L\otimes\Q\bigm|\text{$x\cdot y\in\Z$ for all $y\in L$}\bigr\},
\]
which inherits from~$L$ a $\Q$-valued quadratic form.
Any finite index extension of~$L$ is a subgroup of $L\dual$. In general,
extensions of~$L$ can be described in terms of the \emph{discriminant group}
$\discr L:=L\dual\!/L$: this is a finite abelian group equipped with
%equipped with the induced
the
nondegenerate $(\Q/2\Z)$-valued
quadratic form $(x\bmod L)\mapsto x^2\bmod2\Z$.
Since we omit the details of the computation, we merely refer to the
original paper~\cite{Nikulin:forms}.

%We denote by~$\L$ a fixed representative of the isomorphism class of
%the intersection lattice $H_2(X;\Z)$ of a $K3$-surface~$X$: this is the only
%(up to isomorphism)
%unimodular even lattice of rank~$22$ and signature~$-16$; in other words,
%$\L\cong\bU^3\oplus2\bE_8^2(-1)$.

\subsection{Polarized lattices}\label{s.polarized}
Given an even integer $D\ge2$, a \emph{$2D$-polarized lattice} is a
nondegenerate hyperbolic lattice~$S$ equipped with a distinguished vector
$h\in S$ such that $h^2=2D$.
For such a lattice $(S,h)$ we can define its set
of \emph{lines}
\[*
\Fn(S,h):=\bigl\{v\in S\bigm|v^2=-2,\ v\cdot h=1\bigl\}.
\]
This set is finite. Both~$D$ and~$h$ are often omitted from the notation.
A polarized lattice $(S,h)$ is called a \emph{configuration} if it is
generated over~$\Q$ by $h$ and $\Fn(S,h)$.

Denote by~$\L$ a fixed representative of the isomorphism class of
the intersection lattice $H_2(X;\Z)$ of a $K3$-surface~$X$: this is the only
(up to isomorphism)
unimodular even lattice of rank~$22$ and signature~$-16$; in other words,
$\L\cong\bU^3\oplus2\bE_8^2(-1)$.

\definition\label{def.geometric}
Depending on the geometric problem, we define ``bad'' vectors in a
polarized lattice $(S,h)$ as vectors $e\in S$ with one of the following
properties:
\roster
\item\label{bad.exceptional}
$e^2=-2$ and $e\cdot h=0$ (\emph{exceptional divisor}),
\item\label{bad.hyperelliptic}
$e^2=0$ and $e\cdot h=2$ (\emph{quadric pencil}), or
\item\label{bad.cubic}
$e^2=0$ and $e\cdot h=3$ (\emph{cubic pencil}).
\endroster
Usually, only vectors as in \iref{bad.exceptional}
or~\iref{bad.hyperelliptic} are excluded, whereas vectors as
in~\iref{bad.cubic} are to be excluded if $D=4$ and we consider triquadrics
rather than all octic surfaces.

Respectively, a polarized lattice $(S,h)$ is called
\roster*
\item
\emph{admissible}, if it contains no ``bad'' vectors,
\item
\emph{geometric}, if it is admissible and has a primitive embedding
to~$\L$, and
\item
\emph{subgeometric}, if it admits a geometric finite index extension.
\endroster
%An admissible polarized lattice $(S,h)$ generated over~$\Q$
%by~$h$ and $\Fn(S,h)$ is called a \emph{configuration}.
\enddefinition

Given a smooth $2D$-polarized $K3$-surface $X\subset\Cp{D+1}$,
we denote by $\Fn X$ the set of lines contained in~$X$
(the \emph{Fano graph}) and define
the \emph{Fano configuration} $\Fano_h(X)$
as the primitive sublattice $S\subset H_2(X;\Z)$
generated over~$\Q$ by the polarization~$h$ and the classes $[l]$ of all
lines $l\in\Fn X$.
The next statement is straightforward and well known;
it follows from~\cite{Saint-Donat} (see also
\cite[Theorem 3.11]{DIS} or \cite[Theorem 7.3]{degt:singular.K3}).

\theorem\label{th.K3}
A $2D$-configuration $(S,h)$ is isomorphic to the Fano configuration
$\Fano_h(X)$ of a smooth $2D$-polarized $K3$-surface~$X$ if and only if
$(S,h)$ is geometric\rom; in this case, an isomorphism
$(S,h)\to\Fano_h(X)$
induces a bijection $\Fn(S,h)\to\Fn X$.

Given a primitive embedding $S\into\L$, all
$2D$-polarized
surfaces~$X$ such that
\[
\bigl(H_2(X;\Z),\Fano_h(X),h\bigr)\cong\bigl(\L,S,h\bigr)
\label{eq.family}
\]
constitute one or two \rom(complex conjugate\rom)
equilinear deformation families\rom; modulo
the projective group
$\PGL(D+2,\C)$,
their dimension
equals $20-\rank S$.
\pni
\endtheorem

\remark\label{rem.rigid}
The number of families in \autoref{th.K3} equals one or two depending on
whether $\L$ does or, respectively, does not admit an automorphism
preserving~$h$ and~$S$ (as a set) and reversing the orientation of maximal
positive definite subspaces in $\L\otimes\R$.
A configuration~$S$ is called \emph{rigid} if
$\rank S=20$: such configurations are realized by finitely many projective
equivalence classes of surfaces.
\endremark

\remark\label{rem.D=8}
According to~\cite{Saint-Donat}, if $D=4$, one should also require that
$h\in S$ is a primitive vector. However, this condition holds automatically if
$\Fn S\ne\varnothing$. Besides, in this case one can
also distinguish between all
octic surfaces and triquadrics; as explained above, the difference is in the
definition of the admissibility: cubic pencils as in
\autoref{def.geometric}\iref{bad.cubic} are
either allowed (special octics) or excluded (triquadrics).
\endremark

We use \autoref{th.K3} in conjunction with the following two algorithms.

\algorithm\label{alg.lines}
There is an effective algorithm (using the enumeration of vectors of given
length in a definite lattice, implemented as \texttt{ShortestVectors} in
\GAP~\cite{GAP4}) detecting whether a lattice $(S,h)$ is admissible
and computing the set $\Fn(S,h)$.
\endalgorithm

\algorithm\label{alg.Nikulin}
A nondegenerate lattice admits but finitely many finite index extensions,
which can be effectively enumerated
(see \cite[Proposition 1.4.1]{Nikulin:forms}). Then, one can use
\autoref{alg.lines} to select the admissible extensions
(note that this is not automatic) and, for each such
extension, use \cite[Corollary 1.12.3]{Nikulin:forms} to detect whether it
admits a primitive embedding to~$\L$.
(A necessary, but not sufficient condition is that $\rank S\le20$.)
This computation, using Nikulin's techniques of discriminant forms,
can also be implemented
in \GAP~\cite{GAP4}.
\endalgorithm

\subsection{Real configurations}\label{s.real}
For a real (\ie, invariant under the standard complex conjugation
involution) smooth $K3$-surface $X\subset\Cp{D+1}$, denote by
$\Fn_\R X\subset\Fn X$ the set of real lines contained in~$X$.
As explained in~\cite{DIS}, $X$ can be deformed to a real surface $X'$ such
that $\Fn X'=\Fn_\R X'\cong\Fn_\R X$, \ie, all lines in~$X'$ are real.
Hence, when
computing the bound $\mr(D)$, we can assume \emph{all} lines real.

\theorem[{\cf.~\cite[Lemma 3.10]{DIS}}]\label{th.real}
Consider a geometric configuration $(S,h)$ and a primitive embedding
$S\into\L$. Then, family~\eqref{eq.family}
%the family of $2D$-polarized surfaces~$X$ with
%\[*
%\bigl(H_2(X;\Z),\Fano_h(X),h\bigr)\cong\bigl(\L,S,h\bigr)
%\]
%\rom(\cf. \autoref{th.K3}\rom)
contains a real surface with all lines real if
and only if the generic transcendental lattice $T:=S^\perp$ has a
sublattice isomorphic to $[2]$ or $\bU(2)$.
\pni
\endtheorem

\subsection{Configurations as graphs}\label{s.graphs}
One can easily show that, if a polarized lattice~$S$ is admissible, for any
two lines $u,v\in\Fn S$ one has $u\cdot v=0$ or~$1$; respectively, we say
that $u$ and~$v$ are \emph{disjoint} or \emph{intersect}. Therefore, it is
convenient to regard $\Fn S$ as a graph: the vertices are the lines
$v\in\Fn S$, and two vertices are connected by an edge if and only if the
lines intersect. We adopt the graph theoretic terminology; for example,
we define the \emph{valency} $\val v$ of a line $v\in\Fn S$ as the number of
lines $u\in\Fn S$ such that $u\cdot v=1$.

The term \emph{subgraph} always means an induced subgraph.
To simplify statements, introduce the \emph{relative valency}
\[*
\val_\fiber l:=\#\bigl\{a\in\fiber\bigm|a\cdot l=1\bigr\}
\]
with respect to a subgraph $\fiber\subset\graph$
and, given two subgraphs $\fiber_1,\fiber_2\subset\graph$, define
\[*
\fiber_1\*\fiber_2:=\max\bigl\{\val_{\fiber_i}l_j\bigm|
 \text{$l_j\in\fiber_j$, $(i,j)=(1,2)$ or $(2,1)$}\bigr\}.
\]
For example, $\fiber\*\fiber=0$ if and only if
the subgraph $\fiber\subset\graph$ is
\emph{discrete}.

Conversely, to an abstract graph~$\graph$ (loop free and without multiple
edges) we can associate three lattices
\[*
\Z\graph,\qquad
\Fano(\graph):=\Z\graph/\ker,\qquad
\Fano_h(\graph):=(\Z\graph+\Z h)/\ker.
\]
Here, $\Z\graph$ is freely generated by the vertices $v\in\graph$, so that
$v^2=-2$ and $u\cdot v=1$ or~$0$ if $u$ and~$v$ are, respectively, adjacent
or not. For the last lattice $\Fano_h(\graph)$, the sum is direct, but not
orthogonal; the even integer $2D:=h^2\ge4$ is assumed fixed in advance, and we
have $v\cdot h=1$ for each $v\in\graph$.

A graph~$\graph$ is called
\roster*
\item
\emph{elliptic}, if $\Gs_+(\Z\graph)=\Gs_0(\Z\graph)=0$,
\item
\emph{parabolic}, if $\Gs_+(\Z\graph)=0$ and $\Gs_0(\Z\graph)>0$, and
\item
\emph{hyperbolic}, if $\Gs_+(\Z\graph)=1$.
\endroster
(Since we are interested in hyperbolic lattices only, we do not consider graphs
with $\Gs_+(\Z\graph)>1$.) When the degree $h^2=2D$ is fixed, we extend to
graphs the terms \emph{admissible}, \emph{geometric}, and \emph{subgeometric}
introduced in \autoref{def.geometric}, referring to the corresponding
properties of the lattice $\Fano_h(\graph)$.
By \autoref{th.K3}, subgeometric graphs are those that can appear as induced
subgraphs in the Fano graph $\Fn X$ of a smooth
$2D$-polarized $K3$-surface~$X$.
Conversely, a degree~$D$ is \emph{geometric} for~$\graph$ if
$\graph\cong\Fn X$ for some smooth $2D$-polarized $K3$-surface~$X$; the set
of geometric degrees is denoted by $\gd\graph$. Similarly, for a
subfield $K\subset\C$, we introduce the set $\gd_K\graph$ of the values
of~$D$ for which $\graph\cong\Fn_KX$ for some smooth $2D$-polarized
$K3$-surface $X\into\Cp{D+1}$ defined over~$K$.

The correspondence between graphs and lattices is not exactly one-to-one: in
general, $S\supset\Fano_h(\Fn S)$ is a finite index extension and
$\graph\subset\Fn\Fano_h(\graph)$ is an induced subgraph
(assuming that $\Fano_h(\graph)$ is hyperbolic).

If $\graph$ is elliptic, the lattice $\Fano_h(\graph)$ is
obviously hyperbolic. Parabolic graphs are treated separately in
\autoref{s.parabolic} below.
If $\graph$ is hyperbolic, we can define the
\emph{intrinsic polarization} as a vector $h_\graph\in\Z\graph\otimes\Q$
with the property
that $v\cdot h_\graph=1$ for each vertex $v\in\graph$.
If such a vector exists, it is unique modulo $\ker\Z\graph$; in particular,
$h_\graph^2\in\Q$ is well defined. (Note that $h_\graph^2$ may be negative.)

\proposition\label{prop.hyperbolic}
Given a hyperbolic graph~$\graph$, the $2D$-polarized lattice
$\Fano_h(\graph)$ is hyperbolic if and only if
$h_\graph$ exists and $2D\le h_\graph^2$.
\endproposition

\proof
Assume that $\Fano_h(\graph)$ is hyperbolic. Then, by~\eqref{eq.monotonous},
the image of the projection
$\Gf\:\Z\graph\to\Fano_h(\graph)$
is nondegenerate and
we have the orthogonal projection
\[*
\pr\:\Fano_h(\graph)\to\Gf(\Z\graph)\otimes\Q.
\]
Thus, $h_\graph$ exists (as any pull-back of $\pr h$) and, since
$\Gf(\Z\graph)$ itself is hyperbolic, its orthogonal complement must be
negative definite, \ie, $h^2\le h_\graph^2$
(and, if the equality holds, we must
have $h\in\Gf(\Z\graph)\otimes\Q$).

The converse statement is straightforward:
%Conversely,
if $h_\graph$ exists, we have an orthogonal direct sum
decomposition
$\Fano_h(\graph)\otimes\Q=(\Fano(\graph)\otimes\Q)\oplus\Q v$, where
$v:=h-\Gf(h_\graph)$, in which either $v=0$ (if $2D=h_\graph^2$) or $v^2<0$.
\endproof

\corollary[of the proof]\label{cor.D=h}
Let $\graph$ be a hyperbolic graph and $\rank\Fano(\graph)=20$. Then, the
$2D$-polarized lattice $\Fano_D(\graph)$ can be geometric only if
$2D=h_\graph^2$.
\done
\endcorollary

Thus, each hyperbolic graph~$\graph$ can be contained in the Fano graph of
a smooth $2D$-polarized $K3$-surface~$X$ only for finitely many values of~$D$.
Some of the values allowed by \autoref{prop.hyperbolic} can further
be eliminated by requiring that $\graph$ should be subgeometric and using
Algorithms~\ref{alg.lines} and~\ref{alg.Nikulin}. Below,
for each graph~$\graph$ used, we merely
%state the range of values of~$D$,
describe the set $\gd\graph$,
referring to \autoref{prop.hyperbolic} and omitting further
details of the computation.

\subsection{Valency bounds}\label{s.valency}
Applying
\autoref{prop.hyperbolic} to a star-shaped graph with a ``central''
vertex~$v$ of valency $\val v\ge5$,
we obtain the following
statement.

\proposition\label{prop.valency}
Let $X\subset\Cp{D+1}$ be a smooth $2D$-polarized $K3$-surface,
fix a line $v\in\Fn X$, and assume that all lines intersecting~$v$ are pairwise
disjoint.
Then, we have the following sharp bounds on~$\val v$\rom:
\[*
\minitab\quad
D=       & 2&3&4&5&6&7\dts11&\otherwise\\
\val v\le&12&9&8&7&6&5      &         4
\endminitab
\]
If $D=4$ and $X$ is a triquadric, then $\val v\le7$ \rom(sharp bound\rom).
%Then\rom:
%\roster*
%\item
%$\val v\le12$ if $D=2$,
%\item
%$\val v\le12-D$ if $3\le D\le7$,
%\item
%$\val v\le7$ if $X$ is a triquadric,
%\item
%$\val v\le5$ if $7\le D\le11$, and
%\item
%$\val v\le4$ if $D\ge12$.
%\endroster
%These bounds are sharp.
\done
\endproposition

In \autoref{prop.valency}, if $D=3$ and $v$ intersects
ten disjoint lines $a_1,\ldots.a_{10}$, then there is an eleventh line~$b$
intersecting~$v$ and one of~$a_i$, \cf. \autoref{prop.valency.trig} below.

Without the assumption that the star $\St v\subset\Fn X$ is discrete,
%In general,
the bound for quartics ($D=2$) is $\val v\le20$, and
%the star of~$v$ in $\Fn X$
$\St v$
may be rather complicated (see \cite{DIS,rams.schuett}).
The
%combinatorics of the
star simplifies dramatically if $D\ge3$: the following
statement is also proved by applying \autoref{prop.hyperbolic} to
a few simple test
configurations.

\proposition\label{prop.valency.trig}
Let $X\subset\Cp{D+1}$ be a smooth $2D$-polarized $K3$-surface, $D\ge3$,
and let $v,a_1,a_2\in\Fn X$ be three lines such that
$v\cdot a_1=v\cdot a_2=a_1\cdot a_2=1$.
Then, all other lines intersecting~$v$ are disjoint from
$a_1$, $a_2$, and each other and
we have the following sharp bounds on $\val v$\rom:
\[*
\minitab\quad
D=       & 3&4,5&6\dts11&\otherwise\\
\val v\le&11&  5&      3&         2
\endminitab
\]
Furthermore, if $D=3$, then $\val v\ne10$.
%\roster*
%\item
%$\val v\le9$ or $\val v=11$ if $D=3$,
%\item
%$\val v\le5$ if $D=4$, $5$,
%\item
%$\val v\le3$ if $6\le D\le11$, and
%\item
%$\val v\le2$ if $D\ge12$.
%\endroster
%These bounds are sharp.
\done
\endproposition

\section{Elliptic pencils}\label{S.elliptic}

In this section, we establish a relation between parabolic Fano graphs
$\Fn X$ and elliptic pencils on the $K3$-surface~$X$.
In particular, we establish a certain periodicity for such graphs
(see \autoref{prop.parabolic}) and prove Theorems~\ref{th.D>>0}
and~\ref{th.D>>0.real}.

\subsection{Parabolic graphs}\label{s.parabolic}
It is well known that any parabolic or elliptic
graph~$\graph$ is a disjoint union of
simply laced Dynkin diagrams and affine Dynkin diagrams, with at least one
component affine if $\graph$ is parabolic.
We describe the combinatorial type of~$\graph$ as a
formal sum of the corresponding $\bA$--$\bD$--$\bE$ types of its components.

%Let $\graph$ be a parabolic graph, and
%let $\fiber$ be an affine component of~$\graph$. Then,
If $\fiber$ is a connected parabolic graph (affine Dynkin diagram), the
kernel
$\ker\Z\fiber$ is generated by a single vector $\kk_\fiber=\sum n_vv$,
$v\in\fiber$; this generator can be chosen so that all $n_v>0$ and, under this
assumption, it is unique. We call the coefficient sum
$\deg\fiber:=\sum_vn_v$ the \emph{degree} of~$\fiber$. We have
\[
\gathered
\deg\tA_p=p+1,\ p\ge1,\qquad
\deg\tD_q=2q-2,\ q\ge4,\\
\deg\tE_6=12,\qquad
\deg\tE_7=18,\qquad
\deg\tE_8=30.
\endgathered
\label{eq.degree}
\]

For a connected elliptic graph (Dynkin diagram) $\fiber$, we can
define the \emph{degree set} $\ds\fiber$ as the set $\{\sum_vm_v\}$ of the
coefficient sums of all positive roots $\sum_vm_vv\in\Z\fiber$.
A simple computation shows that
\[
\gathered
\ds\bA_p=\{1\dts p\},\ p\ge1,\qquad
\ds\bD_q=\{1\dts2q-3\},\ q\ge4,\\
\ds\bE_6=\{1\dts11\},\qquad
\ds\bE_7=\{1\dts13,17\},\qquad
\ds\bE_8=\{1\dts16,23\}.
\endgathered
\label{eq.degrees}
\]
The \emph{Milnor number} of a parabolic or elliptic graph~$\graph$
is $\mu(\graph):=\rank\Fano(\graph)$.
We have
\[*
\ls|\graph|=\mu(\graph)+\#(\text{parabolic components of~$\graph$}).
\]

\lemma\label{lem.parabolic}
Let $\fiber$ be a parabolic component of an admissible parabolic
graph~$\graph$.
Then, the image of $\ker\Z\graph$ under the
projection $\Gf\:\Z\graph\to\Fano_h(\graph)$ is generated by the image
$\Gf(\kk_\fiber)\ne0$, which is primitive in any admissible extension
$S\supset\Fano_h(\graph)$.
\endlemma

\proof
Since $\Fano_h(\graph)$ is assumed hyperbolic, by~\eqref{eq.monotonous} one
has $\rank\Gf(\ker\Z\graph)\le1$. On the other hand,
$\Gf(\kk_\fiber)\cdot h=\deg\fiber\ne0$; hence, $\Gf(\kk_\fiber)\ne0$ and this
vector generates $\Gf(\ker\Z\graph)$ over~$\Q$. To show that
$\Gf(\kk_\fiber)\ne0$ is primitive, consider a maximal elliptic subgraph
$\fiber_0\subset\fiber$: it is a Dynkin diagram of the same name.
If $a:=\frac1mk_\fiber\in S$ for some $m\ge2$, then,
by~\eqref{eq.degree} and~\eqref{eq.degrees}, $a\cdot h\in\ds\fiber_0$; hence, there is a root
$e\in\Z\fiber_0$ with $e\cdot h=a\cdot h$, and $e-a$ is
an exceptional divisor as in
\autoref{def.geometric}\iref{bad.exceptional}.
\endproof

\autoref{lem.parabolic} implies that $\Gf$ maps elements $\kk_\fiber\in\Z\fiber$
corresponding to all parabolic components~$\fiber$ of~$\graph$ to the same element of
$\Fano_h(\graph)$, \viz. positive primitive generator of the kernel
$\ker\Gf(\Z\graph)$.

\corollary\label{cor.v.k}
In any admissible graph~$\graph$, one has
$v\cdot\kk_{\fiber'}=v\cdot\kk_{\fiber''}$ for any two disjoint connected
parabolic subgraphs $\fiber',\fiber''\subset\graph$ and any vector
$v\in\Fano_h(\graph)$.
\done
\endcorollary

\corollary\label{cor.degree}
All parabolic components of an admissible parabolic graph~$\graph$ have the
same degree\rom; this common degree $\deg\graph$ is called the \emph{degree}
of~$\graph$.
\done
\endcorollary

An admissible graph~$\graph$ is called \emph{saturated} if
$\graph=\Fn\Fano_h(\graph)$. (\latin{A priori}, this notion depends on the
choice of a degree $h^2=2D$.) The following statement is proved similar to
the primitivity of $\Gf(\kk_\fiber)$ in \autoref{lem.parabolic}.

\corollary\label{cor.elliptic}
If $\graph$ is an admissible parabolic graph and $\fiber$
is an elliptic component of~$\graph$, then
$\deg\graph\notin\ds\fiber$. If $\graph$ is saturated, then also
$(\deg\graph-1)\notin\ds\fiber$.
\done
\endcorollary

\subsection{Parabolic graphs as elliptic pencils}\label{s.elliptic}
Let $X$ be a smooth $2D$-polarized $K3$-surface, and consider a connected
parabolic subgraph $\fiber\subset\Fn X$. (Recall that $\fiber$ is an affine
Dynkin diagram.) The class $\kk_\fiber$, regarded as a divisor, is nef and,
since also $\kk_\fiber^2=0$ and $\kk_\fiber$ is primitive in $H_2(X;\Z)$
(see \autoref{lem.parabolic}), this class is a fiber of
a certain elliptic pencil $\pi\:X\to\Cp1$. We denote by
$\pencil(\fiber)\subset\Fn X$ the subgraph spanned by all linear components
of the reducible fibers of~$\pi$. Clearly,
\[*
\pencil(\fiber)=\bigl\{v\in\Fn X\bigm|v\cdot\kk_\fiber=0\bigr\};
\]
alternatively, $\pencil(\fiber)$ consists of~$\fiber$ and all vertices
$v\in\Fn X$ that are not adjacent to any of the vertices of~$\fiber$.
In this form, the notion of pencil can be extended to any geometric
configuration~$S$ and affine Dynkin diagram $\fiber\subset\Fn S$; we use the
same notation $\pencil(\fiber)$, or $\pencil_S(\fiber)$ if $S$ is to be
specified.

Note that $\pencil(\fiber)$ is a parabolic graph and
$\pencil(\fiber)=\pencil(\fiber')$ for any other parabolic component
$\fiber'\subset\pencil(\fiber)$. In this language, $\deg\pencil(\fiber)$ is
the common projective degree of the fibers of~$\pi$.

Thus, any maximal geometric parabolic graph~$\graph$
is the set of linear components
of the fibers of an elliptic pencil. By \autoref{lem.parabolic}
and~\eqref{eq.monotonous}, we have
\[
\mu(\graph)=\ls|\graph|-\#(\text{parabolic components of~$\graph$})\le18.
\label{eq.mu}
\]
Besides, the usual bound on the topological Euler characteristic yields
\[
\ls|\graph|\le\chi(X)=24,
\label{eq.Euler}
\]
with the inequality strict
whenever $\graph$ has at least one elliptic component.
Under the additional assumption that $\deg\graph\le7$, this can be refined to
\[
\mu(\graph)+b_0(\graph)
=\ls|\graph|+\#(\text{elliptic components of~$\graph$})\le24
\label{eq.Euler+}
\]
by taking into account the fiber components of higher projective degree.
(Indeed, considering singular elliptic fibers one by one, one can easily see that,
under the assumption $\deg\graph\le7$, the number of elliptic components
of~$\graph$ does
not exceed the number of fiber components of~$\pi$ that are not lines and,
thus, are not present in~$\graph$.
On the other hand, the total number of fiber components is at most~24.)
The two latter bounds are less obvious arithmetically; their proof would require
considering graphs satisfying~\eqref{eq.mu}
one by one and using \autoref{alg.Nikulin}
(\cf.~\cite{DIS}).

\remark\label{rem.pencils}
Without the assumption $\deg\graph\le7$, the \latin{a priori} possible
pencils~$\graph$ of a given degree~$d$ can be enumerated as follows.
Consider a collection $\fiber_1,\ldots,\fiber_n$ of affine Dynkin diagrams
such that $\sum_i\ls|\fiber_i|\le24$
% and $\deg\fiber_i\le d$ for each~$i$.
and assign an integral weight
(projective degree)
$w_v\ge1$ to each vertex~$v$ of the union so
that $\sum_{v\in\fiber}n_vw_v=d$ for each component $\fiber=\fiber_i$, where
$\sum_{v\in\fiber}n_vv=\kk_\fiber$.
%is the primitive positive generator of the kernel.
Then, take for~$\graph$ the induced subgraph of
the union $\bigcup_i\fiber_i$
spanned by all vertices of weight~$1$.
(Note that we have $\deg\fiber_i\le d$ for each~$i$ and, for $\graph$ to be
parabolic, the equality must hold for at least one value of~$i$.)
This description does \emph{not} guarantee that $\graph$ is geometric,
but it rules out many graphs that are not.
\endremark

To emphasize the relation between parabolic graphs and elliptic pencils, we
adopt the following terminology. A maximal parabolic subgraph
$\graph\subset\Fn X$ is called a \emph{pencil}, and its components
are called \emph{fibers};
a pencil is uniquely determined by any of its \emph{parabolic} fibers.
(Strictly speaking, only parabolic components of $\graph$ are
whole fibers of~$\pi$, whereas elliptic ones are parts of fibers
containing irreducible curves of higher projective degree; it may even happen
that several elliptic components of~$\graph$ are contained in the same fiber
of~$\pi$.) All other lines in~$X$ are called \emph{sections} of~$\graph$:
they are indeed (multi-)sections of~$\pi$, and we disregard all multisections
of higher projective degree. We denote by $\sec\graph:=\Fn X\sminus\graph$
the set of all sections, and
\[*
\sec_n\graph:=\bigl\{v\in\Fn X\bigm|
 \text{$v\cdot \kk_\fiber=n$ for some/any parabolic fiber~$\fiber$ of~$\graph$}\bigr\}
\]
stands for the set of $n$-sections, $n\ge1$ (see \autoref{cor.v.k}).
For $l\in\graph$, we also use
\[*
\sec l:=\bigl\{v\in\Fn X\sminus\graph\bigm|v\cdot l=1\}
\quad\text{and}\quad
\sec_nl:=\sec l\cap\sec_n\graph.
\]
Sometimes, $1$- and $2$-sections are called \emph{simple sections}
and \emph{bisections}, respectively.
Any section of~$\graph$ intersects each parabolic fiber, but it may miss
elliptic ones.

As an immediate consequence of these definitions, we have
\[
\ls|\Fn X|=\ls|\graph|+\ls|\sec\graph|
\label{eq.count}
\]
for each pencil $\graph\subset\Fn X$.

\convention\label{conv.coord}
A section $s\in\sec\graph$ can be described by its \emph{coordinates}
%$\bs\subset\graph$, \viz. the subset
\[*
\bs:=\bigl\{l\in\graph\bigm|s\cdot l=1\bigr\}\subset\graph.
\]
If $s$ is an $n$-section, then
$\ls|\bs\cap\fiber'|=n$ and $\ls|\bs\cap\fiber''|\le n$ for each parabolic
fiber $\fiber'\subset\graph$ and elliptic fiber $\fiber''\subset\graph$,
respectively.
Often, the map $s\mapsto\bs$ is injective, but we do not assert this in
general.
If a configuration~$S$ is spanned by~$h$, $\graph$, and a number of
\emph{pairwise disjoint} sections $s_1,\ldots,s_n\in\sec\graph$, the Gram
matrix of~$S$ is determined by the coordinates $\bs_1,\ldots,\bs_n$. That is
why, in the rest of the paper,
we pay special attention to finding as many \latin{a priori}
disjoint sections as possible.
\endconvention

\subsection{The periodicity}\label{s.periodicity}
Denote by $\pi_D(\graph)$ the set of connected components of the
equilinear stratum $\Fn X\cong\graph$ of smooth $2D$-polarized $K3$-surfaces.
%(\cf. \autoref{th.K3}).
The following statement
is in sharp contrast with \autoref{prop.hyperbolic}.

\proposition\label{prop.parabolic}
Let $\graph$ be a parabolic graph, $d:=\deg\graph$, and let $D\ge d^2+d$.
Then, for any $D'=D\bmod d$, there is a canonical inclusion
$\pi_{D'}(\graph)\into\pi_D(\graph)$. This inclusion is a bijection if
also $D'\ge d^2+d$.
\endproposition

\proof
Consider the sublattice $H:=\Z h+\Z \kk\subset\Fano_h(\graph)$, where $\kk$ is
the common image of the vectors $\kk_\fiber\in\Z\fiber$
corresponding to all parabolic
components~$\fiber$ of~$\graph$. The change of variables $h\mapsto h\pm \kk$ shows that
the abstract isomorphism type of the pair of lattices
$H\subset\Fano_h(\graph)$ depends on $D\bmod d$ only.
Since any polarized automorphism of $\Fano_h(\graph)$ preserves~$H$
pointwise, the set of isomorphism classes of embeddings
$\Fano_h(\graph)\into\L$ also depends only on $D\bmod d$;
hence, by \autoref{th.K3}, we only need to investigate the dependence on~$D$
%(within the same residue)
of the admissibility of a finite index extension
$S\supset\Fano_h(\graph)$ and the existence of extra lines.
%, \ie, the set
%$\Fn S\sminus\graph$.

Let $E:=H^\perp$; this is a negative definite lattice.
We have $S\subset\Fano_h(\graph)\dual\subset H\dual\oplus E\dual$
(\cf.~\cite{Nikulin:forms}), and any vector $e\in S$ decomposes as
$e=e_H+e_E$, $e_H\in H\dual$, $e_E\in E\dual$.
The sublattice $H$ is primitive: if $a:=\Ga h+\Gb\kk\in\L$, we have
$\Ga=a\cdot l\in\Z$ for any $l\in\graph$, and then $\Gb\in\Z$ by
\autoref{lem.parabolic}.
%Note that $H$ is primitive
%in~$\L$ (\cf. \autoref{lem.parabolic}); hence,
Hence, $e_E^2<0$ unless $e\in H$.
Note also that $H\dual$ is generated by $h^*=\kk/d$
and $\kk^*=(dh-2D\kk)/d^2$.

Let $D\ge d^2+d$.
If $e$ is an exceptional divisor as in
\autoref{def.geometric}\iref{bad.exceptional}, then
\[*
e_H=m\kk^*,\ m\in\Z,\qquad e_H^2=-\frac{2m^2D}{d^2}<-2
\]
unless $m=0$. It follows that $e\in E\dual$ is independent of~$D$.
If $e$ is a quadric pencil as in
\autoref{def.geometric}\iref{bad.hyperelliptic}, then
\[*
e_H=2h^*+m\kk^*,\ m\in\Z,\qquad e_H^2=-\frac{2m(mD-2d)}{d^2}\le0.
\]
Hence, $m=0$ and $e=e_H=2h^*$, but this vector cannot be in~$S$ by the
primitivity of~$H$.
Finally, if $e$ is a line, then
\[*
e_H=h^*+m\kk^*,\ m\in\Z,\qquad e_H^2=-\frac{2m(mD-d)}{d^2}<-2
\]
unless $m=0$ or $m=1$ and $D=d^2+d$. In the latter case, $e_H^2=-2$ and,
hence, $e=e_H=h^*+\kk^*\notin S$. In the former case, $e=h^*+e_E$, where
$e_E\in E\dual$, $e_E^2=-2$ does not depend on~$D$.
Thus, oll ``undesirable'' vectors are independent of the
choice of
$D\ge d^2+d$ within the same class $D\bmod d$. If $D<d^2+d$, these
vectors are still present, but new ones can appear,
ruling out some of the configurations.
%; hence, new
%configurations do not appear, but some of the old ones may be ruled out.
\endproof

\corollary\label{cor.parabolic}
%Listed in
\autoref{tab.parabolic} list all subgeometric pencils~$\graph$ such
that $\ls|\graph|\ge21$.
\table
\def\-{\rlap{$\mdag$}}
\def\*{\rlap{$\mstar$}}
\def\r{\rlap{$\mreal$}}
\def\={\relax\afterassignment\ul\count0=}
\def\ul{\underline{\the\count0}}
\rm
\caption{Large subgeometric parabolic graphs}\label{tab.parabolic}
\def\+{\mplus}
\hbox to\hsize{\hss\vbox{\halign{\strut\quad\hss$#$\hss&
 \quad$#$\hss\quad&#\hss&\quad$#$\hss\quad\cr
\noalign{\hrule\vspace{2pt}}
\ls|\graph|&\graph&Polarizations $D$&T:=\Fano_h(\graph)^\perp\cr
\noalign{\vspace{1pt}\hrule\vspace{2pt}}
24&8\tA_2       &$4\+$ or $\{7\bmod3\}$    &\bU^2(3)\cr
  &6\tA_3       &$2\+$ or $\{6\bmod4\}$    &[8,0,8]\cr
22&7\tA_2+\bA_1 &$4\+$ or $\{7\bmod3\}$    &\bU^2(3)\oplus[-6]\cr
  &5\tA_3+\bA_2 &$2\+$ or $\{6\bmod4\}$    &\bA_2(-4)\oplus[-16]\cr
  &5\tA_3+2\bA_1&$2\+$ or $\{4\bmod2\}$    &[8,0,8]\oplus[-4]\cr
  &4\tA_4+2\bA_1&$2\+$ or $\{7\bmod5\}$    &[10,0,10]\cr
21&7\tA_2       &$\{4\bmod3\}$             &\bU^2(3)\oplus\bA_2(-3)\cr
  &             &$5\+$ or $\{8\bmod3\}$    &\bU^2(3)\oplus\bA_2(-1)\cr
  &6\tA_2+3\bA_1&$3\+$ or $\{6\bmod3\}$\r  &\bA_2\oplus[-6]^3\cr
  &             &$\{4\bmod3\}$             &\bA_2(3)\oplus[-6]^3\cr
  &5\tA_3+\bA_1 &$\{2\bmod2\}$             &[12,4,12]\oplus[-4]^2\cr
  &4\tA_4+\bA_1 &$\{2,5\bmod5\}$           &[50]\oplus\bU(5)\cr
  &             &$6\+$ or $\{11\bmod5\}$\r &[2]\oplus\bU(5)\cr
  &3\tA_5+\bA_3 &$\{4\bmod6\}$             &[24,0,36]\cr
  &             &$2\+$ or $\{8\bmod6\}$    &[4,0,24]\cr
  &3\tA_5+3\bA_1&$\{4\bmod6\}$             &[24,12,24]\cr
  &             &$\{6\bmod6\}$             &[8,4,8]\cr
  &3\tA_6       &$14\+$ or $\{21\bmod7\}$\r&[2,1,4]\cr
  &             &$\{2,4,8\bmod7\}$         &[14,7,28]\cr
  &2\tD_5+\tA_7+\bA_1&$\{2\bmod4\}$        &[32,16,40]\cr
  &3\tD_6       &$\{10\bmod10\}$           &[4,0,4],[20,0,20]\cr
  &             &$\{2,8\bmod10\}$          &[4,0,100]\cr
  &             &$\{4,6\bmod10\}$          &[8,4,52]\cr
  &3\tE_6       &$\{4\bmod12\}$            &[6,0,72]\cr
  &             &$\{6\bmod12\}$\r          &[2,0,24]\cr
  &             &$\{7\bmod12\}$            &[6,0,18],[6,3,6]\cr
  &             &$\{10\bmod12\}$           &[18,0,24]\cr
  &             &$\{12\bmod12\}$           &[6,0,8]\cr
  &             &$\{15\bmod12\}$\r         &[2,0,6],[2,1,2]\cr
\noalign{\vspace{1pt}\hrule}
\crcr}}\hss}
\endtable
\endcorollary

In the table, for each pencil~$\graph$,
we list the values of~$D$ (marking with a~$^+$
those for which any geometric finite index extension of
$\Fano_h(\graph)$ contains sections of~$\graph$) and generic
transcendental lattices $T:=\Fano_h(\graph)^\perp$,
which depend on $D\bmod\deg\graph$ only.
%If $T$ contains $[2]$ or $\bU(2)$, the corresponding series is marked with
%a~$\mreal$; according to \autoref{th.real}, the pencil can be chosen to
%contain real lines only in this case.

\proof[Proof of \autoref{cor.parabolic}]
If $\ls|\graph|=24$, then, by~\eqref{eq.mu}, $\graph$ has at least six
parabolic fibers and, hence, the minimal Milnor number $\mu_0$ of such a
fiber is subject to the inequality $6\mu_0\le18$.
Taking~\eqref{eq.Euler+} into account, we arrive at the two combinatorial
types listed in the table. Similarly, one can show that $\ls|\graph|\ne23$
and classify all pencils with $\ls|\graph|=22$ or~$21$.
(The combinatorial types that are never geometric
are not shown in the table.) Finally, the set $\gd\graph$
for each combinatorial type~$\graph$
is determined using
\autoref{prop.parabolic} and \autoref{alg.Nikulin}.
\endproof

\remark\label{rem.parabolic.vs.hyperbolic}
The bound $d^2+d$ in \autoref{prop.parabolic} is closely related to the
intrinsic polarization in \autoref{prop.hyperbolic}.
If $\graph$ consists of an affine Dynkin diagram of degree~$d$ and a
single $n$-section, one
%obviously
has $h_\graph^2<2(d/n)^2+2(d/n)$.
By the obvious monotonicity of $h_\graph^2$,
%It follows that,
if $D\ge d^2+d$, any (sub-)geometric graph containing an affine Dynkin
diagram of degree~$d$ is parabolic.
(This conclusion can as well be derived from the proof of
\autoref{prop.parabolic}.)
\endremark

\subsection{Proof of \autoref{th.D>>0}}\label{proof.D>>0}
If the Fano graph
$\Fn X$ is hyperbolic, it contains a pencil~$\graph$.
There are but finitely many parabolic graphs satisfying~\eqref{eq.mu}
and~\eqref{eq.Euler}, each of them has finitely many hyperbolic $1$-vertex
extensions, and each extension is geometric for finitely many values of~$D$
by \autoref{prop.hyperbolic}. Thus, for $D\gg0$, the graph $\Fn X$ is either
elliptic, and then $\ls|\Fn X|\le19$ by~\eqref{eq.monotonous}, or parabolic.
In the latter case, the bound $\mc(D)$ on $\ls|\Fn X|$ is given by
\autoref{cor.parabolic}.
\qed

\remark\label{rem.min.parabolic}
The precise lower bound on~$D$ in \autoref{th.D>>0} is rather high: the
conclusion holds to full extent for $D\ge621$.
There exists a
one-vertex
extension of the affine Dynkin diagram~$\tE_8$
(\viz. $\tE_8$ itself and one simple section, \cf.
\autoref{fig.sections} in \autoref{s.taxonomy} below)
which is geometric for all $D\le620$.
\endremark

\remark
\autoref{cor.parabolic} shows also that, for each $D\ge6$, there
exists a smooth $2D$-polarized $K3$-surface~$X$ with $\Fn X$ parabolic and
$\ls|\Fn X|=21$: namely, one can have $\Fn X\cong7\tA_2$ or $6\tA_2+3\bA_1$.
\endremark

\subsection{Proof of \autoref{th.D>>0.real}}\label{proof.D>>0.real}
In view of \autoref{th.D>>0},
the bound $\ls|\Fn_\R X|\le21$ and the values of~$D$ for which it is sharp
follow
immediately from
\autoref{th.real} and \autoref{cor.elliptic} (see \autoref{tab.parabolic},
where real configurations are marked).
For the lower bound $\mr(D)\ge19$ one can use a pencil
$\graph\cong6\tA_2+\bA_1$ or $2\tE_8+\bA_1$.
The values of $D\gg0$ for which $\mr(D)=20$ are given by the classification
of pencils of size~$20$ (using \autoref{rem.pencils}):
the complete list, which is too long, is found in~\cite{degt:Fano.graphs}.
\qed

\section{The taxonomy of hyperbolic graphs}\label{S.locally.elliptic}

In~\cite{DIS} and~\cite{degt:supersingular}, the Fano graphs of spatial
quartics are subdivided into triangular (all large graphs), quadrangular, and
quadrangle free.
In this paper, we adopt a more consistent taxonomy based on the type
of a minimal fiber (see \autoref{def.min.fiber}): this approach
allows a more refined classification and stronger bounds.
%In this section, we
%begin the type-by-type study of the graphs, starting with the simplest
%\emph{locally elliptic} ones.

Till the rest of the paper, we consider hyperbolic graphs only, even though
most statements do not make this assumption.
In the proofs, we implicitly refer to \autoref{cor.parabolic} or, if
necessary, more refined classification of the parabolic graphs of a given
type based directly on \autoref{prop.parabolic} (see
also~\cite{degt:Fano.graphs}).

\subsection{The taxonomy}\label{s.taxonomy}
Recall that any hyperbolic graph~$\graph$ contains a
connected parabolic subgraph: this
obvious statement is part of the classification of elliptic and parabolic
graphs.

\definition\label{def.min.fiber}
A \emph{minimal fiber} in a hyperbolic or parabolic graph~$\graph$ is a
connected parabolic subgraph (affine Dynkin
diagram) $\fiber\subset\graph$ of the minimal possible Milnor number, with
$\tA$ preferred over~$\tD$ over~$\tE$ in the case of equal Milnor numbers.
The graph~$\graph$ itself is referred to as a \emph{$\fiber$-graph}, where
$\fiber$ is the (obviously well-defined) common isomorphism type of the
minimal fibers in~$\graph$.
A \emph{minimal pencil} in~$\graph$ is any pencil of the form
$\pencil(\fiber)\subset\graph$, where $\fiber\subset\graph$ is a minimal
fiber.
This terminology also applies to configurations $(S,h)$, according to the
Fano graph $\Fn(S,h)$.
\enddefinition

The $\tA_2$-, $\tA_3$-, and $\tA_4$-graphs/configurations are alternatively
called \emph{triangular},
\emph{quadrangular}, and \emph{pentagonal}, respectively; these
graphs~$\graph$
are
characterized by the property that $\girth(\graph)=3$, $4$, and~$5$,
respectively.

The $\tD_4$-graphs/configurations are called \emph{astral}; such graphs~$\graph$
are characterized by
the property that $\girth(\graph)\ge6$ and $\graph$ has a vertex of valency
at least four.
Since $\tD_4$ is the only affine Dynkin diagram with a vertex of
valency~$4$,
any $\fiber$-graph~$\graph$ with $\mu(\fiber)\ge5$ is \emph{locally
elliptic}, \ie, one has $\val v\le3$ for any vertex $v\in\graph$.
Analyzing such affine Dynkin diagrams~$\fiber$ one by one, we easily arrive at the
following bound on the number and positions of the sections of a minimal
pencil.

\proposition\label{prop.le.sections}
Let $\fiber\subset\graph$ be a minimal fiber,
and assume that $\mu(\fiber)\ge5$. Then, the
pencil $\pencil(\fiber)$ has at most three parabolic fibers and its
sections are as shown in
\autoref{fig.sections}. Furthermore, all sections are simple and, unless
\smash{$\fiber\cong\tA_5$}, pairwise disjoint.
If \smash{$\fiber\cong\tA_5$},
%a pair of
two sections may intersect only if they are adjacent to opposite corners of the
hexagon~$\fiber$.
\figure
\cpic{le}
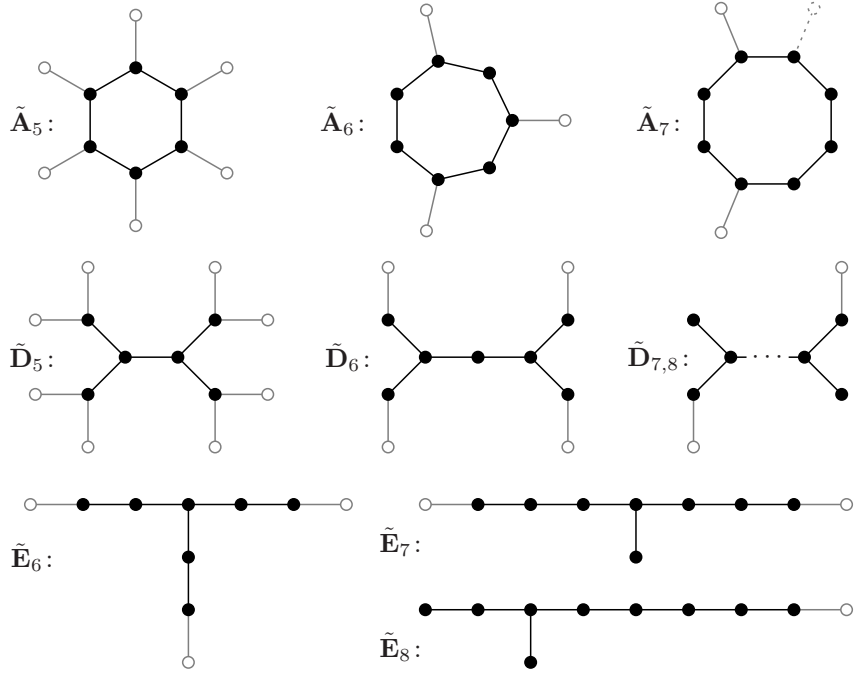
\caption{Sections (white) at a minimal fiber (black)
 (see \autoref{prop.le.sections})}\label{fig.sections}
\endfigure
\done
\endproposition

\remark\label{rem.le.sections}
Certainly, \autoref{fig.sections} shows maximal sets of sections. If
\smash{$\fiber\cong\tA_7$}, the condition is that the ``distance'' between
the two sections in the octagon~$\fiber$ must be at least~$3$
(as otherwise the graph would contain \smash{$\tD_5$} or \smash{$\tD_6$});
hence, there are two distinct maximal sets.
\endremark

%\begingroup
%\setbox22\hbox{$\minitab\quad
%D:  &2\dts10,12\dts16,18\dts28&11,17,29\dts34,36\dts38,40,50\\
%\bm:&                       22&                           21\\
%\endminitab$}
%\setbox36\hbox{$\minitab\quad
%D:  &2\dts38,40\dts45&39&48\dts50,64\\
%\bm:&              22&24&         21\\
%\endminitab$}
%\begin{gather}
%\copy22\tag{$\tD_6$}\\
%\copy36\tag{$\tE_6$}
%\end{gather}
%\endgroup

\subsection{Locally elliptic graphs}\label{s.locally.elliptic}
In view of Theorems~\ref{th.D>>0} and~\ref{th.D>>0.real}, most
interesting are the configurations with twenty or more lines. In the locally
elliptic case,
%(in the sense of the previous section, \ie, $\fiber$-graphs
%with $\mu(\fiber)\ge5$),
one can easily list such configurations. (If $D$ is small, the
number of configurations may be huge and we confine ourselves to those close
to maximal.) The results are collected in the several statements below. In
\autoref{th.mu>5}, we skip the values of~$D$ for which $\bm(D)=20$, as they
are too erratic; these values, as well as a complete description of all
graphs, are found in~\cite{degt:Fano.graphs}.

\theorem[see \autoref{proof.le}]\label{th.a5}
Let $X$ be a smooth $2D$-polarized $K3$-surface, and assume that $\Fn X$ is a
hyperbolic $\tA_5$-graph. Then, with one exception $\Fn X\cong\LE_{25}$ shown
in \autoref{tab.le}, one has a sharp bound $\ls|\Fn X|\le\bm$, where
$\bm:=\bm(D)$ is as follows\rom:
\table
\rm
\def\-{\rlap{$\mdag$}}
\def\*{\rlap{$\mstar$}}
\def\r{\rlap{$\mreal$}}
\def\={\relax\afterassignment\ul\count0=}
\def\ul{\underline{\the\count0}}
\caption{Large locally elliptic graphs (see Theorems~\ref{th.a5}
 and~\ref{th.d5})}\label{tab.le}
\hbox to\hsize{\hss\vbox{\halign{\strut\quad\hss$#$\hss\quad&
 \hss$#$\hss\quad&\hss$#$\hss\quad&\hss$#$\hss\quad&\hss$#$\hss\quad&
 \hss$#$\hss\quad&$#$\hss\quad\cr
\noalign{\hrule\vspace{2pt}}
D&\graph&\ls|\Aut\graph|&\det&(r,c)&\ls|\Aut X|&T:=\Fano_h(\graph)^\perp\cr
\noalign{\vspace{1pt}\hrule\vspace{2pt}}
 2  &\LE_{25}   & 144&780&(2,0)&  12&[4,2,196]\cr
 4  &\LE_{25}   & 144&620&(2,0)&  12&[4,2,156],\ [20,10,36]\cr
 6  &\LE_{25}   & 144&460&(1,0)&  12&[20,10,28]\cr
10\-&\LE_{25}   & 144&140&(1,0)&  12&[12,2,12]\cr
14  &\LE_{28}\r & 336&\=7&(1,0)& 336&[2,1,4]\cr
\noalign{\vspace{1pt}\hrule}
\crcr}}\hss}
\endtable
\[*
\minitab\quad
D:  &2\dts12&13&14&15&16\dts20,22&\otherwise \\
\bm:&     24&22&23&21&         20&\twenty{19}\\
\endminitab
\]
\endtheorem

The notation in \autoref{tab.le} is similar to \autoref{tab.main}
(see \autoref{s.results}), with most fields skipped.

\theorem[see \autoref{proof.le}]\label{th.d5}
Let $X$ be a smooth $2D$-polarized $K3$-surface, and assume that $\Fn X$ is a
hyperbolic $\tD_5$-graph. Then, with one exception $\Fn X\cong\LE_{28}$ shown
in \autoref{tab.le}, one has a sharp bound $\ls|\Fn X|\le\bm$, where
$\bm:=\bm(D)$ is as follows\rom:
\[*
\minitab\quad
D:  &2\dts20&21,26,32&22,24&23,25,27,28,30,34,36&\otherwise \\
\bm:&     24&      21&   22&                  20&\twenty{19}\\
\endminitab
\]
\endtheorem

The extremal graph $\LE_{28}$ is $3$-regular, $3$-arc transitive, and
distance regular; it is
isomorphic to the so-called
\emph{Coxeter graph} (the only symmetric $3$-regular graph on $28$ vertices).
The
configuration
has $21$ minimal pencils, all of type $2\tD_5+\tA_7$.

\theorem[see \autoref{proof.le}]\label{th.mu>5}
Let $X$ be a smooth $2D$-polarized $K3$-surface, and assume that $\Fn X$ is a
hyperbolic
$\fiber$-graph, $\mu(\fiber)\ge6$. Then $\ls|\Fn X|\le\bm:=\bm(D)$, where
\roster*
\item
if \smash{$\fiber\cong\tE_8$} and $D=309$, then $\bm=21$\rom;
\item
if \smash{$\fiber\cong\tE_7$} and $D\in\{2\dts99\}$, then $\bm=21$\rom;
\item
if \smash{$\fiber\cong\tA_6$} and $D\in\{18,21,22\}$, then $\bm=21$\rom;
\item
if \smash{$\fiber\cong\tD_6$}, then $\bm$ is as follows\rom:
\[*
\minitab\quad
D:  &2\dts10,12\dts16,18\dts28&11,17,29\dts34,36\dts38,40,50\\
\bm:&                       22&                           21\\
\endminitab
\]
\item
if \smash{$\fiber\cong\tE_6$}, then $\bm$ is as follows\rom:
\[*
\minitab\quad
D:  &2\dts38,40\dts45&39&48\dts50,64\\
\bm:&              22&24&         21\\
\endminitab
\]
\item
in all other cases, $\bm\le20$.
\endroster
These bounds are sharp.
\rom(See~\cite{degt:Fano.graphs} for the graphs and the values
$\bm(D)=20$.\rom)
\endtheorem

\addendum[see \autoref{proof.le}]\label{ad.le}
There are but six hyperbolic locally elliptic Fano graphs
with $24$ vertices\rom:
\roster*
\item
$\LE_{24}\A$: $\ls|\Aut\LE_{24}\A|=144$,
$3$-regular,
$\gd=\{2\dts12\}$, $\gd_\R=\{3\dts6,8\dts12\}$\rom;
%\rom($3$-regular\rom)\rom;
%bipartite $3$-regular\rom;
\item
$\LE_{24}\E$:
$\ls|\Aut\LE_{24}\E|=72$, $\gd=\gd_\R=\{39\}$\rom;
%bipartite\rom;
\item
$\LE_{24}$: $\ls|\Aut\LE_{24}|=48$,
$\gd=\{2\dts20\}$, $\gd_\R=\{5,10,13,18\}$\rom;
\item
$\LE_{24}'$: $\ls|\Aut\LE_{24}'|=32$,
$\gd=\{2\dts16\}$, $\gd_\R=\{7,12,15,16\}$\rom;
%bipartite\rom;
\item
$\LE_{24}''$: $\ls|\Aut\LE_{24}''|=12$,
$\gd=\{2\dts13\}$, $\gd_\R=\{5,7,12\}$\rom;
\item
$\LE_{24}'''$: $\ls|\Aut\LE_{24}'''|=8$,
$\gd=\{14\}$, $\gd_\R=\varnothing$.
\endroster
As a consequence, for each $D\in\{2\dts20\}$, there exists a
geometric hyperbolic
locally elliptic configuration with $24$ lines, and for each
$D\in\{3\dts13,15,16\}$, there exists a real configuration with these
properties
\rom(\cf. the sharpness in Theorems \ref{th.main}
and~\ref{th.main.real}\rom).
\endaddendum

In the subsequent sections we show that, if $D>16$ and $\mu(\fiber)\le4$, the
number of lines in geometric $\fiber$-configurations $(S,h)$, $h^2=2D$, is
maximized by parabolic graphs. Combining this observation with
Theorems~\ref{th.a5}--\ref{th.mu>5}, we obtain precise sharp bounds on
the number of lines.

\corollary\label{cor.le.mc}
If $D>20$ \rom(\cf. \autoref{th.main}\rom), the maximum $\mc(D)$ is as given
by \autoref{th.D>>0}, with the exception of the following eight values\rom:
$\mc(39)=24$
and $\mc(D)=22$ for $D\in\{ 21, 23, 29, 33, 35, 41, 45\}$.
\done
\endcorollary

\corollary\label{cor.le.mr}
If $D>16$ \rom(\cf. \autoref{th.main.real}\rom),
one has
\[*
\mr(D)=\max\bigl\{\mpar(D),\mhyp(D)\bigl\},
\]
where $\mpar(D)$ is the maximum over the parabolic graphs given by
\autoref{th.D>>0.real}, and $\mhyp(D)$ is the maximum over the hyperbolic
graphs, which is as follows\rom:
%the maximum $\mc(D)$ is as given
%by \autoref{th.D>>0}, with the exception of the following eight values\rom:
%$\mc(39)=24$
%and $\mc(D)=22$ for $D\in\{ 21, 23, 29, 33, 35, 41, 45\}$.
\roster*
\item
$\mhyp(D)=24$ if $D\in\{18,39\}$,
\item
$\mhyp(D)=22$ if
$D\in\{17, 19, 21..25, 27, 30, 33, 34, 37, 41, 42, 45\}$,
\item
$\mhyp(D)=21$ if
$D\in\{20, 26, 28, 29, 31, 32, 35, 36, 38, 40, 43, 44, 46\dts99, 309\}$,
\item
$\mhyp(D)=20$ if
$D\in\{100\dts308, 342, 396, 565, 576, 601\}$,
\endroster
and $\mhyp(D)\le19$ for all other values $D>16$.
\done
\endcorollary

\subsection{Idea of the proof}\label{s.idea}
The approach outlined in this section will be used in most proofs till the
end of the paper.

We attempt to prove a bound $\ls|\Fn X|\le \Mm$ for graphs
$\Fn X$ satisfying certain conditions, subject to certain exceptions that are
to be determined. To this end, we choose an
appropriate pencil $\graph:=\pencil(\fiber)$
with a distinguished fiber~$\fiber$ and try to estimate separately the two
terms in~\eqref{eq.count}. Typically, the bound given by the sum of the two
estimates is too rough. Therefore, we choose appropriate thresholds $\Mg$
and $\Ms$, so that $\Mg+\Ms=\Mm$, and try to \emph{classify},
separately, the following two types of geometric graphs:
\roster*
\item
pencils~$\graph\supset\fiber$
%with more than $\Mg$ lines and
of size $\ls|\graph|>\Mg$ with
at least a few sections
without which the count $\ls|\Fn X|>\Mm$ is not reachable
(see \autoref{s.alg.pencils}), and
\item
section sets
$\sec\graph$ (in the presence of just the distinguished fiber~$\fiber$) of
size $\ls|\sec\graph|>\Ms$ (see \autoref{s.alg.sections}).
\endroster
The thresholds are chosen so that, on the
one hand, the computation is feasible and, on the other hand, the
partial configurations~$S$ obtained have maximal or close to maximal rank. In
the former case,
$\rank S=20$, all geometric extensions of~$S$ are found by
\autoref{alg.Nikulin}; in the latter case, we either add a few
(usually, up to two)
more lines to
increase the rank to~$20$ or prove that $S$ has no further extensions.

The graphs
$\Fn X$ with $\ls|\Fn X|>\Mm$
are found
in the process of the classification, and all others satisfy the
inequality $\ls|\Fn X|\le \Mg+\Ms=\Mm$ given by~\eqref{eq.count}.

\subsection{Proof of Theorems~\ref{th.a5}--\ref{th.mu>5} and \autoref{ad.le}}\label{proof.le}
Given a configuration~$S$ as in the statements, fix a minimal pencil
$\graph:=\pencil(\fiber)\subset\Fn S$.
By \autoref{prop.le.sections},
this pencil has at most three parabolic fibers, hence $\ls|\graph|\le21$
by~\eqref{eq.mu}, and at most eight sections; thus, by~\eqref{eq.count},
we have a bound
$\ls|\Fn S|\le29$, which holds in any degree and over any field (including
positive characteristic).

For the more precise bounds as in the statements, we
classify the configurations, running the modified
version (see \autoref{s.alg.le}) of the algorithm of \autoref{s.alg.pencils}.
For \autoref{th.mu>5}, we take $\Mm=19$, obtaining all graphs with $20$ or
more vertices. For Theorems~\ref{th.a5}, \ref{th.d5} and \autoref{ad.le}, in
order to avoid too large output, we choose $19\le\Mm\le23$
close to the expected maximum, depending on~$D$.
The other thresholds used
are $\Mg=\Mm-\Ms$ and $\Ms=\max{\ls|\sec\graph|}$, where
the latter maximum, depending on~$\fiber$ and~$D$, is computed using
\autoref{prop.hyperbolic}.
(If \smash{$\fiber\cong\bA_5$}, we need to take into
account the possible intersections of pairs of sections adjacent to opposite
corner; up to automorphism, this results in $23$ sets of sections.)
A few technical observations (in terms of the coordinates, see
\autoref{conv.coord}) reduce the number of candidates:
\roster*
\item
$\ls|\bs|\le3$ for each $s\in\sec\graph$
(as the graph is locally elliptic);
\item
$\ls|\bs'\cap\bs''|\le1$ for $s'\ne s''\in\sec\graph$
(as the graph is quadrangle free);
\item
$\ls|\bs|\le2$ for each $s\in\sec\graph$, if $\fiber\cong\tD_7$
(no subgraphs isomorphic to $\tD_6$);
\item
$\ls|\bs|\le1$ for each $s\in\sec\graph$, if $\fiber\cong\tD_8$
(no subgraphs isomorphic to $\tE_7$).
\endroster
Finally, the defining property~$\prop$ is that
$S$ should be a $\fiber$-configuration: since we always have
$\fiber\subset\Fn S$, this property is obviously hereditary.
\qed

\section{Other pencils with disjoint sections}\label{S.astral}

There are two other (unrelated) classes of pencils in Fano graphs whose
sections are
\latin{a priori} known to be disjoint: astral configurations and Fano graphs of
special octics.
At the end of this section,
in \autoref{s.hyperelliptic},
we consider also the toy problem of counting lines in
hyperelliptic models of degree $2D\ge6$.

\subsection{Astral configurations}
Recall that \emph{astral} are the $\tD_4$-graphs. These graphs~$\graph$ are
characterized by the property that $\girth(\graph)\ge6$ and
$\max\val v\ge4$, $v\in\graph$.

\theorem[see \autoref{proof.astral}]\label{th.astral}
Let $X$ be a smooth $2D$-polarized $K3$-surface,
and assume that
the Fano graph
$\Fn X$ is astral.
\table\rm
\def\-{\rlap{$\mdag$}}
\def\*{\rlap{$\mstar$}}
\def\r{\rlap{$\mreal$}}
\def\={\relax\afterassignment\ul\count0=}
\def\ul{\underline{\the\count0}}
\caption{Extremal astral configurations (see \autoref{th.astral})}\label{tab.SC}
\def\*{\mstar}
\hbox to\hsize{\hss\vbox{\halign{\strut\quad\hss$#$\hss\quad&
 \hss$#$\hss\quad&\hss$#$\hss\quad&\hss$#$\hss\quad&\hss$#$\hss\quad&
 \hss$#$\hss\quad&\hss$#$\hss\quad&$#$\hss\quad\cr
\noalign{\hrule\vspace{2pt}}
D&\bm&\graph&\ls|\Aut\graph|&\det&(r,c)&\ls|\Aut X|&T:=\Fano_h(\graph)^\perp\cr
\noalign{\vspace{1pt}\hrule\vspace{2pt}}
 2  &26&\SC_{27}'   &  72& 300&(2,0)&  12&[4,2,76]\cr
 4  &26&\SC_{28}'   & 576&  44&(1,0)&  96&[4,2,12]\cr
      &&\SC_{27}'   &  72& 204&(0,1)&  12&[12,6,20]\cr
      &&\SC_{27}''  &  72& 207&(1,1)&  18&[6,3,36],[12,3,18]\*\cr
 6  &26&\SC_{28}'   & 576&  28&(1,0)&  96&[4,2,8]\cr
 7  &24&\SC_{27}''  &  72&  99&(1,0)&  18&[6,3,18]\cr
 8  &24&\SC_{32}'   &2304&\=12&(1,0)&1152&[4,2,4]\cr
      &&\SC_{26}    &  24&  60&(1,0)&   6&[4,2,16]\cr
      &&\SC_{25}''  &  16&  60&(1,0)&   4&[4,2,16]\cr
10\-&24&\SC_{25}'   &  24&  96&(0,1)&   6&[4,0,24]\cr
12  &23&\SC_{24}'\r &1152&    &\onefam1  &\bU(2)\oplus[4]\cr
14  &21&\SC_{24}'\r &1152&    &\onefam2  &\bU^2(2)\cr
\noalign{\vspace{1pt}\hrule}
\crcr}}\hss}
\endtable
Let $\bm:=\bm(D)$ be as in \autoref{tab.SC} or as follows\rom:
\[*
\minitab\quad
D:  &3,5,9&11,13&15&16&\{17\bmod6\}&\{18\bmod2\}&\{19,21\bmod6\}\cr
\bm:&   24&   22&20&21&          18&          20&             19\cr
\endminitab
\]
Then, with the exceptions
listed in \autoref{tab.SC}, one has a sharp bound
$\ls|\Fn X|\le\bm$.
%, where $\bm:=\bm(D)$ is as given in \autoref{tab.SC}.
\endtheorem

The graph $\SC_{32}'$ in \autoref{tab.SC}
is $4$-regular, $3$-arc transitive, distance regular, and
bipartite.
According to D.~Pasechnik (private communication), this is the
graph of points and lines in
the affine plane $\Bbb{F}_4^2$, with one pencil of
parallel lines removed.
The graph contains eight minimal pencils, all of type \smash{$4\tD_4$}.

Two other sufficiently symmetric graphs are $\SC_{28}'$ and $\SC_{24}'$,
containing three and two
minimal pencils, all of type \smash{$4\tD_4$},
respectively. Each graph~$\graph$ is
bipartite, the bicomponents $C_1$, $C_2$ are orbits of $\Aut\graph$,
and the action of $\Aut\graph$ on $C_1\times C_2$ has two orbits,
distinguished by wether vertices are adjacent or not.

The bipartite graphs $\SC_{28}'$, $\SC_{27}'$, and $\SC_{24}'$
are subgraphs of $\SC_{32}'$.

\subsection{Proof of \autoref{th.astral}}\label{proof.astral}
Consider an astral configuration~$S$, fix a minimal fiber
$\fiber\subset\Fn X$ of type \smash{$\tD_4$}, and let
$\graph:=\pencil(\fiber)$: this pencil may have parabolic fibers of type
\smash{$\tD_4$} or \smash{$\tA_5$} and elliptic fibers of type $\bA_p$,
$p\le4$.
Since we assume that $\girth(\Fn S)\ge6$, all
sections of~$\graph$ are pairwise disjoint, and their number is bounded by
\autoref{prop.hyperbolic} as follows:
\[*
\minitab\quad
D=                &2,4,6,8&3,5,7& 9&10&11&12\dts14&15,16&17\dts20&21\dts34&\otherwise\\
\ls|\sec\graph|\le&     12&   11& 9& 8& 6&       4&    3&       2&       1&         0\\
\endminitab
\]
(Strictly speaking, for large values of~$D$ these bounds
use an additional assumption
that the pencil admits at least one extra line $l\in\graph\sminus\fiber$.)

To simplify the further computation, we eliminate several large sets of
sections. More precisely, we have the following lemma, which is proved in
\autoref{proof.astral.large}, after the necessary terminology is developed.

\lemma[see \autoref{proof.astral.large}]\label{lem.astral.large}
In the settings of the theorem, one has\rom:
\roster
\item\label{astral.large.12}
if $\ls|\sec\graph|=12$, then
$\ls|\Fn S|\le24$ or
$\Fn S\cong\SC_{27}'$ \rom($D=2,4$\rom) or $\SC_{32}'$ \rom($D=8$\rom)\rom;
\item\label{astral.large.11}
if $D=5$ and $\ls|\sec\graph|\ge11$, then
$\ls|\Fn S|\le22$.
\endroster
\endlemma

%The three section sets of size~$12$ are eliminated as in
%\autoref{ss.step.3.sets}: the corresponding configurations~$S_5$ have rank~$18$
%and satisfy the hypotheses of \autoref{lem.disjoint}. For any geometric
%extension $S\supset S_5$, we have $\ls|\Fn S|\le24$ or
%$\Fn S\cong\SC_{27}'$ (for $D=2,4$) or $\SC_{32}'$ (for $D=8$).
%Similarly, arguing as in \autoref{s.alg.quad}, we can eliminate the three
%section sets of size~$11$ for $D=5$: any extension of these configurations
%has at most $21$ lines.

There remains to run the algorithm of \autoref{s.alg.le} and list geometric
configurations~$S$ satisfying the inequality $\ls|\Fn S|>\Mm$, where the
threshold $\Mm$ is as in the statement of the theorem and $\Ms$ is given by
the table above, with $11$ reduced down to~$10$ for $D=5$ and $12$ reduced
down to~$11$, both due to \autoref{lem.astral.large}.
%(This special treatment of several large sets of sections,
%although not necessary, reduces the amount of computations.)
\qed

\subsection{Geometry of special octics}\label{s.special}
In this section, we assume $D=4$. Recall that an octic $X\subset\Cp5$ is
special if and only if its N\'{e}ron--Severi lattice has a cubic pencil as in
\autoref{def.geometric}\iref{bad.cubic}, \ie, a class
$e\in\NS(X)$ such that $e^2=0$ and $e\cdot h=3$.

Let $l\in\Fn X$ be a line. The determinant of the lattice spanned by $h$,
$l$, and~$e$ is $-8t^2+6t+18$, where $t:=l\cdot e$; hence, this lattice is
hyperbolic if and only if $l\cdot e=0,\pm1$. If $l\cdot e=-1$, then
the difference $e-l$
is a quadric pencil as in
\autoref{def.geometric}\iref{bad.hyperelliptic}.
Thus, we conclude that
\[
l\cdot e=0\ \text{or}\ 1\quad
 \text{for any line $l\in\Fn X$}.
\]

Similarly,
%we can show that,
if a cubic pencil $e\in\NS(X)$ exists, it is
unique. Indeed, if there were two such classes $e_1\ne e_2$, the requirement
that $\NS(X)$ should be hyperbolic would imply that $e_1\cdot e_2=1$ or~$2$. In
the former case, $e_1-e_2$ would be an exceptional divisor as in
\autoref{def.geometric}\iref{bad.exceptional}; in the latter, $h-e_1-e_2$
would be a quadric pencil as in
\autoref{def.geometric}\iref{bad.hyperelliptic}.
Thus, below we denote this unique class by $e_X\in\NS(X)$
and refer to it as \emph{the} cubic pencil in~$X$.
As an immediate consequence, the graph $\Fn X$ splits canonically into two sets
\[
C_n=C_n(X):=\bigl\{l\in\Fn X\bigm|l\cdot e_X=n\bigr\},\quad n=0,1.
\label{eq.C12}
\]

\definition\label{def.biquad}
A \emph{biquadrangle} in a graph~$\graph$ is a complete bipartite subgraph
$K_{2,3}\subset\graph$. In other words, a biquadrangle consists of a
quadrangle (type~$\tA_3$ fiber) $\fiber=\{l_1,l_2,l_3,l_4\}$
and a bisection~$s_{13}$ adjacent to two opposite corners
$l_1,l_3\in\fiber$.
\enddefinition

\definition\label{def.principal}
A bipartite graph~$\graph$ is called \emph{principal} if it has a
bicomponent, also called \emph{principal}, containing all vertices of valency
greater than two.
\enddefinition

\proposition\label{prop.special}
Let $X\subset\Cp5$ be a smooth octic. If $X$ is special, then either
\roster
\item\label{special.trig}
$\girth(\Fn X)=3$, or
\item\label{special.biquad}
$\girth(\Fn X)=4$ and $\Fn X$ has a biquadrangle
\rom(then, $\Fn X$ is bipartite\rom), or
\item\label{special.principal}
$\girth(\Fn X)\ge6$ and $C_1(X)$ is a principal bicomponent of $\Fn X$.
%$\Fn X$ is a principal bipartite graph\rom; moreover,
%the set $C_1:=\{l\in\Fn X\,|\,l\cdot e_X=1\}$ is a principal bicomponent.
\endroster
Conversely, if $\Fn X$ is as in items~\iref{special.trig}
or~\iref{special.biquad} above, then $X$ is special.
\endproposition

\proof
Assume that $X$ is special and let $e:=e_X$.

Consider two lines $l_1,l_2\in\Fn X$, $l_1\cdot l_2=1$. If
$l_1\cdot e=l_2\cdot e=0$, then the class $l_3:=e-l_1-l_2$ is also a line, and
$\{l_1,l_2,l_3\}$ is a triangle. Conversely, if $\{l_1,l_2,l_3\}$ is a
triangle, then the sum $e:=l_1+l_2+l_3$ is a cubic pencil and $X$ is
special.

In the rest of the proof, we assume that $\girth(\Fn X)\ge4$.

If $l_1\cdot e=l_2\cdot e=1$, then $l_1+l_2+2e-h$ is an exceptional divisor,
contradicting to the assumption that $X$ is smooth. Thus,
$C_0$ and~$C_1$ are a pair of distinguished complimentary
bicomponents and $\Fn X$ is bipartite.
Since $h$ is fixed, in the relation
\[*
h=e+l_1+l_2+l_3+l_4+s_{13},
\]
the presence of any five classes implies the presence of the sixth.
%(The notation is explained below.)
(This relation holds whenever $l_1,l_2,l_3,l_4,s_{13}$ are as in
\autoref{def.biquad} and $l_1,l_3\in C_0$.)
Hence, the following statements are equivalent:
\roster*
\item
the graph $\Fn X$ has a biquadrangle,
consisting of a quadrangle $\{l_1,l_2,l_3,l_4\}$
and a bisection~$s_{13}$ intersecting~$l_1$ and~$l_3$;
\item
$X$ is special and the graph $\Fn X$ has a quadrangle $\{l_1,l_2,l_3,l_4\}$;
here, we can assume that $l_1,l_3\in C_0$ and $l_2,l_4\in C_1$;
\item
$X$ is special and there is a vertex $l_1\in C_0$ such that
%valency greater than two,
$\val l_1>2$,
\ie, incident to at least three
other vertices $l_2$, $l_4$, $s_{13}$ (and thus, $C_1$ is not principal).
\endroster
These observations complete the proof.
\endproof

\corollary[of the proof and the uniqueness of $e_X$]\label{cor.special}
In
%case~\iref{special.biquad} of the statement of \autoref{prop.special},
\autoref{prop.special}\iref{special.biquad},
one
has $\val v\le3$ for each vertex $v\in C_0$. In particular,
each minimal pencil has exactly one bisection
\rom(\cf. \autoref{def.biquad}\rom).
\done
\endcorollary

\proposition\label{prop.C12}
Assume that $X$ is special and $\girth(\Fn X)\ge4$, so that
the graph $\Fn X=C_0\cup C_1$ is bipartite.
Then the following statements hold\rom:
\roster
\item\label{special.C1}
$\ls|C_1|\le7$ or $\ls|C_1|=9$\rom;
in the latter case, $\val v=3$ for each vertex $v\in C_0$\rom;
\item\label{special.C0}
$\ls|C_0|\le12$ and, if $\ls|C_0|\ge9$,
then $\val v\le6$ for each vertex $v\in C_1$.
\endroster
\endproposition

\proof
For Statement~\iref{special.C1}, it suffices to observe that any nine lines
$l_1,\ldots,l_9\in C_1$ are subject to the relation
$3h=l_1+\ldots+l_9+5e_X$.

The bound $\ls|C_0|\le12$ follows from the fact that each line $l\in C_0$ is
a component of a type~\smash{$\tA_1$} (or \smash{$\tA_1^*$})
fiber of the elliptic pencil
$\pi\:X\to\Cp1$ defined by~$e_X$ (the other component of this fiber being the
conic $e_X-l$). The valency bound is mainly given by the inertia index of the
lattice; in the border case $C_0=\{l_1,\ldots,l_9\}$ and a line $l\in C_1$
intersecting $l_1,\ldots,l_7$, the class
\[*
l_1+\ldots+l_7-l_8-l_9+3l-h
\]
is an exceptional divisor as in
\autoref{def.geometric}\iref{bad.exceptional}.
(In fact, the valency bound does not use the assumption that $X$ should be
special.)
\endproof

\subsection{Proof of \autoref{th.special}}\label{proof.special}
In view of \autoref{prop.special}, there only remains to prove the bounds on
the number of lines in a special octic~$X$.

If $\girth(\Fn X)=3$, we pick a triangle $\fiber\subset\Fn X$ and consider
$\graph:=\pencil(\fiber)$. The sections of~$\graph$ are pairwise disjoint by
\autoref{prop.valency.trig}, and
\autoref{prop.hyperbolic} yields $\ls|\sec\graph|\le9$;
hence, $\ls|\Fn X|\le33$, as stated. For the more refined statement,
%as in the theorem,
we run the algorithm of \autoref{s.alg.le} and list all
configurations with more than $\Mm:=26$ lines.

\figure
\hbox to\hsize{\hss$
\makeatletter
\setbox0\hbox{$\m@th\joinrel\relbar\joinrel\relbar\joinrel$}%
\def\zbox#1{\kern-\wd0\hbox to\wd0{\hss$#1$\hss}}%
\def\1{\zbox\bullet}\def\.{}\let\\\cr
\vcenter{\offinterlineskip\halign{&\copy0\zbox|#\cr
\1&\1&\1&\1&\.&\.&\.&\.&\.&\.&\.&\.\\
\.&\.&\.&\.&\1&\1&\1&\1&\.&\.&\.&\.\\
\.&\.&\.&\.&\.&\.&\.&\.&\1&\1&\1&\1\\
\1&\1&\.&\.&\1&\1&\.&\.&\1&\1&\.&\.\\
\1&\1&\.&\.&\.&\.&\1&\1&\.&\.&\1&\1\\
\.&\.&\1&\1&\1&\1&\.&\.&\.&\.&\1&\1\\
\.&\.&\1&\1&\.&\.&\1&\1&\1&\1&\.&\.\\
\.&\.&\.&\.&\.&\.&\.&\.&\.&\.&\.&\.\\
\.&\.&\.&\.&\.&\.&\.&\.&\.&\.&\.&\.\\}}
$\hss}
\makeatother

\caption{The configuration $\QC_{21}$ (see \autoref{th.special})}\label{fig.QC21}
\endfigure

Assume that $\girth(\Fn X)\ge4$ and, hence, $\Fn X=C_0\cup C_1$ is
bipartite. By \autoref{prop.C12}, we have either $\ls|\Fn X|\le19$ or
$\ls|C_1|=9$ and $\ls|C_0|\ge11$; in the latter case,
$\ls|\Fn X|=20$ or~$21$ and, since $C_1$ is not a principal
bicomponent, $\Fn X$ has a biquadrangle by \autoref{prop.special}.
For the refined classification, we run an algorithm similar to
\autoref{s.alg.le}, using a ``pencil''
$C_0=\{l_1,\ldots,l_n\}$, $n=11$, $12$, and adding ``sections'' $s\in C_1$
(regarded as subsets of $C_0$,
\cf. \autoref{conv.coord}).
To form a biquadrangle, we start with a triple
of sections $s_1,s_2,s_3$ such that
\[*
\bs_1\cap\bs_2=\bs_2\cap\bs_3=\bs_3\cap\bs_1=\{l_1,l_2\};
\]
there are but three such
triples. Then, at most three more sections can be added increasing the
rank of the lattice, and we arrive at the dichotomy
$\ls|\Fn X|\le19$ or $\Fn X\cong\QC_{21}$
(see \autoref{fig.QC21}).
\qed

\remark
The computation at the end of the proof can be taken to a more combinatorial
level if one takes into account the intrinsic structure on the set $C_0$ of
size~$12$ or~$11$ arising from the
imprimitive embedding $\Z C_0+\Z h+\Z c_X\into\L$
(\cf. pencils of type $(0,12)$ and $(0,11)$
in~\cite{DIS}).
We leave this exercise to the reader.
\endremark

One can use Propositions~\ref{prop.special} and~\ref{prop.C12} to take the
classification one step further and find the extremal configurations
that are neither triangular nor quadrangular.
For such a graph $\graph=C_0\cup C_1$, we have $\ls|\graph|<19$
unless $\ls|C_0|=12$ and $\ls|C_1|=7$.
This time, the class $c_X$ does not need to belong to $\Fano_h(X)$;
hence, we replace this lattice with the primitive hull of
$\Fano_h(X)+\Z c_X$.
The result of the straightforward computation
using \autoref{s.alg.le} is as follows.

\proposition\label{prop.dual}
Let $X\subset\Cp5$ be a smooth special octic such that
$\girth(\Fn X)\ge5$.
Then $\Fn X\cong\SC_{19}'$, $\SC_{19}''$
\rom(see \autoref{fig.SC19}\rom)
or $\ls|\Fn X|\le18$.
The extremal Fano graphs $\SC_{19}'$, $\SC_{19}''$
can also be realized by smooth triquadrics.
\figure
\hbox to\hsize{\hss$
\makeatletter
\setbox0\hbox{$\m@th\joinrel\relbar\joinrel\relbar\joinrel$}%
\def\zbox#1{\kern-\wd0\hbox to\wd0{\hss$#1$\hss}}%
\def\1{\zbox\bullet}\def\.{}\let\\\cr
\vcenter{\offinterlineskip\halign{&\copy0\zbox|#\cr
\1&\1&\1&\1&\.&\.&\.&\.&\.&\.&\.&\.\\
\.&\.&\.&\.&\1&\1&\1&\1&\.&\.&\.&\.\\
\1&\.&\.&\.&\1&\.&\.&\.&\1&\1&\.&\.\\
\.&\1&\.&\.&\.&\1&\.&\.&\.&\.&\1&\1\\
\.&\.&\1&\.&\.&\.&\1&\.&\1&\.&\1&\.\\
\.&\.&\.&\1&\.&\.&\.&\1&\.&\1&\.&\1\\
\.&\.&\.&\.&\.&\.&\.&\.&\.&\.&\.&\.\\}}
\qquad
\vcenter{\offinterlineskip\halign{&\copy0\zbox|#\cr
\1&\1&\1&\1&\.&\.&\.&\.&\.&\.&\.&\.\\
\.&\.&\.&\.&\1&\1&\1&\1&\.&\.&\.&\.\\
\1&\.&\.&\.&\1&\.&\.&\.&\1&\.&\.&\.\\
\.&\1&\.&\.&\.&\1&\.&\.&\1&\.&\.&\.\\
\.&\.&\1&\.&\.&\.&\1&\.&\.&\1&\.&\.\\
\.&\.&\.&\1&\.&\.&\.&\1&\.&\1&\.&\.\\
\.&\.&\.&\.&\.&\.&\.&\.&\.&\.&\1&\1\\}}
$\hss}
\makeatother

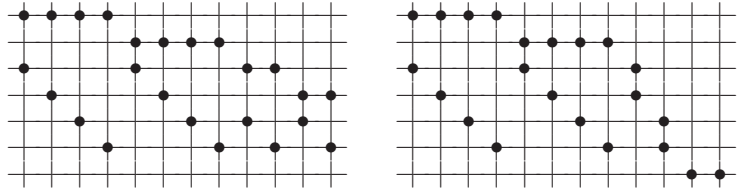
\caption{The configurations $\SC_{19}'$ and $\SC_{19}''$ (see \autoref{prop.dual})}\label{fig.SC19}
\endfigure
\done
\endproposition

\subsection{Hyperelliptic models}\label{s.hyperelliptic}
A \emph{hyperelliptic model} of degree~$2D$ of a $K3$-surface $X$ is a
two-to-one map $X\to Y\into\Cp{D+1}$ splitting through a rational
surface~$Y$. Due to~\cite{Saint-Donat}, a degree~$2D$ model $X\to\Cp{D+1}$
defined by a class $h\in\NS(X)$ is hyperelliptic if and only if
one of the following conditions holds:
\roster
\item\label{he.D=1}
$D=1$: any model $X\to\Cp2$ is hyperelliptic;
\item\label{he.Veronese}
$D=4$ and $h$ is divisible by~$2$ in $\NS(X)$:
the model is the composition of a double covering $X\to\Cp2$ and the Veronese
embedding $\Cp2\into\Cp5$;
\item\label{he.quadric}
$D\ge2$ and there is a quadric pencil $e\in\NS(X)$ as in
\autoref{def.geometric}\iref{bad.hyperelliptic}: the rational surface
$Y\subset\Cp{D+1}$ is a scroll
%$\Sigma_{D-2}$
in this case (see also \cite{Reid}).
\endroster
A hyperelliptic model $X\to Y\into\Cp{D+1}$ is considered \emph{smooth} if so
is its ramification locus in~$Y$. A \emph{line} in a hyperelliptic model is a
rational curve $L\subset X$ that is mapped one-to-one onto a line
in~$\Cp{D+1}$.
Thus, arithmetically, a line
%is still
can still be described as
a class $l\in\NS(X)$ such that $l^2=-2$ and $l\cdot h=1$.

We will treat the case $D=1$ in a subsequent paper:
the conjectural bound for such models is $144$
lines (see~\cite{degt:singular.K3}).

In Case~\iref{he.Veronese},
where the model splits through the
Veronese embedding, it has no lines, as
$l\cdot h=0\bmod2$ for any $l\in\NS(X)$.

Thus, we are left with
Case~\iref{he.quadric}, where the image~$Y$
%$Y\cong\Sigma_{D-2}$
is a scroll.
%with an exceptional section~$E$ of self-intersection $(2-D)$.
To study such models, we
need to change the definition of the admissibility of a configuration~$S$ (see
\autoref{def.geometric}): the exceptional divisors as
in~\iref{bad.exceptional} are still excluded, but this time we \emph{assume}
the presence of a quadric pencil $e\in S$ as in~\iref{bad.hyperelliptic},
requiring instead that there should be no class $c\in S$ such that
\roster
\item[$(2')$]
$c^2=0$ and $c\cdot h=1$ (\emph{fixed component}, see \cite{Nikulin:Weil}).
\endroster
Now, arguing as in \autoref{s.special}
and using an obvious analogue of \autoref{prop.hyperbolic},
we can easily prove the following statement.

\lemma\label{lem.hyperelliptic}
Consider a smooth hyperelliptic model $X\to\Cp{D+1}$, $D\ge2$, and let
$e\in\NS(X)$ be a quadric pencil. Then, for any line $l\in\Fn X$, we have
$l\cdot e=0$ or~$1$.
%\rom; the latter is possible only if $D\in\{2,3,4,6\}$.
Furthermore, if $D\ge3$, the quadric pencil $e:=e_X\in\NS(X)$ is unique.
\done
\endlemma

If $D=2$, the obvious maximum of $48$ lines
($24$ over each of the two rulings) is realized by a connected
equilinear
family containing, in particular, the discriminant minimizing
singular surfaces (see~\cite{degt:singular.K3}).

Thus, from now on, we assume that $D\ge3$,
and \autoref{lem.hyperelliptic}
gives us a natural splitting $\Fn X=C_0\cup C_1$ as
in~\eqref{eq.C12}.
The set~$C_0$ consists of the components of the
reducible fibers of an elliptic pencil ---the pullback of the ruling of~$Y$;
these fibers are of type~\smash{$\tA_1$} (or~\smash{$\tA_1^*$}), we have
$l'\cdot l''=2$ for two lines $l'$, $l''$ in the same fiber, and each line
$l\in C_1$ intersects exactly one of~$l'$, $l''$.
Further computation shows that the set~$C_1$ is empty
with the following few exceptions:
\roster*
\item
$D=3$ and $C_1=\{l',l''\}$, $l'\cdot l''=1$ and $h=l'+l''+2e$;
\item
$D=4$ and $C_1=\{l',l''\}$, $l'\cdot l''=0$ and $h=l'+l''+3e$;
\item
$D=6$ and $C_1=\{l\}$ so that $h=2l+5e$; hence, $C_0=\varnothing$.
\endroster
Geometrically, $C_1$ projects to the exceptional section
$E\subset Y$; if $D=6$ ($E^2=-4$), this section is a component of the
ramification locus.

\theorem\label{th.hyperelliptic}
Let $X\to\Cp{D+1}$ be a hyperelliptic model, $D\ge3$. Then, unless $D=6$ and
$\ls|\Fn X|=1$,
we have $\ls|\Fn X|=0\bmod2$ and $\ls|\Fn X|\le24-4\epsilon$,
where we let $\epsilon:=D\bmod2\in\{0,1\}$.
This bound is sharp.
%\done
\endtheorem

\proof
If $C_1=\varnothing$, the configuration consists of a single pencil~$C_0$
and, by an obvious analogue of \autoref{prop.parabolic}, the maximal number
of lines depends on $D\bmod2$ only. Using \autoref{alg.Nikulin}, we arrive at
$\ls|\Fn X|\le24-4\epsilon$.
The other case $\ls|C_1|=2$ is treated similarly: for each value $D=3,4$ and
each cardinality $\ls|C_0|$, there is a single candidate, which is
analyzed using \autoref{alg.Nikulin}.
%Accidentally, we obtain the same bound
%$\ls|C_0|\le22-4\epsilon$.
\endproof

\section{Pentagonal configurations}\label{S.pent}

In the remaining three sections, we consider configurations with minimal
pencils whose sections may intersect each other.
%Recall that a graph $\graph$ is \emph{quadrangle free} if
%$\girth(\graph)\ge5$, \ie, $\graph$ contains neither triangles nor quadrangles.

\subsection{Statements}\label{s.quad-free.statements}
%In this section, we consider quadrangle free configurations that are not
%locally elliptic, \ie, have a vertex of valency at least~$4$.
%The principal result is the following theorem.
Recall that a graph~$\graph$ is \emph{pentagonal} if $\girth(\graph)=5$.
All large geometric pentagonal graphs are described by the following theorem.

\theorem[see \autoref{proof.a4}]\label{th.a4}
Let $X$ be a smooth $2D$-polarized $K3$-surface,
and assume that
the Fano graph
$\Fn X$ is pentagonal, \ie, $\girth(\Fn X)=5$.
\table\rm
\def\-{\rlap{$\mdag$}}
\def\*{\rlap{$\mstar$}}
\def\r{\rlap{$\mreal$}}
\def\={\relax\afterassignment\ul\count0=}
\def\ul{\underline{\the\count0}}
\caption{Extremal pentagonal configurations (see \autoref{th.a4})}\label{tab.a4}
\def\*{\mstar}
\hbox to\hsize{\hss\vbox{\halign{\strut\quad\hss$#$\hss\quad&
 \hss$#$\hss\quad&\hss$#$\hss\quad&\hss$#$\hss\quad&\hss$#$\hss\quad&
 \hss$#$\hss\quad&\hss$#$\hss\quad&$#$\hss\quad\cr
\noalign{\hrule\vspace{2pt}}
D&\bm&\graph&\ls|\Aut\graph|&\det&(r,c)&\ls|\Aut X|&T:=\Fano_h(\graph)^\perp\cr
\noalign{\vspace{1pt}\hrule\vspace{2pt}}
 2  &29&\QF_{30}'   & 240& 220&(2,0)&  60&[4,2,56],\ [16,6,16]\cr
      &&\QF_{30}''  &  40& 255&(1,1)&  10&[4,1,64]\*,\ [16,1,16]\cr
 3  &29&\QF_{30}'   & 240& 180&(1,0)&  60&[14,4,14]\cr
      &&\QF_{30}''  &  40& 195&(1,0)&  10&[14,1,14]\cr
 4  &29&\QF_{30}'   & 240& 140&(1,1)&  60&[12,2,12]\cr
      &&\QF_{30}''  &  40& 135&(1,1)&  10&[12,3,12]\cr
 5\-&28&\QF_{30}'   & 240& 100&(1,0)&  60&[10,0,10]\cr
      &&\QF_{30}''  &  40&  75&(1,0)&  10&[10,5,10]\cr
      &&\QF_{30}''' &  24&  36&(1,0)&  12&[4,2,10]\cr
 6  &28&\QF_{30}'   & 240&  60&(1,0)&  60&[8,2,8]\cr
      &&\QF_{36}''\r&1440&\=15&(1,0)& 720&[2,1,8]\cr
 7  &26&\QF_{30}'   & 240&\=20&(1,0)& 120&[4,2,6]\cr
 8  &24&\QF_{25}    &  80&    &\onefam1  &\bU(5)\oplus[4]\cr
 9  &23&\QF_{25}\r  &  80&\=15&\onefam1  &\bU(5)\oplus[2]\cr
10\-&21&\QF_{25}\r  &  80&    &\onefam2  &\bU(5)\oplus\bU\cr
\noalign{\vspace{1pt}\hrule}
\crcr}}\hss}
\endtable
Let $\bm:=\bm(D)$ be as in \autoref{tab.a4} or as follows\rom:
\[*
\minitab\quad
D:  &\{10,11\bmod5\}&\{12\bmod5\}&\otherwise\cr
\bm:&                21&       22&        20\cr
\endminitab
\]
Then, with the exceptions
listed in \autoref{tab.a4}, one has a sharp bound
$\ls|\Fn X|\le\bm$.
\endtheorem

The notation used in
%the main part of
\autoref{tab.SC} is similar to
\autoref{tab.main} (with some columns skipped);
it is explained in \autoref{s.results}.
It is remarkable that the same graph $\QF_{30}'$ is a Fano graph
of $2D$-polarized $K3$-surfaces for all $2\le D\le7$. The graph
$\QF_{30}''$ works for $2\le D\le5$; if $D=6$, this graph is still
geometric, but not saturated (pencils acquire extra fibers
and/or sections), and
$\QF_{36}''\supset\QF_{30}''$ is its saturation.

The graph $\QF_{36}''$
is $5$-regular, $2$-arc transitive, and distance regular;
according to D.~Pasechnik (private communication),
it is the \emph{Sylvester graph}. The graph contains $36$ minimal pencils,
all of type \smash{$4\tA_4+\bA_1$}.

The graph $\QF_{30}'$ is $4$-regular and $1$-arc transitive;
it contains six copies of \smash{$4\tA_4$}.

%\roster*
%\item
%$\QF_{36}''$:\enspace
%$\pents{[ [ 4, 0, 0, 1 ], 36 ]}$;
%this graph is $5$-regular, $2$-arc transitive, and distance regular;
%it is the \emph{Sylvester graph}
%(D.~Pasechnik);
%\item
%$\QF_{30}'$:\enspace
%$\pents{[ [ 4, 0, 0, 0 ], 6 ]}$;
%the graph is $4$-regular and $1$-arc transitive;
%\item
%$\QF_{30}''$:\enspace
%$\pents{[ [ 4, 0, 0, 0 ], 1 ], [ [ 3, 0, 1, 1 ], 20 ]}$;
%\item
%$\QF_{25}$:\enspace
%$\pents{[ [ 4, 0, 0, 0 ], 1 ]}$.
%\endroster

\subsection{Sections at a fiber~\pdfstr{Sigma}{$\fiber$} of type~\pdfstr{A\_4}{$\tA_4$}}\label{s.a4.sections}
Consider a geometric pentagonal configuration~$S$ and pick a distinguished
pentagon $\fiber:=\{l_1,\ldots,l_5\}\subset\Fn S$.

\convention\label{conv.numbering}
Here and below, we always assume that the $n$ lines $l_1,\ldots,l_n$ of a
type~\smash{$\tA_{n-1}$} fiber
in a configuration are numbered cyclically, so that
\[*
l_{i\pm n}=l_i,\qquad
l_i\cdot l_j=1\quad \text{if and only if $i-j=\pm1\bmod n$}.
\]
Thus, we regard the indices as elements of the cyclic group $\Z/n$.
\endconvention

Consider the pencil $\graph:=\pencil(\fiber)$. Since $\Fn S$ has no triangles or
quadrangles, we immediately obtain the following lemma.

\lemma\label{lem.pent.disjoint}
For $\fiber\subset\graph\subset\Fn S$ as above,
the following statements hold\rom:
\roster*
\item
all sections of $\graph$ are simple\rom;
\item
$(\sec l_i)\*(\sec l_j)=0$ unless $i-j=\pm2\bmod5$\rom;
\item
$(\sec l_i)\*(\sec l_j)\le1$ if $i-j=\pm2\bmod5$.
\endroster
In particular, $\fiber$ is suitable for enumeration in the sense of
\autoref{s.alg.sections}.
\endlemma

\lemma\label{lem.pent}
For $\fiber\subset\graph\subset\Fn S$ as above, one has the following bounds
on $\ls|\sec\graph|$\rom:
\[*
\minitab\quad
D=                &2,3& 4& 5& 6& 7&8&9,10&11,12&13\dts15&16\dts25&\otherwise\\
\ls|\sec\graph|\le& 16&17&18&15&10&7&   5&    3&       2&       1&         0\\
\endminitab
\]
Furthermore, the following statements hold\rom:
\roster
\item\label{pent.2}
if $2\le D\le4$ and $\ls|\sec\graph|\ge14$, then $\ls|\Fn S|\le29$\rom;
\item\label{pent.5}
if $D=5$ and $\ls|\sec\graph|\ge14$, then
$\Fn S\cong\QF_{30}'''$ or $\ls|\Fn S|\le28$\rom;
\item\label{pent.6}
if $D=6$ and $\ls|\sec\graph|\ge12$, then
$\Fn S\cong\QF_{36}''$ or $\ls|\Fn S|\le28$\rom;
\item\label{pent.7}
if $7\le D\le9$, then there is at most one pair of intersecting sections\rom;
\item\label{pent.10}
if $D\ge10$, then all sections of~$\graph$ are pairwise disjoint.
\endroster
Besides, if $D=5$, then $\val l\in\{0\dts3,5\}$ for each line $l\in\fiber$.
\endlemma

\proof
Since $\fiber$ is suitable for enumeration, we can run the algorithm of
\autoref{s.alg.sections}.

If $D\ge7$, we can easily enumerate \emph{all} sets
of sections, obtaining, in particular, Statements~\iref{pent.7}
and~\iref{pent.10}. (These lists are also used in \autoref{proof.a4} below.)

If $D\le6$, we only list the sets $\sec\graph$ with $\ls|\sec\graph|\ge14$
($\ls|\sec\graph|\ge12$ if $D=6$). Then,
Statements~\iref{pent.2}--\iref{pent.6} are proved
by applying Step~$3$ (see \autoref{ss.step.3.sets})
in order to classify large configurations. The technical details are
explained in \autoref{s.alg.pent}.

The last statement is straightforward.
\endproof

%\subsection{Pencils with a fiber of type \pdfstr{A\_4}{$\tA_4$}}\label{s.a4.pencils}

\subsection{Proof of \autoref{th.a4}}\label{proof.a4}
Following the idea explained in \autoref{s.idea}, we can still use the
algorithm of \autoref{s.alg.le} based on discrete sets of sections.

Let $\Mm$ be as in the statement of the theorem. If $D\ge14$, there are
relatively few sets of sections $\sec\graph$ (see \autoref{lem.pent} and its
proof), and we can use the algorithm directly, taking into account
(if present) the one pair of sections that intersect.

If $D\le12$, we use the following trick. Let $\Ms=13$ ($\Ms=11$ if $D=12$).
Then, by \autoref{lem.pent}, it suffices to consider pencils~$\graph$ with
$\ls|\graph|>\Mg:=\Mm-\Ms$.
Pick such a pencil~$\graph$ and let
$n_i:=\val l_i$, $i\in\Z/5$, be a desired set of goals, such that
%$\Mm-\ls|\graph|<\sum_in_i\le\Ms$.
\[
\Mm-\ls|\graph|<\sum_{i\in\Z/5}n_i\le\Ms.
\label{eq.pent.sum}
\]
(In addition, we can require that each pair
$(n_i,n_{i+1})$, $i\in\Z/5$, is realized by the valencies of a pair of
adjacent lines of~$\fiber$ in a geometric configuration: this condition is
easily checked as all sections involved are pairwise disjoint, see
\autoref{lem.pent.disjoint}).
Let, further, $m:=\max\{n_i+n_{i+1}\,|\,i\in\Z/5\}$, and consider the set
\[*
\CN(n_i):=\bigl\{(p,q)\bigm|
 \text{$p\ge q$, $p+q=m$, $(p,q)=(n_i,n_{i\pm1})$ for some $i\in\Z/5$}\bigr\}.
\]
Then, denoting by~$\CN$ the union of $\CN(n_i)$ over all sets of
goals~$(n_i)$ satisfying~\eqref{eq.pent.sum},
we can assert that, in any geometric pentagonal extension
$S\supset\Fano_h(\graph)$ such that $\ls|\Fn X|>\Mm$, there is a pair
$(p,q)\in\CN$ such that $\graph$ has at least $(p+q)$ disjoint sections:
%furthermore,
up to automorphism of~$\fiber$, we can assume that
$\val l_1\ge p$ and $\val l_2\ge q$.

If $D\in\{2\dts4,6\}$, the computation based on this observation completes
the proof, as all configurations obtained by adding to the pencil $(p+q)$
disjoint sections as above have the maximal rank~$20$ and can be analyzed
using \autoref{alg.Nikulin}.

If $D=5$, then, occasionally, we may need to add an extra section.
Namely, let $(n_i)$ be a set of goals such that $2m\le\sum_in_i$:
this is the case for all configurations of low rank.
Then, a simple argument shows that
%, up to the action of the dihedral group~$\DG{10}$, we can assume
we can reindex the lines so
that $n_1\ge n_2$, $n_1+n_2=m$ is maximal, and
$n_3>0$, \ie, $\graph$ has at least one section adjacent to~$l_3$.
Taking this extra section into account and allowing it to intersect one of
those adjacent to~$l_1$, we always obtain a configuration of the maximal
rank~$20$, and there remains to apply \autoref{alg.Nikulin}.
\qed

\section{Quadrangular configurations}\label{S.quad}

Recall that a graph $\graph$
and the corresponding configuration are \emph{quadrangular} if
$\girth(\graph)=4$, \ie, $\graph$ has no triangles,
but it does contain a quadrangle.

\subsection{Statements}\label{s.quad.statements}
In this section, we mostly ignore the case $D=2$, although
the maximal cardinality of
a triangle free configuration of lines in a
spatial quartic remains an open problem (see the discussion in
\autoref{s.quartics}).
The next statement is the principal result of the section.

\theorem[see \autoref{proof.quad}]\label{th.quad}
Let $X\subset\Cp{D+1}$ be a smooth $2D$-polarized $K3$-surface,
$D\ge3$, and assume that the Fano graph
$\Fn X$ is quadrangular, \ie, $\girth(\Fn X)=4$.
Then, with few exceptions, we have the following sharp bounds\rom:
\[*
\minitab\enspace
D=           & 3& 4&5,6& 7& 9&\{8\bmod4\}&\{10\bmod4\}&\{11\bmod4\}&\{13\bmod4\}\cr
\ls|\Fn X|\le&36&30& 28&22&20&         22&          24&          20&          18\cr
\endminitab
\]
The exception is the case $D=4$, $\Fn X\cong\QC_*^*$
\rom(see \autoref{tab.main} and comments below\rom).
\endtheorem

\addendum[see \autoref{proof.quad}]\label{ad.quad}
The complete set $\{\ls|\Fn X|\}$ of values taken by the line count
of a quadrangular configuration is as follows\rom:
\roster*
\item
if $D=3$, then $\{\ls|\Fn X|\}=\{4\dts36\}$\rom;
\item
if $D=4$, then $\{\ls|\Fn X|\}=\{4\dts30,32\dts34,36\}$.
\endroster
%
%If $D=4$, all extremal Fano graphs $\QC_*^*$ are represented by
%triquadrics.\mnote{\todo: to reconsider}
%
%If $D\ge24$, any configuration maximizing the number of lines is parabolic.
\endaddendum

According to \autoref{th.special}, the exceptional Fano
graphs~$\QC_*^*$ in \autoref{tab.main}
are realized by triquadrics only.
The next statement is proved by a simple computation.

\addendum\label{ad.quad.max}
There are inclusions $\QC_{32}'\subset\QC_{34}'\subset\QC_{36}'$ and, hence,
respective specializations of the corresponding families of triquadrics.
%octic $K3$-surfaces.
All other graphs~$\QC_*^*$ in \autoref{th.quad} are maximal with respect to
inclusion.
\done
\endaddendum

The graphs
$\QC_{36}''$, $\QC_{32}'$, $\QC_{32}''$, $\QC_{32}'''$, and~$\QC_{32}\K$ are
$6$-regular and vertex transitive. The last one, $\QC_{32}\K$, is the famous
\emph{Kummer configuration $16_6$}: it is $2$-arc transitive,
distance regular, and bipartite.
The graphs $\QC_{36}'$ and~$\QC_{32}\K$ can be represented by the discriminant
minimizing singular $K3$-surface $X([4,0,8])$ (see~\cite{degt:singular.K3}).
More details on these graphs and corresponding triquadrics
%(\eg, automorphism groups and transcendental lattices)
are found in \autoref{tab.main} (see also~\cite{degt:Fano.graphs}).

\subsection{Sections at a fiber~\pdfstr{Sigma}{$\fiber$} of type
 \pdfstr{A\_3}{$\tA_3$}}\label{s.quad.sections}
Till the end of this section, we
consider a geometric quadrangular configuration~$S$ and pick a
quadrangle $\fiber:=\{l_1,l_2,l_3,l_4\}$.
We assume the lines ordered cyclically (see
\autoref{conv.numbering}) and, usually, so that
\[
\val l_1\ge\val l_3,\quad
\val l_1\ge\val l_2\ge\val l_4,\quad
\val l_3\ge\val l_4\ \text{if $\val l_1=\val l_2$}.
\label{eq.quad.order}
\]
Consider also the pencil $\graph:=\pencil(\fiber)$;
since $\kk_\fiber=l_1+l_2+l_3+l_4$ and $\deg\graph=4$, all
fibers of~$\graph$ are of types $\tA_3$, $\bA_2$, or $\bA_1$.

We
subdivide quadrangular configurations according to whether they do or do not
have a biquadrangle (see \autoref{def.biquad}).

\lemma\label{lem.quad}
For $\fiber\subset\Fn S$ as above, one has\rom:
\roster
\item\label{quad.discrete}
each subgraph $\sec l_i$, $1\le i\le4$, is discrete\rom;
\item\label{quad.i+1}
$\sec l_i\cap\sec l_j=\varnothing$ whenever $j-i=1\bmod2$\rom;
\item\label{quad.D=3}
if $D=3$, then
$(\sec l_i)\*(\sec l_{i\pm1})\le2$ and $(\sec l_i)\*(\sec l_{i+2})\le3$\rom;
\item\label{quad.D>3}
if $D\ge4$, then $(\sec l_i)\*(\sec l_j)\le1$ for any pair $(i,j)$\rom;
\item\label{quad.D>7}
if $D\ge8$, then the subgraph $\sec\graph$ is discrete.
\endroster
If, in addition, $S$ has no biquadrangle, then also
\roster[\lastitem]
\item\label{quad.i+2}
$\sec l_i\cap\sec l_j=\varnothing$ for any pair $1\le i<j\le4$\rom;
\item\label{quad.D=2,3}
$(\sec l_i)\*(\sec l_{i\pm1})\le1$ and $(\sec l_i)\*(\sec l_{i+2})\le2$.
\endroster
\rom(In view of Statement~$\iref{quad.D>3}$, the last restriction is meaningful only
for $D=2,3$.\rom)
\endlemma

\proof
Most assertions are immediate from the assumption that $S$ should have no
triangle or, for the last two statements, biquadrangle.
Statements~\iref{quad.D=3}, \iref{quad.D>3}, and~\iref{quad.D>7}
are proved by applying \autoref{prop.hyperbolic} to
simple test configurations.
\endproof

\lemma\label{lem.quad.sections}
If $S$ has no biquadrangle,
one has the following
bounds on $\ls|\sec\graph|$\rom:
\[*
\minitab\quad
D=                &2,3& 4& 5&  6&7,8&9,10&11\dts16&\otherwise\cr
\ls|\sec\graph|\le& 21&20&12&\-8&  4&   2&       1&          0\cr
\endminitab
\]
Furthermore, the following statements hold\rom:
\roster*
\item
if $D=2,3$
%, then $\ls|\sec\graph|\le21$\rom; if
and $\ls|\sec\graph|=21$, then $\ls|\Fn S|=34$
\rom(a unique configuration\rom)\rom;
\item
if $D=4$
%, then $\ls|\sec\graph|\le20$\rom; if
and $\ls|\sec\graph|\ge15$, then either $\ls|\Fn S|\le30$
or $\Fn S$ is isomorphic to
one of the exceptional graphs $\QC_{36}'$, $\QC_{36}''$,
$\QC_{34}'$, $\QC_{32}'$, $\QC_{32}''$, $\QC_{32}'''$, $\QC_{32}\K$.
%\item
%if $D=5$, then $\ls|\sec\graph|\le12$\rom;
%\item
%if $D=6$, then $\ls|\sec\graph|\le8$\rom;
%\item
%if $D=7,8$, then $\ls|\sec\graph|\le4$\rom;
%\item
%if $D\ge9$, then $\ls|\sec\graph|\le2$.
\endroster
\endlemma

\proof
If $D\ge4$, then, due to \autoref{lem.quad}, the fiber~$\fiber$ is suitable
for enumeration in the sense of \autoref{s.alg.sections}. Applying the
algorithm
%described in \autoref{s.alg.sections},
%assuming~\eqref{eq.quad.order},
%processing the lines in the order $l_1,l_3,l_2,l_4$, and taking for the
%defining property
of \autoref{s.alg.sections}
(in the order $l_1,l_3,l_2,l_4$),
assuming~\eqref{eq.quad.order} and taking for the defining property
\[
\prop\:S\longmapsto
 \text{$S$ has neither a triangle nor a biquadrangle},
\label{eq.quad.prop}
\]
we arrive at the
statement of the lemma.

If $D=4$ and $\ls|\sec\graph|\ge15$,
the configurations are classified using Step~$3$ of the algorithm (see
\autoref{ss.step.3.sets}); we postpone the proof until \autoref{s.alg.quad},
after the algorithm is explained in more details.
%then we have $\rank S_4\ge18$ and all
%configurations obtained satisfy the hypotheses of
%\autoref{lem.disjoint}.

If $D=2,3$, the situation is more complicated, see
\autoref{lem.quad}\iref{quad.D=2,3}. However,
large configurations can still be
enumerated similar to \autoref{s.alg.sections}. Assuming that $\val l_1\ge6$ and
$\val l_1+\val l_3\ge11$, see~\eqref{eq.quad.order}, we arrive at a unique
configuration~$S$ satisfying the inequality $\ls|\sec\graph|>20$.
This configuration has rank~$20$, contains $34$ lines, and has no nontrivial
$\prop$-geometric finite index extensions.
\endproof

\subsection{Configurations with a biquadrangle}\label{s.biquad}
Assume that the pencil $\graph=\pencil(\fiber)$ has a bisection~$s_{13}$
intersecting~$l_1$ and~$l_3$. Using \autoref{prop.hyperbolic}, one can
%easily
show that $D\le4$ and, unless $D=2$, such a bisection is unique. We ignore
the case $D=2$ (\cf. \autoref{ex.quad.D=2} below),
% and consider the remaining
%two cases $D=3,4$ separately.
whereas the case $D=4$ was considered in \autoref{s.special}, resulting in the
bound $\ls|\Fn S|\le21$ (see \autoref{th.special}).

%\lemma\label{lem.biquad.D=4}
%Assume that $D=4$ and $\graph$ has a bisection~$s_{13}$
%intersecting~$l_1$ and~$l_3$. Then all other
%sections of~$\graph$ are simple and one has
%\[*
%\sec_1\graph=\sec l_2\cup\sec l_4,\qquad
%\ls|\sec_1\graph|\le8,\qquad
%\ls|\Fn S|\le29.
%\]
%%In particular, this implies that $\ls|\Fn S|\le33$.
%Furthermore, such a configuration~$S$ cannot be realized by a triquadric.
%\endlemma
%
%\proof
%The surface cannot be a triquadric as the class
%\[*
%e:=h-l_1-l_2-l_3-l_4-s_{13}
%\]
%is a cubic pencil as in \autoref{def.geometric}\iref{bad.cubic}.
%
%The equality $\sec_1\graph=\sec l_2\cup\sec l_4$
%(equivalent to $\sec l_1=\sec l_3=\varnothing$)
%and the bound $\ls|\sec_1\graph|\le8$ are proved by a simple enumeration as in
%\autoref{s.alg.sections}. (The assumption that the configuration is triangle
%free is crucial.)
%%
%%Finally, using \autoref{prop.hyperbolic}, one can show that
%By \autoref{prop.valency},
%the presence of a bisection implies that $\graph$ has at most four parabolic
%fibers. Then, $\ls|\graph|\le20$ by \autoref{cor.parabolic} (see
%\autoref{tab.parabolic}) and, hence, $\ls|\Fn S|\le29$ by \eqref{eq.count}.
%\endproof

%The rest of this section deals with the case $D=3$, \ie, sextics in~$\Cp4$.

\lemma\label{lem.biquad.D=3}
Assume that $D=3$ and $\graph$ has a bisection~$s_{13}$
intersecting~$l_1$ and~$l_3$.
Then $\ls|\Fn S|\le36$.
\endlemma

\proof
It is immediate that $\graph$ has another bisection,
\viz.
\[*
s_{24}:=h-l_1-l_2-l_3-l_4-s_{13},
\]
intersecting~$l_2$ and~$l_4$,
so that $s_{13}\cdot s_{24}=1$, and
any other line $l\in\Fn S$ intersects exactly one of
%$l_1$, $l_2$, $l_3$, $l_4$, $s_{13}$, $s_{24}$.
$l_1,l_2,l_3,l_4,s_{13},s_{24}$.
\figure
\cpic{biquad}
\caption{A biquadrangle for $D=3$}\label{fig.biquad}
\endfigure
Thus, renaming the lines to $a_1,a_2,a_3,b_1,b_2,b_3$
as in \autoref{fig.biquad}, we have
\[*
\ls|\Fn S|
 =\sum_{i=1}^3\val a_i+\sum_{j=1}^3\val b_j-12
 =\frac13\sum_{i,j=1}^3(\val a_i+\val b_j)-12.
%\label{eq.biquad.count}
\]
Since $\val l\le9$ for each line $l\in\Fn S$ (see
\autoref{prop.valency}), we have $\ls|\Fn S|\le42$.

From now on, assume that $\ls|\Fn S|>36$.
Reindexing the lines (and choosing, if necessary, another quadrangle~$\fiber$),
we can also assume that $\val a_1\ge\val b_1$ and the pair
$(i,j)=(1,1)$ maximizes the sum $\val a_i+\val b_j$, $i,j=1,2,3$.
%Then, in view of~\eqref{eq.biquad.count},
By the count above, $\val a_1=9$ and
$\val b_1=9$ or~$8$; hence, $\ls|\graph|=\val a_1+\val b_1-2=16$ or~$15$, respectively.
In the former case, \eqref{eq.quad.order} implies $\val l_1\ge8$. In the
latter case, we can order the lines so that $\val l_4\le\val l_2\le8$ and
$\val l_3\le\val l_1$; then again $\val l_1\ge8$.

Now, we can run an algorithm similar to
that of \autoref{s.alg.pencils},
taking ``triangle free'' for the defining property~$\prop$ and
starting
with a pencil~$\graph$ of size $\ls|\graph|=16$ or~$15$
and a pair of bisections
$s_{13}$, $s_{24}$. Recall that each line in $\Fn S$ intersects exactly one
of the bisections~$s_{13}$, $s_{24}$; together with the assumption that $S$
should be
triangle free, this observation
determines the pair of bisections uniquely up to
$\Aut\graph$. Furthermore, all other sections of~$\graph$ are simple and
disjoint from $s_{13}$, $s_{24}$. Thus, taking~$l_1$ for~$l$,
at Steps~1 and~2 we add at least
five disjoint sections $s_i\in\sec l$ (recall that $\val l_1\ge8$).
Configurations of low rank (Step~$3$) are analyzed as described in
\autoref{ss.quad}.
\endproof

\subsection{Proof of \autoref{th.quad} and \autoref{ad.quad}}\label{proof.quad}
Let $S:=\Fano_h(X)$ be as in the theorem.
If $S$ has a biquadrangle, the statement is given by
%Lemmas~\ref{lem.biquad.D=4} and~\ref{lem.biquad.D=3}.
\autoref{lem.biquad.D=3} and \autoref{th.special}.
Thus, assume that $S$ is biquadrangle free.

We can use~\eqref{eq.mu} and~\eqref{eq.Euler+} to list all combinatorial types
of pencils~$\graph$ of degree~$4$;
then, Propositions~\ref{prop.parabolic} and~\ref{prop.hyperbolic} imply the
following lemma.

\lemma\label{lem.quad.pencil}
The cardinality $\ls|\graph|$ is bounded as follows\rom:
\[*
\minitab\quad
D\in              &\{6\bmod4\}&\{4\bmod4\}&\{7\bmod4\}&\{5\bmod4\}\cr
\max{\ls|\graph|}=&         24&         22&         20&         18\cr
\endminitab
\]
%\roster*
%\item
%$\ls|\graph|\le24$ if $D\in\{6\bmod4\}$
%\rom(\eg, $\graph\cong6\tA_3$\rom)\rom;
%\item
%$\ls|\graph|\le22$ if $D\in\{4\bmod4\}$
%\rom(\eg, $\graph\cong5\tA_3+2\bA_1$\rom)\rom;
%\item
%$\ls|\graph|\le20$ if $D\in\{7\bmod4\}$
%\rom(\eg, $\graph\cong5\tA_3$\rom)\rom;
%\item
%$\ls|\graph|\le18$ if $D\in\{5\bmod4\}$
%\rom(\eg, $\graph\cong4\tA_3+2\bA_1$\rom).
%\endroster
If $\graph$ has at least one simple section,
the values of~$D$ are bounded as follows\rom:
\roster*
\item
if $\ls|\graph|\ge22$, then $D=4$\rom;
\item
if $\ls|\graph|=21$, then $D\in\{2,4,6\}$\rom;
\item
if $\ls|\graph|=19$, $20$, then $D\in\{2,4,6,8\}$\rom;
\item
if $\ls|\graph|=17$, $18$, then $D\in\{2\dts10\}$.
\done
\endroster
\endlemma

Lemmas~\ref{lem.quad.sections} and~\ref{lem.quad.pencil} imply
the statement of \autoref{th.quad} for $D\ge9$.
%and the last assertion of \autoref{ad.quad}.

For the remaining few values of~$D$, including $D=2$,
we run the algorithm of
\autoref{s.alg.pencils}, taking~\eqref{eq.quad.prop} for the defining
property~$\prop$.
The threshold is
%are $\ls|\graph|>\Mg:=16$ and
$\ls|\Fn S|>\Mm$, where $\Mm:=\bm(D)$ is the bound stated in the theorem; then, in
view of \autoref{lem.quad.sections}, we can assume that
$\ls|\graph|>\Mg:=16$. (For $D=6$ or~$8$ we can
%further
increase~$\Mg$
to~$20$ or~$18$, respectively.)
Furthermore, assuming~\eqref{eq.quad.order} and
taking~$l_1$ for~$l$, we have
%a lower bound
\[*
\ls|\sec l|=\ls|\sec_1l|\ge v_{\min}:=
\bigl\lceil\tfrac14(\Mm+1-\ls|\graph|)\bigr\rceil.
\]
Configurations of low rank (Step~$3$) are analyzed as described in
\autoref{ss.quad}.
%Each configuration~$S_0$ to be considered at Step~$3$ of the algorithm (see
%\autoref{ss.step.3}) has $\rank S_0\ge18$, and we process these
%configurations as explained in \autoref{rem.quad.extra}.
The result
%of the computation
is the conclusion that $\ls|\Fn S|\le\bm(D)$ unless
$D=4$ and
$\Fn S\cong\QC_{36}'$, $\QC_{36}''$, $\QC_{34}'$, $\QC_{33}$, $\QC_{32}\0$,
as stated.
As a by-product (\cf. \autoref{s.validation}), we obtain \autoref{ad.quad}.
%The lists of values taken by $\ls|\Fn X|$ in \autoref{ad.quad} are a
%by-product of the computation (\cf. \autoref{s.validation}), and the
%assertion on triquadrics is checked directly.
%\autoref{ad.quad} is a
%by-product of the computation (\cf. \autoref{s.validation}).
\qed

%\medbreak
\subsection{Spatial quartis}\label{s.quartics}
If $X\subset\Cp3$ is a smooth quartic and $\Fn X$ is triangle free, the best
known bound is $\ls|\Fn X|\le52$, whereas the best example so far had
$33$ lines.
(\autoref{ex.quad.D=2} below has $37$ lines.)
The last proof applies to $D=2$, and we state this fact separately
(combining it with the results of the previous sections). This leaves
configurations with a biquadrangle as the only case still open.
%Recall that, in this section, we ignored the case $D=2$ (spatial quartics),
%even though the question on the
%maximal number of lines in a triangle free configuration in a quartic is still
%open: the best known bound is 52, and the best example so far had
%33 lines. However, the last proof applies to quartics as well, and we
%state this fact separately, taking also into account the results of the previous
%sections.

\proposition\label{prop.quad.D=2}
Let $X\subset\Cp3$ be a smooth quartic, and assume that $\Fn X$ has
no triangles or biquadrangles. Then $\ls|\Fn X|\le36$,
and this bound is sharp.
\done
%A remark on quartics ($D=2$):\mnote{\todo: to do}
%at most 36 if biquadrangle free; otherwise, there is an example with 37.
\endproposition

\example\label{ex.quad.D=2}
There exists a smooth spatial quartic $X\subset\Cp3$ with the graph
$\Fn X$ triangle free and $\ls|\Fn X|=37$. This configuration, found
in~\cite{degt:Fano.graphs},
was obtained
during the experiments leading to the proof of \autoref{lem.biquad.D=3}.
\endexample

\section{Triangular configurations}\label{S.trig}

A graph~$\graph$ and the corresponding configuration
are \emph{triangular} if $\girth(\graph)=3$.
Large triangular configurations of lines in smooth quartics
%(the case $D=2$)
%although quite complicated, have already
have been studied in~\cite{DIS} (see $D=2$ in \autoref{tab.main}).
Therefore,
%in this section
we confine ourselves to the values $D\ge3$.

\subsection{Statements}\label{s.trig.statements}
The following theorem is the principal result of this section.

\theorem[see \autoref{proof.trig}]\label{th.trig}
Let $X\subset\Cp{D+1}$ be a smooth $2D$-polarized $K3$-surface, $D\ge3$,
and assume that
the Fano graph $\Fn X$ is triangular, \ie, $\girth(\Fn X)=3$.
Then, with few exceptions, we have the following sharp bounds\rom:
\[*
\minitab\quad
D=           & 3& 4& 5&\{7\bmod3\}&\otherwise\cr
\ls|\Fn X|\le&36&29&24&         24&         21\cr
\endminitab
\]
The exceptions are $D=3$, $\Fn X\cong\TC_{42},\TC_{38}$ and
$D=4$, $\Fn X\cong\TC_{33}$ \rom(see \autoref{tab.main}\rom).
%Then, we have the following sharp bounds on
%the cardinality $\ls|\Fn X|$\rom:
%\roster*
%\item
%if $D=3$, then $\ls|\Fn X|=42$, $\ls|\Fn X|=38$, or $\ls|\Fn X|\le36$\rom;
%\item
%if $D=4$, then $X$ is not a triquadric and
%$\Fn X\cong\TC_{33}$ or $\ls|\Fn X|\le29$\rom;
%\item
%if $D=5$ or
%%$D\ge7$ and $D=1\bmod3$
%$D\in\{7\bmod3\}$, then $\ls|\Fn X|\le24$\rom;
%\item
%otherwise, $\ls|\Fn X|\le21$.
%\endroster
%If $D=3$ and $n=38$ or~$42$, then
%there is a unique configuration $\TC_{n}$ with
%$\ls|\TC_{n}|=n$\rom; it
%is realized by a single connected $1$-parameter family
%of sextic surfaces.
\endtheorem

\addendum[see \autoref{proof.trig}]\label{ad.trig}
If $D=3$, the number of lines in triangular configurations takes
all values in $\{3\dts36,38,42\}$.
\endaddendum

The graph $\TC_{42}$ is $9$-regular; it contains two minimal pencils,
which are both of type \smash{$3\tA_2+3\bA_1$}. The corresponding equilinear
family contains models of the two discriminant minimizing surfaces $X([2,1,20])$ and
$X([6,3,8])$ (see~\cite{degt:singular.K3}).

%As in the previous sections, we describe the combinatorial structure of large
%triangular graphs by listing the types of all pencils
%with a fiber of type~$\tA_2$. The
%extremal graphs can be described as follows:
%\roster*
%\item
%$\TC_{42}$:\enspace
%$\trigs{[ [ 6, 3 ], 2 ]}$; the graph is $9$-regular;
%it is represented by the model $6_{42}$ of
%the discriminant minimizing singular $K3$-surfaces
%$X([2,1,20])$ and $X([6,3,8])$ and the model $6_{42}'$ of the surface $X([6,0,8])$
%(see~\cite{degt:singular.K3});
%\item
%$\TC_{38}$:\enspace
%$\trigs{[ [ 5, 4 ], 2 ]}$; the graph is represented by the
%model $6_{38}$ of $X([6,0,8])$;
%\item
%$\TC_{33}$ is a single pencil of type $8\tA_2$ with $9$ simple sections;
%\item
%if $D\ge5$, the extremal configurations (unique if $D\ne6$) are
%$\trigs{[ [ 6, 3 ], 1 ]}$, $\trigs{[ [ 8, 0 ], 1 ]}$, or
%$\trigs{[ [ 7, 0 ], 1 ]}$ for $D=0$, $1$, or $2\bmod3$, respectively
%(see \autoref{tab.parabolic});
%\item
%if $D=6$, there is another extremal configuration described as
%$\trigs{[ [ 6, 0 ], 1 ]}$.
%\endroster
%More details are found in \autoref{tab.main} (see
%also~\cite{degt:Fano.graphs}).

\subsection{Lines intersecting a triangle}
Consider a geometric configuration~$S$ and a
triangle
$\fiber=\{l_1,l_2,l_3\}\subset\Fn S$. Let $\graph:=\pencil(\fiber)$.
The next lemma is an application of \autoref{prop.hyperbolic}
to simple test configurations
(\cf. also \autoref{prop.valency.trig}).
%consisting of~$\fiber$ and up to three other lines.

\lemma\label{lem.trig}
Assume that $D\ge3$. Then any section of~$\graph$ is simple.
%, so that
%\[*
%\sec\graph=\sec_1\graph=\sec l_1\cup\sec l_2\cup\sec l_3
%\]
%is a disjoint union.
Furthermore,
\roster
\item\label{trig.disjoint}
each subgraph $\sec l_i$, $1\le i\le3$, is discrete\rom;
\item\label{trig.<=1}
$(\sec l_i)\*(\sec l_j)\le1$ for $1\le i\ne j\le3$\rom;
\item\label{trig.D>=8}
if $D\ge4$, then $\sec\graph$ is discrete.
\done
\endroster
\endlemma

%Note also that, if $D=4$, the surface $X\subset\Cp5$ cannot be a triquadric,
%as the class $e:=\kk_\fiber=l_1+l_2+l_3$ is a cubic pencil as in
%\autoref{def.geometric}\iref{bad.cubic}.

\lemma\label{lem.trig.sections}
We have the following sharp bounds on the cardinality $\ls|\sec\graph|$\rom:
\[*
\minitab\quad
D=                & 3& 4&5,6&7\dts11&\otherwise\cr
\ls|\sec\graph|\le&21& 9&  3&      1&         0\cr
\endminitab
\]
If $D=3$
and
$\ls|\sec\graph|\ge20$, then either $\Fn S\cong\TC_{42}$ or
$\ls|\Fn S|\le36$.
\endlemma

\proof
If $D\ge4$, then, in view of \autoref{lem.trig}\iref{trig.D>=8}, all
statements easily follow from \autoref{prop.hyperbolic}.
%applied to simple test configurations.
If $D=3$, then \autoref{lem.trig}\iref{trig.disjoint},
\iref{trig.<=1} means that $\fiber$ is suitable for enumeration
(see \autoref{s.alg.sections}),
and we use the
algorithm of \autoref{s.alg.sections}.
\endproof

\subsection{Proof of \autoref{th.trig} and \autoref{ad.trig}}\label{proof.trig}
Using \autoref{prop.hyperbolic} and bounds~\eqref{eq.mu}, \eqref{eq.Euler+},
we can compile a list of large pencils $\graph:=\pencil(\fiber)$ admitting a
section; the result is shown in \autoref{tab.trig}.
\table
\caption{Pencils of degree~$3$ with a section}\label{tab.trig}
\hbox to\hsize{\hss\vbox{\halign{\strut\quad\hss$#$\hss\quad&&
 \quad$#$\hss\quad\cr
\noalign{\hrule\vspace{2pt}}
\ls|\graph|&\graph&D,\ \text{remarks}\cr
\noalign{\vspace{1pt}\hrule\vspace{2pt}}
24&8\tA_2       &4,\ s\in\Fano_h(\graph)\otimes\Q\cr
22&7\tA_2+\bA_1 &4,\ s\in\Fano_h(\graph)\otimes\Q\cr
21&7\tA_2       &5,\ s\in\Fano_h(\graph)\otimes\Q\cr
  &6\tA_2+3\bA_1&3,4\cr
20&6\tA_2+2\bA_1&2,3,4\cr
19&6\tA_2+\bA_1 &2,3,4,5\cr
  &5\tA_2+4\bA_1&3,4\cr
\noalign{\vspace{1pt}\hrule}
\crcr}}\hss}
\endtable
(In the first three lines, we indicate that any section lies in
$\Fano_h(\graph)\otimes\Q$ and, hence, $\Fano_h(\graph)$ has no geometric
extensions of positive corank.) Comparing this with \autoref{cor.parabolic}
and \autoref{lem.trig.sections}, we conclude that, for $D\ge6$,
the maximal cardinality $\ls|\Fn X|$ is realized by a parabolic
configuration and is as stated in the theorem.
Similarly, if $D=5$, the maximum $\ls|\Fn X|=24$
is given by the only geometric
finite index extension of $7\tA_2$, which has three sections.

The case $D=4$ is considered in \autoref{proof.special}, resulting in
\autoref{th.special}.
Thus, we
assume that $D=3$ and $\ls|\Fn X|>\Mm:=36$. In view of
\autoref{lem.trig.sections}, we can also assume that
$\ls|\sec\graph|\le \Ms:=19$
and, hence, $\ls|\graph|>\Mg:=17$. Indexing the lines
so that $\ls|\sec l_1|\ge\ls|\sec l_2|\ge\ls|\sec l_3|$, we apply the
algorithm of \autoref{s.alg.pencils}, taking for~$l$ the
line~$l_1$ and
estimating the number of sections \via
\[*
\ls|\sec l|\ge v_{\min}:=\bigl\lceil\tfrac13(37-\ls|\graph|)\bigr\rceil.
\]
The final step of the algorithm (see \autoref{ss.step.3}) is explained in
\autoref{ss.trig}.
\qed

\subsection{Proof of \autoref{th.main} and \autoref{ad.main}}\label{proof.main}
The case of spatial quartics is treated in \cite{DIS}. For the other
degrees, \autoref{th.main} is proved by comparing the bounds given
by
Theorems~\ref{th.a5}--\ref{th.mu>5} (locally elliptic graphs),
%see also \autoref{cor.le.mc}),
\ref{th.d5} (astral graphs), and
\ref{th.a4}, \ref{th.quad}, and~\ref{th.trig} (pentagonal, quadrangular, and
triangular graphs, respectively).

\autoref{ad.main} is given by~\cite{DIS} (for $D=2$) or \autoref{ad.quad}
(for $D=3,4$).
\qed

\subsection{Proof of \autoref{th.main.real}}\label{proof.main.real}
We use \autoref{th.real} and detect the
extremal graphs in \autoref{tab.main} realized by configurations of real
lines in real surfaces
(marked with a~$\mreal$ in the table). If such a graph does not exist, we
merely state that $\mr(D)=\bm(D)$ is the auxiliary bound in
\autoref{th.main}. The sharpness of this bound
for $D=8$ and $D>10$ is given by \autoref{ad.le}.

If $D=5$, the Fano graph~$\QF_{28}$ of
the discriminant minimizing surface $X([2,0,16])$ (see
\autoref{th.min.det} and~\cite{degt:Fano.graphs}) is real.
If $D=7$, a certain real
configuration $\QF_{26}$ (see~\cite{degt:Fano.graphs})
is found in the course of the proof of
\autoref{th.a4}.
\qed

\appendix

\section{Algorithms}\label{S.algorithms}

This appendix outlines a few technical details concerning the algorithms
used to enumerate large pencils or sets of sections. All computations were
done in \GAP~\cite{GAP4}.

\subsection{Defining properties}\label{s.validation}
In the enumeration routines, we are interested in the
configurations satisfying a certain
hereditary defining property~$\prop$ (such as
locally elliptic,
triangle free, quadrangle free, \etc.)
%which is inherited by subconfigurations.
We incorporate this property as part of \autoref{alg.lines}:
upon establishing the admissibility of a configuration~$S$, we
compute the graph $\Fn S$ and disregard~$S$ as
invalid if $\prop(S)$ does not hold.

Since \autoref{alg.Nikulin} relies
upon \autoref{alg.lines}, we will essentially speak about
\emph{$\prop$-\hbox{\rom(sub-\rom)}\penalty0geometric}
configurations, \ie, those admitting a geometric finite index extension
satisfying the defining property~$\prop$.

For statistical purposes (\eg, for establishing the sharpness of the bounds),
in \autoref{alg.Nikulin}
%on the number of lines),
we keep track of the sizes of all geometric configurations:
%This task is also built into \autoref{alg.Nikulin}:
as soon as a
$\prop$-geometric configuration~$S$ is
discovered, its cardinality
$\ls|\Fn S|$ is recorded.

\subsection{Pencils}\label{s.alg.pencils}
As an essential part of most proofs,
we enumerate all $\prop$-geometric
configurations~$S$
spanned by a pencil $\graph:=\pencil(\fiber)$
with a distinguished fiber~$\fiber$ and a number of sections
$s_i\in\sec_1\graph$ such that
\[
\ls|\graph|>\Mg,\qquad\ls|\Fn S|>\Mm,
\label{eq.thresholds}
\]
where $\Mg$ and~$\Mm$ are certain thresholds fixed in advance.
The types of the fibers of~$\graph$ are given by Corollaries~\ref{cor.degree}
and~\ref{cor.elliptic}, and the possible combinatorial types of
pencils are listed using \eqref{eq.mu}--\eqref{eq.Euler+} (see also
\autoref{rem.pencils}).
%all
%pencils~$\graph$ satisfying $\ls|\graph|>\Mg$ in~\eqref{eq.thresholds}
%are listed
%using~\eqref{eq.mu} and~\eqref{eq.Euler}.

Thus, we fix a combinatorial type $\graph\supset\fiber$ and change the
defining property~$\prop$ to
\[*
\prop_\graph\:S\longmapsto
 \text{$\prop(S)$ and $\graph\subset\Fn S$ is a maximal parabolic subgraph}.
\]
We also fix a line $l\in\fiber$ and,
assuming that $\ls|\sec_1l|\ge\ls|\sec_1l'|$, $l'\in\fiber$, in the
resulting
configuration~$S$, determine a lower bound $v_{\min}\le\ls|\sec_1 l|$ on the
number of simple sections through~$l$
that is necessary
for the inequality $\ls|\Fn S|>\Mm$ in~\eqref{eq.thresholds}.

\subsubsection{Step~$1$}\label{ss.step.1}
We start with the lattice $S_0$ and stabilizer~$G_0$:
\[*
S_0:=\Fano_h(\graph),\qquad
G_0:=\stab\{l\}\subset\Aut\graph.
\]
Then,  sections
$s_1,s_2,\ldots$ intersecting~$l$ are added to~$\graph$ one by one,
and we consider
consecutive overlattices and subgroups
\[*
S_i:=(S_{i-1}+\Z s_i)/\ker\supset S_{i-1},\qquad
G_i:=\stab\{s_1,\ldots,s_i\}\subset G_0.
\]
At each step, we require that
\[
\text{$S_i$ should be $\prop_\graph$-subgeometric (\cf. \autoref{s.validation})
 and $\rank S_i>\rank S_{i-1}$}.
\label{eq.step}
\]
This process continues until no other section satisfying~\eqref{eq.step}
can be added.

%\remark\label{rem.section=set}
In all cases,
%all sections $s_i\in\sec l$ are pairwise disjoint
the set $\sec l$ is discrete
(\cf. \autoref{s.alg.sections}\iref{suitable.same} below) and,
hence, a section~$s$ is described by its coordinates $\bs\subset\graph$
(see \autoref{conv.coord}), so that $l\in\bs$.
At each step, we try for~$s_i$ a single
representative of each $G_{i-1}$-orbit. We also check partially
the defining property~$\prop$, leaving the
ultimate validation to \autoref{alg.Nikulin},
where the group~$G_i$ is used to reduce the number
of candidates
to be analyzed.
%when checking whether $S_i$ is subgeometric.
%\endremark

\remark\label{rem.pool}
As a technical tool reducing the computation and overcounting, we fix an
integer $2\le p\le4$ (depending on the size $\ls|G_0|$) and store, during the
$i$-th step,
the $G_0$-orbits of all
$i$-tuples $\{s_1,\ldots,s_i\}$, $i\le p$, of sections
satisfying~\eqref{eq.step}. These pools~$P_i$ are used two-fold:
first, during
the $i$-th step (if $i\le p$), we disregard the orbits that have
already been stored; second, when
%preparing the $(i+1)$-st step and
collecting the
candidates~$s_{i+1}$ to extend a set $\{s_1,\ldots,s_i\}$,
we
select only those $s_{i+1}\in P_1$
for which each subset $s\subset\{s_i,\ldots,s_{i+1}\}$
of cardinality $n\le\min(i,p)$ is in the pool~$P_n$.
\endremark

\subsubsection{Step~$2$}\label{ss.step.2}
We use \autoref{alg.Nikulin} to compute all
$\prop_\graph$-geometric finite index extensions
of all configurations~$S$ obtained \emph{at all stages} of Step~$1$.
Most configurations violating the hypothesis
$\ls|\Fn S|\le \Mm$ are found at this point.
Then, disregarding the configurations of rank~$20$ (as admitting no further
extensions) and those in which $\ls|\sec l|<v_{\min}$ (as not meeting the
goal),
we replace each remaining geometric extension with the respective lattice
$\Fano_h(\graph\cup\sec l)$. The sets
$\sec l$ are sorted again, retaining a
single representative of each $G_0$-orbit.

\subsubsection{Step~$3$}\label{ss.step.3}
The remaining configurations have rank~$19$ or $18$. We pick
another line $l'\in\fiber$ and try to add up to two independent extra
sections $e_i\in\sec l'$. (The choice of~$l'$
%or even pair of lines $(l',l'')$
is explained and justified in \autoref{s.trig.quad} below.)
Thus, we start with
\[*
S_0:=\Fano_h(\graph\cup\{s_1,\ldots,s_n\}),\qquad
G_0:=\stab\{s_1,\ldots,s_n\}\subset\stab(l,l')\times\Bbb S(\sec l)
%G_0:=\stab(l,l',s_1,\ldots,s_n)\subset\Aut\graph
\]
and add one or two extra sections as in Step~$1$, obtaining consecutive
lattices~$S_i$ and groups~$G_i$. In addition to~\eqref{eq.step}, we require that
%ranks should increase, $\rank S_i>\rank S_{i-1}$, but
the set $\sec_1l$ should remain
constant (as otherwise $S_i$ is among the configurations considered at
Steps~1 and~2). Finally, an analogue of Step~$2$ is applied to all new lattices
obtained.

\remark\label{rem.extra=set}
The extra sections~$e_i$ may intersect the original sections~$s_k$; hence, an
extra section~$e$ is determined by its coordinates
$\be\subset\graph$ (\cf. \autoref{conv.coord}) and
a subset of $\{s_1,\ldots,s_n\}$.
%(This fact is the reason why we require that $G_0$ stabilize $\sec l$ as a
%tuple rather than as a set.)
The cardinality of the latter
%subset of $\{s_1,\ldots,s_n\}$
is usually
bounded by the defining property~$\prop$
(\cf. \autoref{s.alg.sections}\iref{suitable.distinct} below).
\endremark

\subsection{Pencils with disjoint sections}\label{s.alg.le}
If \emph{all} sections of interest are
known to be
pairwise disjoint (\eg, in
the treatment of the locally elliptic configurations, see
\autoref{proof.le}), the algorithm of \autoref{s.alg.pencils}
can be modified. Namely, we start with a pencil~$\graph$ and a distinguished
fiber $\fiber=\{l_1,\ldots,l_n\}$ and fix the \emph{goals}
$v_i\le\ls|\sec l_i|$, $i=1,\ldots,n$. Then, disjoint sections are added
one by one as in \autoref{ss.step.1}, except that we do \emph{not} require
that the rank of the configuration should increase at each step.
As soon as a configuration~$S$ of the maximal rank~$20$ is obtained, it is
analyzed as in \autoref{ss.step.2}.
Configurations of lower rank \emph{that do not meet the goals} are discarded
as small.
%are treated by other means, \eg, by changing the
%goals or observing that no more sections can be added.

\subsection{Triangular and quadrangular configurations}\label{s.trig.quad}
\latin{A posteriori}, it turns out that it is only triangular
(for $D=3$) and
quadrangular configurations that need the algorithm of
\autoref{s.alg.pencils} to full extent. Below, we explain and justify the
choice of the additional line(s) $l',l''\in\fiber$ for Step~$3$ (see
\autoref{ss.step.3}).

\subsubsection{Triangular configurations\pdfstr{}{
 \rm(see \autoref{proof.trig})}}\label{ss.trig}
Let $\fiber:=(l_1,l_2,l_3)$ be the triangle used in the proof, so that
$l=l_1$.
Left to Step~$3$ are four configurations~$S_0$ of rank~$19$ and three
configurations of rank~$18$. In the former case, a single extra section is to be
added and, by the obvious symmetry, we can use $l'=l_2$.

If $\rank S_0=18$, then,
in all three cases, $\sec l_2=\sec l_3=\varnothing$ in each
geometric finite
index extension $S_0'\supset S_0$.
In two cases,
$\sec l_2\ne\varnothing$ and $\sec l_3\ne\varnothing$
in each geometric corank~$1$ extension
$S_1'\supset S_0$; hence, any geometric corank~$2$ extension is
generated over~$S_0$ and~$\Q$
by two sections adjacent to the same line, which, by symmetry, can be chosen
to be~$l_2$.
In the exceptional case, we have $\graph\cong5\tA_2+3\bA_1$ and
\[*
\ls|\sec l_1|=7,\quad \ls|\sec l_2|=1,\quad \ls|\sec l_3|=0;
\]
%and $\graph$ is of type $5\tA_2+3\bA_1$;
%by the above argument, in any
hence, if a geometric corank~$2$ extension
$S_2'\supset S_0$ is not generated by two sections in the same set, we have
$\ls|\sec\graph|\le9$ and $\ls|\Fn S|\le27$.
\qed

\subsubsection{Quadrangular configurations\pdfstr{}{
 \rm(see \autoref{s.biquad} and \autoref{proof.quad})}}\label{ss.quad}
%Consider the quadrangle $\fiber:=(l_1,l_2,l_3,l_4)$
Let $\fiber:=(l_1,l_2,l_3,l_4)$ be the quadrangle
used in the proofs, with the numbering
satisfying~\eqref{eq.quad.order} except that, for \autoref{lem.biquad.D=3},
we lift the restriction $\val l_1\ge\val l_2$.
Each configuration~$S_0$ to be considered at Step~$3$
%of the algorithm (see \autoref{ss.step.3})
has $\rank S_0\ge18$, but the number of configurations
is quite large.
In order to avoid
tedious case-by-case analysis,
we
merely
add up to two extra
sections $e'\in\sec l'$, $e''\in\sec l''$
trying \emph{all} pairs $(l',l'')$;
since $\val l_2\ge\val l_4$ by~\eqref{eq.quad.order},
it suffices to consider
%all pairs
\[*
(l',l'')=(l_2,l_2),\quad (l_2,l_3),\quad (l_2,l_4),\quad (l_3,l_3).
\]
(If $\rank S_0=19$, \ie, only one extra section is to be added,
we try $l'=l_2$ or~$l_3$.)
The maximal number of original sections $s_i\in\sec l_1$ that an extra section
$s^*\in\sec l^*$ may intersect (see \autoref{rem.extra=set}) is given by
\autoref{lem.quad}\iref{quad.D=3}--\iref{quad.D>7}, \iref{quad.D=2,3},
depending on the degree $h^2=2D$, line $l^*=l_2$ or~$l_3$, and whether the
resulting
configuration is or is not allowed to have biquadrangles.
If $l\ne l''$, we also need to take into account the intersection
$e'\cdot e''=0$ or~$1$. These complications slow the computation down, but it
still remains reasonably feasible.
\qed

\subsection{Sets of sections}\label{s.alg.sections}
Another part of the proofs is the classification of sets of sections in the
presence of a fixed parabolic subgraph (affine Dynkin diagram)~$\fiber$ and
a hereditary defining property~$\prop$.

Fix an ordering $\fiber=\{l_1,\ldots,l_n\}$ of the lines
constituting~$\fiber$. We say that $\fiber$ is \emph{suitable for enumeration}
if its sections in any $\prop$-geometric configuration
have the following properties:
\roster
\item\label{suitable.same}
%any two sections $s',s''\in\sec l_i$ are disjoint, $1\le i\le n$, and
each subgraph $\sec l_i$, $1\le i\le n$, is discrete, and
\item\label{suitable.distinct}
%any section $s_i\in\sec l_i$ intersects at most one section $s_j\in\sec l_j$,
%$i\ne j$.
$(\sec l_i)\*(\sec l_j)\le1$ for $1\le i\ne j\le n$.
\endroster

\subsubsection{Steps $1_1,1_2$}\label{ss.step.12}
Assuming these properties, we fix the \emph{goals} $v_i=\ls|\sec l_i|$ and
start with the lattice
\[*
S_2:=S_{2,0}:=\Fano_h(\fiber\cup\sec l_1\cup\sec l_2),
\]
where $\sec l_i=\{s_{i1},\ldots,s_{iv_i}\}$ consists of $v_i$ pairwise
disjoint sections, $i=1,2$, and $s_{1k}\cdot s_{2k}=1$ for
$1\le k\le r$ (another parameter fixed in advance) and
$s_{1p}\cdot s_{2q}=0$ for all other
pairs $(p,q)$. Let $G_2=\SG{r}\times\SG{v_1-r}\times\SG{v_2-r}$ be the group
of symmetries of this configuration.

\subsubsection{Step $1_m$, $3\le m\le n$}\label{ss.step.1m}
At this step, the previously constructed configuration
$S_{m-1}=S_{m-1,0}$ is
considered frozen, \ie, we change the defining property~$\prop$ to
\[*
\prop_{m-1}\:S\longmapsto
 \text{$\prop(S)$ and $\ls|\sec l_k|=v_k$ for $k<m$}.
\]
We add to~$S_{m-1}$, one by one, \emph{exactly} $v_m$ pairwise disjoint
sections $s_{mi}\in\sec l_m$ and consider consecutive overlattices
$S_{m-1,i}:=(S_{m-1,i-1}+\Z s_{m,i})/\ker$.
Unlike~\eqref{eq.step}, we do \emph{not} require that the ranks should
increase. Each section $s_{mi}$ is determined by an $(m-1)$-tuple of at most
one-element sets $[s_{mi}]_k\subset\sec l_k$, $k<m$.
Since the symmetry groups involved are typically much smaller than
those in
\autoref{s.alg.pencils}, we use a more aggressive algorithm: at each step, we
collect the candidates $s_{mi}$ extending all previously constructed sets
$\{s_{m1},\ldots,s_{m,i-1}\}$, compute $G_{m-1}$-orbits of \emph{all}
$i$-tuples $\{s_{m1},\ldots,s_{mi}\}$ obtained, and choose one
representative of each orbit. We also use pools $P_i$, $i\le4$, as
in \autoref{rem.pool}.

Upon completion of this step, we obtain a collection of lattices
$S_m:=S_{m-1,v_m}$ and stabilizers~$G_m$,
\[*
S_m:=(S_{m-1}+\Z\sec l_m)/\ker,\qquad
G_m:=\stab\sec l_m\subset G_{m-1},
\]
to be used at Step~$1_{m+1}$. (We retain only those configurations~$S_m$
which admit a $\prop$-geometric extension with \emph{exactly} $v_k$ sections
intersecting~$l_k$, $k\le m$.)

If only sets of sections are to be classified (\eg, in order to determine
their maximal number), we stop at this point.
Otherwise, if we are interested in all $\prop_n$-geometric extensions~$S$ of
the lattices~$S_n$ obtained satisfying an inequality
\[
\ls|\Fn S|>\Mm
\label{eq.threshold}
\]
fixed in advance,
we proceed similar to \autoref{s.alg.pencils}.

\subsubsection{Step~$2$}\label{ss.step.2.sets}
For each lattice~$S_n$ obtained at the \emph{final} step~$1_n$,
compute all
$\prop_n$-geometric finite index extensions
(by \autoref{alg.Nikulin}); then, select those
satisfying \eqref{eq.threshold}. After this step,
we disregard all lattices of rank~$20$.

\subsubsection{Step~$3$}\label{ss.step.3.sets}
For the remaining configurations $S_n:=S_{n,0}$,
we increase the
rank by adding, one by one, $(20-\rank S_n)$ disjoint
\emph{extra lines} (\ie, lines~$e_i$ disjoint from
all $l_1,\ldots,l_n$) and considering the configurations
\[*
S_{n,i}:=(S_{n,i-1}+\Z e_i)/\ker
\]
and all their $\prop_n$-geometric finite index extensions
computed by \autoref{alg.Nikulin}.
At each step, we require that $\rank S_{n,i}>\rank S_{n,i-1}$.
An extra line~$e$ can be regarded as a subset $[e]\subset\sec\fiber$
(\cf. \autoref{conv.coord};
the cardinality of this subset and/or its intersections with $\sec l_i$
is usually bounded by geometric arguments),
and we use an aggressive algorithm similar to \autoref{ss.step.1m}.

In the rest of this section, we justify the fact that, if
$\rank S_n=18$, it suffices to add pairs of \emph{disjoint} extra
lines~$e_1$, $e_2$, keeping the computation feasible. The other cases,
where $\rank S_n\le17$, are considered below when they appear.

We say that a $\prop_n$-geometric corank~$r$ extension $S'\supset S_n$ is
\emph{well generated} if it has a \emph{good basis}, \ie, a collection of
\emph{pairwise disjoint} lines $e_1,\ldots,e_r\in\pencil_{S'}(\fiber)\sminus\fiber$
generating~$S'$ over $S_n$ and~$\Q$.
Any line appearing in a good basis is called \emph{generating}.
Given a good basis $\{e_1,\ldots,e_r\}$,
we denote by $E_i:=E(e_i)$ the connected component
of the pencil $\pencil_{S'}(\fiber)$ containing the vertex~$e_i$, $i=1,\ldots,r$.
If $r=1$, we consider also the union
$\Eg(S')=\bigcup_e E(e)$ over all generating lines~$e$.

Let $\Cal T(\fiber)$ be the set of all
types of fibers that can occur in the pencil $\pencil(\fiber)$ (in any
$\prop_n$-geometric configuration), see Corollaries~\ref{cor.degree}
and~\ref{cor.elliptic}.
%Then, we have the following obvious statement.

\lemma\label{lem.disjoint}
Given $S_n$ as above and a
$\prop_n$-geometric corank~$1$ extension $S'\supset S_n$,
assume that, for any inclusion $\Eg(S')\into\fiber'$
to a graph
$\fiber'\in\Cal T(\fiber)$, we have
\[*
\ls|\fiber'|-\ls|\Eg(S')|\le \Mm-\ls|\Fn S'|.
\]
Then, any $\prop_n$-geometric corank~$2$ extension
$S''\supset S'\supset S_n$
satisfying~\eqref{eq.threshold} is
well generated over~$S_n$.
%\done
\endlemma

\proof
An extension $S''\supset S'$ is well generated over~$S_n$
if and only if there is a line
$e_2\in\pencil_{S''}(\fiber)\sminus\pencil_{S'}(\fiber)$ disjoint from at
least one generating line~$e_1$ of~$S'$. If this is \emph{not} so, all new
lines~$e_2$ and all generating lines~$e_1$ must be in the same
fiber~$\fiber'$ of $\pencil_{S''}(\fiber)$, and we have
$\ls|\Fn S''|\le\ls|\Fn S'|+\ls|\fiber'|-\ls|\Eg(S')|\le\bm$.
\endproof

\corollary\label{cor.disjoint}
Under the assumptions of \autoref{lem.disjoint}, if $\Eg(S')$
does not admit an embedding into $\fiber'\in\Cal T(\fiber)$,
then any $\prop_n$-geometric extension
$S''\supset S'\supset S_n$ of corank~$2$ is
well generated over~$S_n$.
\done
\endcorollary

\subsection{Quadrangular configurations}\label{s.alg.quad}
In the proof of \autoref{lem.quad.sections}, for $D=4$ we need to use
Step~$3$
of the previous algorithm to list all
geometric configurations~$S$
containing a pencil $\graph:=\pencil(\fiber)$, $\fiber\cong\tA_3$, and such
that $\ls|\Fn S|>30$ and $\ls|\sec\graph|\ge15$.
In this case, the rank $\rank S_4$ can be as low as~$15$, so that we may need
to add up to five extra lines, and \autoref{lem.disjoint} does not always
apply.

Fix a configuration~$S_4$ obtained at Step~$1_4$.
Recall that we consider triangle free
configurations and insist that the set of sections of the pencil does not
increase; in other words, we take for the defining property
(depending on~$S_4$ or, more precisely, on the fixed goals
$v_i=\ls|\sec l_i|$ in~$S_4$, $1\le i\le4$)
\[*
\prop_4\:S\longmapsto
 \text{$S$ is triangle free and $\sec l_i=v_i$ in~$S$, $1\le i\le4$}.
\]
In the notation
introduced prior to \autoref{lem.disjoint},
%of the previous section,
we have $E_i\cong\tA_3$, $\bA_2$, or~$\bA_1$
for any good basis $\{e_1,\ldots.e_n\}$, and the following
two statements are proved by a simple analysis of the modifications of the
components under $\prop_4$-geometric extensions. (The hypotheses of both
statements are to be checked for \emph{geometric} extensions, \ie, those
obtained after applying \autoref{alg.Nikulin}.)

\lemma\label{lem.quad.large}
Assume that a $\prop_4$-geometric
corank~$r$
extension $S'\supset S_4$ has a
good basis $e_1,\ldots,e_r$, and let
$\br:=\#\bigl\{i\le r\bigm|E_i\cong\bA_1\}$.
Then any $\prop_4$-geometric corank~$1$ extension $S''\supset S'$
such that $\ls|\Fn S''|>\ls|\Fn S'|+2\br$ is well generated over~$S_4$.
\endlemma

\proof
In the worst case scenario, a component $E_i\cong\bA_1$ extends to
$E_i''\cong\tA_3$, so that the line $e_i''\in E''$ opposite to~$e_i$ is already
in $\pencil_{S'}(\fiber)$, and all new components are of this form. In this case,
$\ls|\Fn S''|\le\ls|\Fn S'|+2\br$.
\endproof

\corollary\label{cor.quad.even}
Assume that a $\prop_4$-geometric corank~$r$ extension $S'\supset S_4$ has a
good basis such that $\ls|E_i|\ge2$, $1\le i\le r$.
Then any $\prop_4$-geometric corank~$1$ extension $S''\supset S'$ is well
generated over~$S_4$.
\done
\endcorollary

\proof[End of the proof of \autoref{lem.quad.sections}]
%Now, we can complete the proof of \autoref{lem.quad.sections},
We start from a configuration~$S_4$, $\rank S_4\ge 15$, and add extra
lines one by one, verifying, at each step, that all subsequent corank~$1$
extensions are still well generated.

If $\rank S_4=19$, there is nothing to prove: all extensions are well
generated.

If $\rank S_4=18$, then
every $\prop_4$-geometric corank~$1$ extension $S'\supset S_4$ either
satisfies the hypotheses of
%\autoref{lem.disjoint} (or rather \autoref{cor.disjoint})
Corollary~\ref{cor.disjoint}
or \ref{cor.quad.even} or has
$\ls|\Fn S'|\le28$,
and then we can use \autoref{lem.quad.large};
thus, any $\prop_4$-geometric
corank~$2$ extension $S''\supset S_4$ that satisfies the inequality
$\ls|\Fn S''|>30$ is well generated.

If $\rank S_4=15$ or~$16$, then, step-by-step, one shows that,
for \emph{each} $\prop_4$-geometric extension
$S\rr\supset S_4$, one has $\ls|E|\ge2$ for each connected component
$E\subset\pencil_{S\rr}(\fiber)$;
then, by \autoref{cor.quad.even}, all extensions at the next step
are also well generated.

Finally, let $\rank S_4=17$ and pick a
$\prop_4$-geometric corank~$1$ extension
$S'\supset S_4$. The computation shows that, if $S'$ does \emph{not} satisfy
the hypotheses of \autoref{lem.disjoint} or \autoref{cor.quad.even}, then
$\Eg(S')\cong\bA_1$ or $2\bA_1$ and $\ls|\Fn S'|\le26$.
%Mimicking the proof of \autoref{lem.quad.large}, one concludes that
By \autoref{lem.quad.large},
any
extension $S\rr\supset S'$ such that $\ls|\Fn S\rr|>30$
contains a well-generated $\prop_4$-geometric
corank~$2$ extension $S''\supset S'\supset S$
(\ie, there is a line
$l\in\pencil_{S\rr}(\fiber)\sminus\pencil_{S'}(\fiber)$
disjoint from at least one good generator of~$S'$).

Thus,
let $S''\supset S_4$ be a well-generated
$\prop_4$-geometric corank~$2$ extension.
%consider a $\prop_4$-geometric corank~$2$ extension $S''\supset S_4$.
In all but a few cases, \autoref{lem.quad.large}
implies that any further extension $S'''\supset S''$ is well generated whenever
$\ls|\Fn S'''|>30$.
%rules
%out the existence of a further extension $S'''\supset S''$ that is not well
%generated and has more than $30$ lines.\mnote{\todo: more details}
In each exceptional case, the
%subgraph~$\gen$
graph $\gen\subset\graph_{S''}(\fiber)$
of generating lines splits,
$\gen=\gen_1\cup\gen_2$, $\gen_1\cap\gen_2=\varnothing$, so that
\[*
\gen_1=\{e_1',e_1''\}\cong2\bA_1,\quad
\gen_2\cong2\bA_1,\ \tA_3+\bA_2,\ \mbox{or}\ 2\tA_3
\]
and good are the bases of the form $\{e_1,e_2\}$, where $e_1\in\gen_1$,
$e_2\in\gen_2$.
(If $\gen_2\cong2\bA_1$, then we also have $\ls|\Fn S''|=28$.)
A simple argument shows that
any offending extension
$S'''\supset S''$ is generated by an extra line~$e_2$
%intersecting~$e_1'$ and~$e_1''$,
such that $e_2\cdot e_1'=e_2\cdot e_1''=1$,
so that $e_1',e_2,e_1''$ are part of a quadrangle. This
observation limits the number of choices for~$e_2$;
we do list all such extensions
(with the intersections of the extra lines prescribed)
and find out that they are never
$\prop_4$-geometric.
%such lines
%can easily be listed, and one finds out that an extension~$S'''$ as above
%cannot be $\prop_4$-geometric.
\endproof
%\qed

\subsection{Pentagonal configurations}\label{s.alg.pent}
In the proof of \autoref{lem.pent},
where the fiber $\fiber\cong\tA_4$ is a pentagon,
we encounter configurations $S_5$ of rank $\rank S_5\ge17$.
%The case where
%$\rank S_5=18$ is completely covered by the following
The following statement is a
refinement of
\autoref{lem.disjoint} using the known structure of~$\fiber$.
We assume that the defining property is $\prop\:S\mapsto(\girth(\Fn S)=5)$
and consider a
subgeometric configuration~$S_5$ obtained after Step~$1_5$ of the algorithm
of~\autoref{s.alg.sections}.
The other notations and terminology are introduced prior to
\autoref{lem.disjoint}.

\lemma\label{lem.pent.large}
Let $S_5$ be as above and $S'\supset S_5$ a
$\prop_5$-geometric corank~$1$ extension.
Define $\Gd(S')$ as follows\rom:
\roster*
\item
$\Gd(S')=2$ if $\Eg\cong\bA_1$\rom;
\item
$\Gd(S')=1$ if $\Eg\cong2\bA_1$ or $\Eg\cong\bA_2$ with only one line
generating\rom;
\item
$\Gd(S')=0$ in all other cases.
\endroster
Then, any $\prop_5$-geometric corank~$2$ extension
$S''\supset S'\supset S_5$
satisfying the inequality $\ls|\Fn S''|>\ls|\Fn S'|+\Gd(S')$
is well generated over~$S_5$.
\done
\endlemma

\proof[Proof of \autoref{lem.pent}\iref{pent.2}--\iref{pent.6}]
For the computation, we need to show that Step~$3$ of the algorithm applies,
\ie, that any sufficiently large $\prop_5$-geometric extension of any
configuration~$S_5$ is well generated. If $\rank S_5=18$,
we apply \autoref{lem.pent.large},
%can be
checking the hypotheses case by case.
If $D=5$, there are four configurations~$S_5$ of rank~$17$, and
\autoref{lem.pent.large} implies that any sufficiently
large extension $S''\supset S_5$ of rank~$19$ is well generated.
%One can see that a
A further extension $S'''\supset S''$ is well generated whenever
\[
\ls|\Fn S'''|>\ls|\Fn S''|+4,
\label{eq.pent.4}
\]
the worst case scenario being that of $S''$
having but two disjoint generating lines. The computation shows that
$\ls|\Fn S''|\le25$ and, if $\ls|\Fn S''|=25$,
the pencil
%$\graph\subset S''$
$\pencil_{S''}(\fiber)$
has a type~$\bA_2$ fiber with both lines generating,
%the threshold above
reducing~\eqref{eq.pent.4} to
$\ls|\Fn S'''|>27$.
%$\ls|\Fn S'''|>\ls|\Fn S''|+2=27$.
It follows that any
$\prop_5$-geometric extension $S\supset S_5$ with
$\ls|\Fn S|>28$ is well generated, and Step~$3$ of the algorithm
%(see \autoref{ss.step.3.sets})
%shows that the only such extension is
results in a unique such extension, \viz.
$\QF_{30}'''$.
\endproof

%Similar arguments can be used in the treatment of astral configurations.

\subsection{Astral configurations}\label{proof.astral.large}
Let \smash{$\fiber\cong\tD_4$} and $\graph:=\pencil(\fiber)$.
The \emph{central line} of a type
\smash{$\tD_4$} fiber $\fiber'\subset\graph$ is the only line
$l\in\fiber'$ of relative valency~$4$.
The following simple observations follow from the assumption that
$\girth(\graph)\ge6$:
\roster*
\item
all sections of~$\graph$ are simple or double (\ie, bisections);
\item
a section $s\in\sec\graph$
is double if and only if it intersects the central line of any
(equivalently, each)
type \smash{$\tD_4$} fiber;
\item
if $\graph$ has a bisection, then it has no fibers of
type~\smash{$\tA_5$};
\item
if $\graph$ has two bisections, then
it has no parabolic fibers other than~$\fiber$.
%$\fiber$ is the only parabolic fiber.
\endroster

\proof[Proof of \autoref{lem.astral.large}]
%\subsection{Proof of \autoref{lem.astral.large}}\label{proof.astral.large}
We proceed as in \autoref{ss.step.3.sets}, mainly
using \autoref{lem.disjoint}. Below, we consider extensions of a
configuration $S_5$ spanned by~$h$ and $\fiber\cup\sec\graph$.

Each of the three sets $\sec\graph$ of cardinality~$12$ is covered by
\autoref{lem.disjoint}: one has $\rank S_5=18$ and any $\prop_5$-geometric extension
$S''\supset S_5$ such that $\ls|\Fn S''|>22$ is well generated.
Using \autoref{ss.step.3.sets}, we immediately obtain
Statement~\iref{astral.large.12}.

There are several sets of cardinality~$11$, but we only consider the
three that are geometric for $D=5$. Two of these sets have more than one
bisection; hence, $\graph$ has no parabolic fibers other than~$\fiber$.
We have $\rank S_5=16$ and the computation shows that each
well-generated
geometric extension $S''\supset S'\supset S_5$ of corank $c=2$ or~$1$
has exactly $c$ extra lines.
Since all fibers of~$\graph$ other than~$\fiber$ are of types~$\bA_p$,
$p\le4$,
%it is clear that
each generating line can intersect at most two other lines $l\in\graph$.
Hence,
any further extension
$S\supset S_5$ such that $\ls|\Fn S|>\ls|\Fn S_5|+2+4=22$ must be well
generated;
%however,
using \autoref{ss.step.3.sets}, we show that
such extensions do not exist.

In the third case, there is a single bisection~$s$
%(hence, all parabolic fibers are of type~\smash{$\tD_4$})
and $\rank S_5=17$. Well-generated geometric extensions
$S''\supset S'\supset S_5$ of corank $c=2$ or~$1$ have up to $c+1$ extra
lines, which are all generating, pairwise disjoint, and \emph{disjoint
from~$s$}, implying that an extra
line~$e$ cannot be the central line of a type \smash{$\tD_4$} fiber
of~$\graph$; hence, as in the previous case, $e$ intersects at most two
other lines $l\in\graph$.
As above, we conclude that any geometric extension $S\supset S_5$ such that
$\ls|\Fn S|>21$ must be well generated; the classification of such
extensions using \autoref{ss.step.3.sets} proves
that $\ls|\Fn S|\le21$.
%Statement~\iref{astral.large.11}.
%\qed
\endproof

{
\let\.\DOTaccent
\def\cprime{$'$}
\bibliographystyle{amsplain}
\bibliography{degt}
}

\end{document}